\documentclass{article}
\title{Hopf Algebras and Congruence Subgroups}
\author{Yorck Sommerh\"auser \qquad Yongchang Zhu}
\date{}

\usepackage{amsmath}
\usepackage{theorem}
\usepackage{amssymb}
\usepackage{latexsym}
\usepackage{diagramb}
\usepackage{endnotes}

\newcounter{num}

\makeatletter
\renewcommand{\subsection}{\@startsection{subsection}{2}{0em}%
{\baselineskip}{-0em}{\bfseries\normalsize}}
\makeatother

\makeatletter
\newcommand{\listofdefinitions}{\@starttoc{def}}

\newcommand{\l@definition}[2]{\par\noindent#1 {\itshape #2}}
\makeatother

\theoremstyle{plain}
\theorembodyfont{\rmfamily}

\newtheorem{thm}{Theorem}

\newtheorem{prop}[thm]{Proposition}
\newtheorem{lemma}[thm]{Lemma}

\newtheorem{corollary}[thm]{Corollary}

\newtheorem{pf}{Proof.}

\newtheorem{defn}[thm]{Definition}

\theoremstyle{break}
\theorembodyfont{\rmfamily}
\newtheorem{propb}[thm]{Proposition}

\newtheorem{lemb}[thm]{Lemma}

\newtheorem{corb}[thm]{Corollary}

\newcommand{\jac}[2]{\genfrac{(}{)}{}{}{#1}{#2}}

\newcommand{\qed}{$\Box$}

\newcommand{\Ch}{\operatorname{Ch}}
\newcommand{\End}{\operatorname{End}}

\newcommand{\Gal}{\operatorname{Gal}}
\newcommand{\SL}{\operatorname{SL}}
\newcommand{\GL}{\operatorname{GL}}
\newcommand{\PSL}{\operatorname{PSL}}
\newcommand{\PGL}{\operatorname{PGL}}

\newcommand{\Hom}{\operatorname{Hom}}
\newcommand{\Ind}{\operatorname{Ind}}
\newcommand{\Res}{\operatorname{Res}}
\newcommand{\id}{\operatorname{id}}

\newcommand{\Tr}{\operatorname{tr}}

\newcommand{\op}{\scriptstyle \operatorname{op}}
\newcommand{\cop}{\scriptstyle \operatorname{cop}}

\def\1{{(1)}}
\def\2{{(2)}}
\def\3{{(3)}}
\def\4{{(4)}}
\def\5{{(5)}}
\def\6{{(6)}}
\def\7{{(7)}}
\def\8{{(8)}}
\def\9{{(9)}}

\def\o{\otimes}

\def\da{\Delta}

\def\dh{\Delta}
\def\dd{\Delta_D{}}
\def\df{\Delta_F}
\def\ea{\varepsilon}

\def\eh{\varepsilon}
\def\ed{\varepsilon_D{}}
\def\sa{S}

\def\sd{S_D}
\def\sh{S}

\def\lma{\lambda}
\def\rha{\rho}
\def\Lma{{\Lambda}{}}
\def\ga{\gamma}
\def\Ga{\Gamma}

\def\lmh{\lambda{}}
\def\rhh{\rho{}}
\def\Lmh{\Lambda{}}
\def\Gh{\Gamma{}}
\def\lmd{\lambda_D{}}
\def\rhd{\rho_D{}}
\def\Lmd{\Lambda_D{}}

\def\H{1}
\def\A{1}

\def\D{1}

\def\ua{u}
\def\ud{u_D}
\def\ude{u_{D(1)}}
\def\udz{u_{D(2)}}

\def\N{{\mathbb N}}
\def\C{\mathbb C}
\def\Z{{\mathbb Z}}
\def\Q{{\mathbb Q}}

\def\T{{\mathbf T}}
\def\V{{\mathbf S}}
\def\C{{\mathbf C}}

\def\t{{\mathfrak T}}
\def\v{{\mathfrak S}}
\def\r{{\mathfrak R}}

\def\gt{{\mathfrak t}}
\def\gv{{\mathfrak s}}
\def\gr{{\mathfrak r}}
\def\ga{{\mathfrak a}}
\def\gd{{\mathfrak d}}

\def\bt{{\bar{\mathfrak t}}}
\def\bv{{\bar{\mathfrak s}}}
\def\br{{\bar{\mathfrak r}}}

\begin{document}

\maketitle

\begin{abstract}
\hspace{-5mm}We prove that the kernel of the natural action of the modular group on the center of the Drinfel'd double of a semisimple Hopf algebra is a congruence subgroup. To do this, we introduce a class of generalized Frobenius-Schur indicators and endow it with an action of the modular group that is compatible with the original one.
\end{abstract}
\thispagestyle{empty}

\newpage

\tableofcontents

\newpage

\section*{Introduction} \label{Sec:Introd}
\addcontentsline{toc}{section}{Introduction}
At least since the work of J.~L.~Cardy in~1986, the importance of the role of the modular group has been emphasized in conformal field theory, and it has been extensively investigated since then.\endnote{\cite{Cardy}; \cite{Verl}; \cite{MoSeib}.} This importance stems from the fact that the characters of the primary fields, which depend on a complex parameter, are equivariant with respect to the action of the modular group on the upper half plane on the one hand and a linear representation of the modular group on the other hand, which is finite-dimensional in the case of a rational conformal field theory. It was soon noticed in the course of this development that under quite general assumptions a frequently used generator of the modular group has finite order in this representation.\endnote{\cite{Vafa}.} Since this generator and one of its conjugates together generate the modular group, this leads naturally to the conjecture that the kernel of the before-mentioned representation is a congruence subgroup. After an intense investigation, this conjecture was finally established by P.~Bantay.\endnote{\cite{BanModGal}; \cite{CostGann}; \cite{Eholz1}.}

In a different line of thought, Y.~Kashina observed, while investigating whether the antipode of a finite-dimensional Yetter-Drinfel'd Hopf algebra over a semisimple Hopf algebra has finite order, that certain generalized powers associated with the semisimple Hopf algebra tend to become trivial after a certain number of steps.\endnote{\cite{KashAnti}; \cite{KashPow}.} She established this fact in several cases and conjectured that in general this finite number after which the generalized powers become trivial, which is now called the exponent of the Hopf algebra, divides the dimension of the Hopf algebra. This conjecture is presently still open. However, P.~Etingof and S.~Gelaki, realizing the connection between these two lines of thought, were able to establish the finiteness of the exponent and showed that it divides at least the third power of the dimension.\endnote{\cite{EG4}, Thm.~4.3, p.~136.} They also explained the connection of the exponent to the order of the generator of the modular group by showing that the exponent of the Hopf algebra is equal to the order of the Drinfel'd element of the Drinfel'd double of the Hopf algebra. In this context, it should be noted that this connection between Hopf algebras and conformal field theory has been intensively investigated by many authors; we only mention here the modular Hopf algebras and modular categories of N.~Reshetikhin and V.~G.~Turaev on the one hand and the modular transformations considered by V.~Lyubashenko and his coauthors on the other hand.\endnote{\cite{ResTur}; \cite{Tur}; \cite{Kerl}; \cite{KerlLyub}; \cite{Ly1}; \cite{Ly2}; \cite{LyuMaj}.}

It is the purpose of the present work to unite these two lines of thought further by establishing an analogue of Bantay's results for semisimple Hopf algebras. We will show in Theorem~\ref{CongDrinfDoubl} that the kernel of the action of the modular group on the center of the Drinfel'd double of a semisimple Hopf algebra is a congruence subgroup of level~$N$, where~$N$ is the exponent of the Hopf algebra discussed above. The proof of this theorem becomes possible by the use of a new tool, a further generalization of the higher Frobenius-Schur indicators studied earlier by Y.~Kashina and the authors.\endnote{\cite{YYY2}.} These new indicators, which we call equivariant Frobenius-Schur indicators, are functions on the center of the Drinfel'd double and carry an action of the modular group that is equivariant with respect to the action of the modular group on the center. This equivariance in particular connects, via the action of the Verlinde matrix that arises from the other frequently used generator of the modular group, the first formula for the higher Frobenius-Schur indicators with the second resp.~third formula, whose interplay is crucial for the proof of Cauchy's theorem for Hopf algebras.\endnote{\cite{YYY2}, Cor.~2.3, p.~17; Prop.~3.2, p.~23; Cor.~6.4, p.~48; Thm.~3.4, p.~26.} 

The Drinfel'd double is an example of a factorizable Hopf algebra, and the results for the Drinfel'd double can be partially generalized to this more general class. However, in the case of a factorizable semisimple Hopf algebra, the modular group acts in general only projectively on the center of the Hopf algebra. This phenomenon also occurs in conformal field theory, and also in the general framework of modular categories, of which the representation category of a semisimple factorizable Hopf algebra is an example.\endnote{\cite{BantProjKern}; \cite{Tur}, Sec.~II.3.9,  p.~98.}
But it is still possible to talk about the kernel of the projective representation, i.e., the subgroup of the modular group that acts as the identity on the associated projective space of the center. We will also show, in Paragraph~\ref{ProjCong}, that in this more general case this so-called projective kernel is a congruence subgroup of level~$N$.

However, if the Drinfel'd element of the factorizable Hopf algebra has
the same trace as its inverse in the regular representation, then the 
projective representation just discussed is in fact an ordinary linear representation. This happens in particular in the case of a Drinfel'd double, where both of these traces are equal to the dimension of the doubled Hopf algebra. If these traces coincide, it is therefore meaningful to talk about the kernel of the linear representation, and we show in Theorem~\ref{HopfSymbCong} that this kernel is also a congruence subgroup of level~$N$.

The article is organized as follows: In Section~\ref{Sec:ModGroup}, after briefly recalling some facts about the modular group, we describe a relation that characterizes the orbits of the principal congruence subgroups and plays an important role in the proof of the orbit theorem in Paragraph~\ref{OrbitThm}. In Section~\ref{Sec:QuasitriHopf}, we recall some basic facts about quasitriangular Hopf algebras and the Drinfel'd double construction, and prove some lemmas about the Drinfel'd element and the evaluation form. In Section~\ref{Sec:FactorizeHopf}, we prove some facts about factorizable Hopf algebras that are important for the equivariance properties that we will discuss later. In Section~\ref{Sec:ActModGroup}, we construct the action of the modular group on the center of a factorizable Hopf algebra. It must be emphasized that this construction is not new; on the contrary, it is discussed in abundance in the literature we have already quoted, especially in V.~G.~Turaev's monograph on the one hand and in two closely related articles V.~Lyubashenko on the other hand.\endnote{\cite{Tur}; \cite{Ly1}; \cite{Ly2}.} What we do in this section is to translate Lyubashenko's graphical proof of the modular identities into the language of quasitriangular Hopf algebras, thereby offering a presentation of these results that is not yet available in the literature in this form.\endnote{\cite{Kerl}; \cite{LyuMaj}.}

In Section~\ref{Sec:FaktSemisim}, we specialize to the semisimple case. We can then use the centrally primitive idempotents as a basis and therefore get explicit matrices for the action of the modular group constructed in Section~\ref{Sec:ActModGroup}. In the case of a Drinfel'd double, there is a different construction for the action of the modular group based on the evaluation form and using a slightly less frequently used set of generators of the modular group. This description of the action, which is crucial for the proof of the equivariance theorem in Paragraph~\ref{EquivarThm}, is given in Section~\ref{Sec:CaseDrinfDoub}.

For two modules~$V$ and~$W$ of a semisimple Hopf algebra~$H$, the modules
$V \o W$ and~$W \o V$ are in general not isomorphic. However, as we show in Section~\ref{Sec:IndMod}, the corresponding induced modules of the Drinfel'd double~$D(H)$ are isomorphic. The constructed isomorphism is the essential element for the definition of the equivariant Frobenius-Schur indicators~$I_V((m,l),z)$ in Section~\ref{Sec:EquiFrobSchur}, which depend on an $H$-module~$V$, two integers~$m$ and~$l$, and a central element~$z$ in the Drinfel'd double~$D(H)$. We then prove the equivariance theorem
$I_V((m,l)g,z) = I_V((m,l),g.z)$ for an element~$g$ of the modular group. 
In Paragraph~\ref{OrbitThm}, we prove the orbit theorem, which asserts that the equivariant indicators only depend on the orbit of~$(m,l)$ under the principal congruence subgroup determined by the exponent. This is applied in Section~\ref{Sec:CongSubgrThm} to prove the congruence subgroup theorem, which asserts that~$g.z=z$ for all~$z$ in the center of the Drinfel'd double~$D(H)$ and all~$g$ in the principal congruence subgroup. Note that the orbit theorem is an immediate consequence of the equivariance theorem and the congruence subgroup theorem. Finally, in the case of an arbitrary factorizable Hopf algebra, we prove the projective congruence subgroup theorem, which asserts that the kernel of the projective representation is a congruence subgroup.

The Wedderburn components of the character ring of a semisimple factorizable Hopf algebra are isomorphic to subfields of the cyclotomic field determined by the exponent.\endnote{\cite{YYY2}, Prop.~6.2, p.~44.}
As in conformal field theory,\endnote{\cite{BoerGoer}, App.~B, p.~302.} we therefore get an action of the Galois group of the cyclotomic field on the character ring. As we explain in Section~\ref{Sec:GaloisGroup}, this linear action of the Galois group arises naturally as the composition of the two semilinear actions that preserve the characters resp.\ the primitive idempotents of the character ring. In Section~\ref{Sec:GalInd},
we relate these actions of the Galois group to the equivariant Frobenius-Schur indicators, which enables us to show in Theorem~\ref{GalMod} that in the case of a Drinfel'd double the action of the Galois group coincides with the action of the diagonal matrices in the reduced modular group~$\SL(2,\Z_N)$. This is again confirming the parallels with conformal field theory, where the analogous result was known in many cases.\endnote{\cite{CostGann}, \S~2.3, Thm.~2, p.~7.} However, this theorem does not hold for a general semisimple factorizable Hopf algebra, as we see in Section~\ref{Sec:GalCong}: Under the assumption that the character of the regular representation takes the same value on the Drinfel'd element and on its inverse, which happens for Drinfel'd doubles, the action of the modular group, which is in general only projective, becomes an ordinary linear representation. Generalizing the congruence subgroup theorem from Paragraph~\ref{CongDrinfDoubl}, we show in Theorem~\ref{HopfSymbCong} that the kernel of this linear representation
is again a congruence subgroup of level~$N$, so that we again get an action of the reduced modular group~$\SL(2,\Z_N)$. But this time the action of the Galois group may differ from the action of the diagonal matrices by a certain Dirichlet character, which, as it generalizes the Jacobi symbol to Hopf algebras, we call the Hopf symbol.

Throughout the whole exposition, we consider an algebraically closed base field that is denoted by~$K$. From Section~\ref{Sec:FaktSemisim} on until the end, we assume in addition that~$K$ has characteristic zero. All vector spaces considered are defined over~$K$, and all tensor products without subscripts are taken over~$K$. The dual of a vector space~$V$ is denoted by $V^*:=\Hom_K(V,K)$, and the transpose of a linear 
map~$f: V \rightarrow W$ is denoted by $f^*: W^* \rightarrow V^*$. 
Unless stated otherwise, a module is a left module. Also, we use the so-called Kronecker symbol~$\delta_{ij}$, which is equal to~$1$ if $i=j$ and zero otherwise. The set of natural numbers is the set~$\N:=\{1,2,3,\ldots\}$; in particular, $0$ is not a natural number. The symbol~$\Q_m$\index{$\Q_m$} denotes the $m$-th cyclotomic field, and not the field of $m$-adic numbers, and~$\Z_m$\index{$\Z_m$} denotes the set~$\Z/m\Z$ of integers modulo~$m$, and not the ring of $m$-adic integers.
The greatest common divisor of two integers~$m$ and~$l$ is denoted by~$\gcd(m,l)$ and is always chosen to be nonnegative.

Furthermore, $H$ denotes a Hopf algebra of finite dimension~$n$ with coproduct~$\dh$\index{$\Delta$}, counit~$\eh$\index{$\varepsilon$}, and antipode~$\sh$\index{$S$}. We will use the same symbols to denote the corresponding structure elements of the dual Hopf algebra~$H^*$, except for the antipode, which is denoted by~$\sh^*$. The opposite Hopf algebra, in which the multiplication is reversed, is denoted by~$H^{\op}$, and the coopposite Hopf algebra, in which the comultiplication is reversed, is denoted by~$H^{\cop}$. If $b_1,\ldots,b_n$ is a basis of~$H$ with dual basis $b^*_1,\ldots,b^*_n$, we have the formulas\endnote{\cite{ResSem}, Prop.~2.3', p.~542.}
\begin{align*}
\sum_{i=1}^n &b^*_i \o b_{i\1} \o b_{i\2} \o \ldots \o b_{i(m)} = \\
&\sum_{i_1,i_2,\ldots,i_m=1}^n b^*_{i_1} b^*_{i_2}  \cdots  b^*_{i_m}\o b_{i_1} \o b_{i_2} \o \ldots \o b_{i_m} 
\end{align*}
and
\begin{align*}
\sum_{i=1}^n & b^*_{i\1} \o b^*_{i\2} \o \ldots \o b^*_{i(m)} \o b_i= \\
&\sum_{i_1,i_2,\ldots,i_m=1}^n b^*_{i_1} \o  b^*_{i_2} \o \ldots \o b^*_{i_m} \o b_{i_1} b_{i_2} \cdots b_{i_m} 
\end{align*}
which we will refer to as the dual basis formulas. We use the letter~$A$ instead of~$H$ if the Hopf algebra under consideration is quasitriangular. 
With respect to enumeration, we use the convention that propositions, definitions, and similar items are referenced by the paragraph in which they occur; they are only numbered separately if this reference is ambiguous.

The essential part of the present work was carried out when the first author held a visiting research position at the Hong Kong University of Science and Technology. He thanks the university, and in particular his host, for the hospitality. He also thanks the University of Cincinnati
for a follow-on visiting position during which part of the manuscript was written. The second author would like to express his appreciation for the support by the RGC Competitive Earmarked Research Grant HKUST 6059/04.

\newpage
\section{The modular group} \label{Sec:ModGroup}
\subsection[Generators and relations]{} \label{GenRel}
In this article, the modular group is defined as the group $\Gamma := \SL(2,\Z)$ of $2 \times 2$-matrices with integer entries and determinant~$1$; note that many authors define it as the quotient group $\PSL(2,\Z)$ instead. The modular group is generated by the two matrices\endnote{\cite{Apos}, Sec.~2.2, Thm.~2.1, p.~28; \cite{KK}, \S~II.2, p.~108.}
$$\gv := \begin{pmatrix} 0 & -1 \\ 1 & 0 \end{pmatrix} 
\qquad \text{and} \qquad
\gt := \begin{pmatrix} 1 & 1 \\ 0 & 1 \end{pmatrix}$$
It is easy to see that these matrices satisfy the relations
$$\gv^4 = 1 \qquad (\gt\gv)^3 = \gv^2$$
however, it is a nontrivial result that these are defining relations for the modular group.\endnote{\cite{CoxMos}, \S~7.2, p.~85; \cite{KleinFricke}, \S~II.9.1, p.~454; \cite{Maass}, Sec.~II.1, Thm.~8, p.~53; \cite{Magnus}, Thm.~3.1, p.~108.}

It is possible to replace the generator~$\gv$ by the generator
$$\gr := \gt^{-1} \gv^{-1} \gt^{-1} = 
\begin{pmatrix} 1 & 0 \\ -1 & 1 \end{pmatrix}$$
The generators~$\gr$ and~$\gt$ satisfy the relations
$$\gt \gr \gt = \gr \gt \gr \qquad (\gr \gt)^6 = 1$$ 
and it follows from the corresponding result for the preceding generators that this also constitutes a presentation of the modular group in terms of generators and relations. From this, we get that 
$\gv^{-1} \gr = \gt \gr \gt \gr = \gt \gv^{-1}$, which means that $\gr=\gv \gt\gv^{-1}$, so that the generators~$\gr$ and~$\gt$ are conjugate.

The matrix
$$\ga := \begin{pmatrix}1 & 0 \\ 0 & -1 \end{pmatrix}$$
is not contained in $\Gamma$, but conjugation by~$\ga$ induces an automorphism of~$\Gamma$, for which we introduce the following notation:
\begin{defn}
For $g \in \Gamma$, we define~$\tilde{g} := \ga g \ga^{-1} = \ga g \ga$.
\end{defn}
Note that we have
$$\begin{pmatrix}1 & 0 \\ 0 & -1 \end{pmatrix}
\begin{pmatrix}a & b \\ c & d \end{pmatrix}
\begin{pmatrix}1 & 0 \\ 0 & -1 \end{pmatrix} = 
\begin{pmatrix}a & -b \\ -c & d \end{pmatrix}$$
so that $\tilde{g} = g^{-1}$ whenever~$a=d$. In particular, we find for the special matrices that we have used above as generators that
$$\tilde{\gv} = \gv^{-1} \qquad \tilde{\gt} = \gt^{-1} \qquad 
\tilde{\gr} = \gr^{-1}$$

\subsection[Congruence subgroups]{} \label{CongSubgr}
If $N$ is a natural number, the quotient map from~$\Z$ to
$\Z_N:=\Z/N\Z$ induces a group homomorphism
$$\SL(2,\Z) \rightarrow \SL(2,\Z_N)$$
by applying the quotient map to every component of the matrix. The kernel of this map is denoted by~$\Gamma(N)$ and called the principal congruence subgroup of level~$N$. In other words, we have
$$\Gamma(N) := \{\begin{pmatrix}
a & b \\
c & d
\end{pmatrix} \in \SL(2,\Z) \mid 
a \equiv d \equiv 1, \; b \equiv c \equiv 0 \pmod{N}\}$$
In particular, we have~$\Gamma(1)=\Gamma$. A subgroup of the modular group is called a congruence subgroup if it contains~$\Gamma(N)$ for a suitable~$N$, and the smallest such~$N$ is called the level of the congruence subgroup.

The modular group acts naturally on the lattice~$\Z^2 := \Z \times \Z$.
The orbits of the principal congruence subgroups can be described as follows:
\begin{prop}
Two nonzero lattice points $(m,l), (m',l') \in \Z^2$ are in the same $\Gamma(N)$-orbit if and only if
$t := \gcd(m,l) = \gcd(m',l')$ and
$$m/t \equiv m'/t \pmod{N} \qquad \qquad l/t \equiv l'/t \pmod{N}$$
\end{prop}
\begin{pf}
If $(m,l)$ and $(m',l')$ are in the same $\Gamma(N)$-orbit,
so that
$$\begin{pmatrix}
m' \\ l'
\end{pmatrix} = 
\begin{pmatrix}
a & b \\
c & d
\end{pmatrix}
\begin{pmatrix}
m \\ l
\end{pmatrix}$$
then we have for the ideals of~$\Z$ generated by $m,l$ resp.\ $m',l'$ that
$$(m',l') = (am + bl, cm + dl) \subset (m,l) = (t)$$
and vice versa, so that the first assertion holds. If we divide the above relation by~$t$, we see that $(m/t,l/t)$ and $(m'/t,l'/t)$ are still in the same $\Gamma(N)$-orbit, and if we reduce this relation modulo~$N$, we see that they are componentwise congruent.

For the converse, we can assume that $t=1$. If now two pairs~$(m,l)$ and $(m',l')$ of relatively prime integers are componentwise congruent modulo~$N$, this also holds for the pairs $g(m,l)$ and $g(m',l')$
for any $g \in \Gamma$, and if we can show that $g(m,l)$ and $g(m',l')$
are in the same $\Gamma(N)$-orbit, then this also holds for the original pair $(m,l)$ and $(m',l')$, as $\Gamma(N)$ is a normal subgroup.

Now as $m$ and $l$ are relatively prime, we can find integers 
$n$ and~$k$ satisfying $mn+lk=1$, so that
$$\begin{pmatrix}
m \\ l
\end{pmatrix}=
\begin{pmatrix}
m & -k \\ l & n
\end{pmatrix}
\begin{pmatrix}
1 \\ 0
\end{pmatrix}$$
In other words, $(m,l)$ and $(1,0)$ are in the same $\Gamma$-orbit,
so that we can in fact assume $(m,l) = (1,0)$. This means that we only have to show that a pair~$(m',l')$ of relatively prime integers of the form $(m',l')=(1+aN,bN)$ is in the same $\Gamma(N)$-orbit as~$(1,0)$. But this means that we have to find integers $c,d \in \Z$ so that
$$\begin{pmatrix}
m' \\ l'
\end{pmatrix}=
\begin{pmatrix}
1+aN & cN \\ bN & 1+dN
\end{pmatrix}
\begin{pmatrix}
1 \\ 0
\end{pmatrix}$$
subject to the determinant condition
$$1 = (1+aN)(1+dN) - (bN)(cN) = 1 + aN + dN + ad N^2 - bc N^2$$
or alternatively 
$$0 = a + d(1 + a N) - bc N= a + dm' - cl'$$
As $m'$ and $l'$ are relatively prime, this equation is solvable. 
\qed
\end{pf}
Note that this proposition shows that the condition to be componentwise congruent modulo~$N$ is not sufficient for two lattice points to be in the same \mbox{$\Gamma(N)$-orbit}, as the pairs $(2,4)$ and $(5,7)$ illustrate for $N=3$. Furthermore, it should be noted that in the case $N=1$ it yields the following fact:
\begin{corollary}
Two nonzero lattice points $(m,l), (m',l') \in \Z^2$ are in the same $\Gamma$-orbit if and only if $\gcd(m,l) = \gcd(m',l')$.
\end{corollary}

\subsection[Orbits and congruence relations]{} \label{OrbCongRel}
The group homomorphism from~$\SL(2,\Z)$ to~$\SL(2,\Z_N)$ discussed at the beginning of Paragraph~\ref{CongSubgr} is surjective. The proof of this fact uses the following lemma, which we will use below for a different purpose:\endnote{\cite{KK}, Kap.~II, \S~3, p.~116.} 
\begin{lemma}
Suppose that $m$, $l$, and~$N$ are relatively prime integers and 
that $l \neq 0$. Then there exists an integer~$k \in \Z$ such that
$m+kN$ is relatively prime to~$l$. 
\end{lemma}

We need to introduce another subgroup of the modular group. We denote
by~$\Delta(N)$ the subgroup of~$\Gamma$ that is generated by all conjugates~$g\gt^N g^{-1}$ of~$\gt^N$ for~$g \in \Gamma$. This subgroup is obviously normal, and it follows from the discussion in Paragraph~\ref{GenRel} that it contains~$\gr^N$. Since $\Gamma(N)$ is a normal subgroup that contains~$\gt^N$, we have that~$\Delta(N)$ is contained in~$\Gamma(N)$. However, $\Delta(N)$ is strictly smaller than~$\Gamma(N)$ if~$N \ge 6$, and it is not even a congruence subgroup in this case.\endnote{\cite{KleinFricke}, \S~II.7.5, p.~397; \cite{Kno}, p.~96.} To deal with this difficulty, we adapt the following notion from the theory of monoids to our situation:\endnote{\cite{JacBasI}, Def.~1.4, p.~54.}
\begin{defn}
An equivalence relation on the lattice~$\Z^2$ is called a congruence relation if, for all $g \in \Gamma$, the lattice points~$g.(m,l)$ and $g.(n,k)$ are equivalent whenever the lattice points~$(m,l)$ and~$(n,k)$ are equivalent.
\end{defn}
From every normal subgroup of the modular group, we get a congruence relation by defining that two lattice points are equivalent if they are in the same orbit under the action of the normal subgroup. In this way, both~$\Gamma(N)$ and~$\Delta(N)$ give rise to congruence relations.

Considering relations as sets of pairs, one can show as in the case of monoids that the intersection of congruence relations is again a congruence relation.\endnote{\cite{JacBasI}, Par.~1.8, Exerc.~8, p.~57.} Therefore, for every relation there is a smallest congruence relation that contains this relation, namely the intersection of all congruence relations that contain the given relation. In this sense, we now consider the smallest congruence relation~$\sim$ on the lattice~$\Z^2$ that
has the following two properties:
\begin{enumerate}
\item
We have $(m,l) \sim \gt^N.(m,l)$.

\item 
We have $(m,l) \sim (m,kl)$ for every $k \in \Z$ that satisfies
$k \equiv 1 \pmod{N}$ and $\gcd(m,kl) = \gcd(m,l)$.
\end{enumerate}

The second property appears to be asymmetrical with respect to the two components. This is, however, not the case, because if~$k$ satisfies
$k \equiv 1 \pmod{N}$ and $\gcd(km,l) = \gcd(m,l)$, we have
$$\begin{pmatrix}
m \\ l
\end{pmatrix} = 
\begin{pmatrix}
0 & -1 \\
1 & 0
\end{pmatrix}
\begin{pmatrix}
l \\ -m
\end{pmatrix} \sim
\begin{pmatrix}
0 & -1 \\
1 & 0
\end{pmatrix}
\begin{pmatrix}
l \\ -km
\end{pmatrix} = \begin{pmatrix}
km \\ l
\end{pmatrix}$$
so that $(m,l) \sim (km,l)$.

We will need another property of the congruence relation~$\sim$:
\begin{prop}
For every integer~$n \in \Z$, we have 
$(nm,nl) \sim (nm',nl')$ whenever $(m,l) \sim (m',l')$.
\end{prop}
\begin{pf}
This is obvious if $n=0$, so let us assume that $n \neq 0$. 
Suppose that~$\approx$ is an arbitrary congruence relation that satisfies the two defining properties of~$\sim$, i.e., that satisfies  $(m,l) \approx \gt^N.(m,l)$ and $(m,l) \approx (m,kl)$ for every $k \in \Z$ with the properties $k \equiv 1 \pmod{N}$ and $\gcd(m,kl) = \gcd(m,l)$. Recall that~$\sim$ is the intersection of all such congruence relations.
We define a new relation~$\approx_n$ by setting
$$(m,l) \approx_n (m',l') :\Leftrightarrow (nm,nl) \approx (nm',nl')$$
It is immediate that this is again a congruence relation. It also satisfies the first defining property, namely that 
$(m,l) \approx_n \gt^N.(m,l)$.  For the
second property, note that if $k \in \Z$  satisfies
$k \equiv 1 \pmod{N}$ and $\gcd(m,kl) = \gcd(m,l)$, it also satisfies
$\gcd(nm,knl) = \gcd(nm,nl)$, so that $(nm,nl) \approx (nm,knl)$ and therefore~$(m,l) \approx_n (m,kl)$.

This shows that $(m,l) \sim (m',l')$ implies $(m,l) \approx_n (m',l')$,
which means that $(nm,nl) \approx (nm',nl')$. As this holds for all such congruence relations~$\approx$, we get $(nm,nl) \sim (nm',nl')$, as asserted.
\qed
\end{pf}

It is not hard to see that, if we had dropped the second defining property above, the congruence relation that would have arisen would have been exactly the one determined by the group~$\Delta(N)$ as described above. The following theorem asserts that, by incorporating the second property, we get exactly the congruence relation determined by the group~$\Gamma(N)$:
\begin{thm}
Two lattice points~$(m,l)$ and~$(m',l')$ are in the same $\Gamma(N)$-orbit
if and only if $(m,l) \sim (m',l')$.
\end{thm}
\begin{pf}
\begin{list}{(\arabic{num})}{\usecounter{num} \leftmargin0cm \itemindent5pt}
\item
Let us first show that equivalent lattice points are in the same $\Gamma(N)$-orbit. For this, we need to look at the two defining properties of our congruence relation. For the first property, 
it is obvious that~$(m,l)$ and~$\gt^N.(m,l)$ are in the same $\Gamma(N)$-orbit. For the second property, suppose that $k \in \Z$ satisfies $k \equiv 1 \pmod{N}$ and $t:= \gcd(m,l) = \gcd(m,kl)$. If $(m,l)$ is nonzero, we have that~$(m/t,l/t)$ and~$(m/t,kl/t)$ are componentwise congruent modulo~$N$. By Proposition~\ref{CongSubgr}, this implies that $(m,l)$ and~$(m,kl)$ are in the same $\Gamma(N)$-orbit.  Clearly, this is also the case if~$(m,l)=(0,0)$.

This shows that the congruence relation determined by~$\Gamma(N)$, for which the equivalence classes are exactly the $\Gamma(N)$-orbits, takes part in the intersection that was used to define the relation~$\sim$. In other words, if $(m,l) \sim (m',l')$, then~$(m,l)$ and~$(m',l')$ are in the same $\Gamma(N)$-orbit.

\item
Now suppose that~$(m,l)$ and~$(m',l')$ are in the same $\Gamma(N)$-orbit.
In the case $N=1$, we have $\Gamma(N)=\Delta(N)=\Gamma$, and we have already pointed out above that the two lattice points are then equivalent. We will therefore assume in the sequel that~$N>1$.

We first consider the case where~$m$ and~$l$ are relatively prime; by Corollary~\ref{CongSubgr}, we then also have that~$m'$ and~$l'$ are relatively prime. We need a couple of reductions.
The first reduction is that we can assume in addition that all the components~$m$, $l$, $m'$, and~$l'$ are also relatively prime to~$N$. To see this, choose two distinct primes~$p$ and~$q$ that do not divide~$N$.
By Corollary~\ref{CongSubgr}, we can then find $g \in \Gamma$ such that
$g.(m,l)=(p,q)$. If we define $(p',q') := g.(m',l')$, then $(p,q)$ and~$(p',q')$ are also in the same~$\Gamma(N)$-orbit, because~$\Gamma(N)$ is a normal subgroup. By Proposition~\ref{CongSubgr}, this implies that
$p'$ and $q'$ are relatively prime and that
$$p \equiv p' \pmod{N} \qquad q \equiv q' \pmod{N}$$
so that in particular also $p'$ and $q'$ are relatively prime to~$N$.
But if we could establish that $(p,q)$ and~$(p',q')$ are equivalent, then also~$(m,l)$ and~$(m',l')$ would be equivalent, because~$\sim$ is a congruence relation. Therefore, we can assume from the beginning that all the components~$m$, $l$, $m'$, and~$l'$ are also relatively prime to~$N$.
Note that this implies in particular that the components are nonzero.
This completes our first reduction.

\item
The second reduction is that we can assume in addition that~$m$ is relatively prime to~$l'$ and that~$m'$ is relatively prime to~$l$. Now the numbers $m$ and~$Nl$ are relatively prime, which obviously implies that the numbers $m$, $Nl$, and~$l'$ are relatively prime.
We can therefore apply the lemma stated at the beginning of the paragraph to find an integer $k \in \Z$ such that~$m+kNl$ is relatively prime to~$l'$. Note that $m+kNl$ is still relatively prime to~$N$ and~$l$, and that $(m+kNl,l)=\gt^{kN}.(m,l)$ is equivalent to~$(m,l)$ even if~$k$ is negative.
By replacing $(m,l)$ by $(m+kNl,l)$, we can therefore assume from the beginning that in addition~$m$ and~$l'$
are relatively prime. Using the same argument with the lattice points interchanged, we can furthermore assume that~$m'$ and~$l$ are relatively prime.

\item
The third reduction is that we can assume in addition that~$m=m'$. We have assumed that $m'$ and $N$ are relatively prime, which implies that there is a number~$k' \in \Z$ such that $m'k' \equiv 1 \pmod{N}$. The numbers~$k'$ and~$N$ are then relatively prime, which obviously implies that the numbers~$k'$, $N$, and~$ll'$ are relatively prime. As all components are nonzero, we can apply the above lemma again to find an integer~$k \in \Z$ such that $n:=k'+kN$ is relatively prime to~$ll'$. As we also have $nm' \equiv k'm' \equiv 1 \pmod{N}$, we get by the variant of the second defining property of our congruence relation discussed above that 
$(m,l) \sim (nm'm,l)$. But $m \equiv m' \pmod{N}$ by Proposition~\ref{CongSubgr}, so that
$$nm \equiv nm' \equiv 1 \pmod{N}$$
and furthermore $nm$ and $l'$ are relatively prime. Again by the variant of the second defining property, we therefore see that
$(m',l') \sim (nmm',l')$. By replacing~$(m,l)$ by~$(nm'm,l)$
and~$(m',l')$ by~$(nmm',l')$, we can therefore reduce to the situation where $m=m'$.

\item
We have now two lattice points $(m,l)$ and $(m,l')$ with relatively prime components~$m$ and~$l$ resp.~$m$ and~$l'$. Moreover, all of these components are relatively prime to~$N$, and in particular nonzero.
By assumption, they are in the same $\Gamma(N)$-orbit, so that~$l \equiv l' \pmod{N}$ by Proposition~\ref{CongSubgr}. We have to establish that they are equivalent. 

For this, we argue as in the preceding step: Choose
$k$ such that $kl \equiv 1 \pmod{N}$. Then the numbers $k$ and $N$ are relatively prime, which clearly implies that the numbers $k$, $N$, and~$m$ are relatively prime. Therefore, again by the above lemma,
we can find an integer~$n' \in \Z$ such that $n:=k+n'N$ is relatively prime to~$m$. We then have $n \equiv k \pmod{N}$ and therefore
$$nl \equiv nl' \equiv 1 \pmod{N}$$
and $m$ and $nl$ resp.~$nl'$ are relatively prime. By the second defining property of our congruence relation, we have that $(m,l) \sim (m,nl'l)$, and similarly that
$(m,l') \sim (m,nll')$. As equal pairs are clearly equivalent, this finishes the proof in the case of lattice points with relatively prime components.

\item
We now consider the general case, in which we have two lattice points~$(m,l)$ and~$(m',l')$ in the same $\Gamma(N)$-orbit, but~$m$ and~$l$ are not necessarily relatively prime. We have to establish that they are equivalent, and we can clearly assume that they are different from the origin. By Proposition~\ref{CongSubgr}, we have
$t := \gcd(m,l) = \gcd(m',l')$ and
$$m/t \equiv m'/t \pmod{N} \qquad \qquad l/t \equiv l'/t \pmod{N}$$
Now $m/t$ and~$l/t$ are relatively prime, and $m'/t$ and~$l'/t$ are relatively prime as well. Furthermore, $(m/t,l/t)$ and~$(m'/t,l'/t)$ are in the same $\Gamma(N)$-orbit. By the facts already established, we therefore get $(m/t,l/t) \sim (m'/t,l'/t)$, which implies 
$(m,l) \sim (m',l')$ by the above proposition.
\qed
\end{list}
\end{pf}

\subsection[Presentations of the factor groups]{} \label{PresFac}
\vspace{-2mm}
The groups $\SL(2,\Z_N) \cong \SL(2,\Z)/\Gamma(N)$ are obviously 
generated by the images of the generators under the canonical map, which are
$$\bv := \begin{pmatrix} \bar{0} & -\bar{1} \\ \bar{1} & \bar{0} \end{pmatrix} 
\qquad \text{and} \qquad
\bt := \begin{pmatrix} \bar{1} & \bar{1} \\ \bar{0} & \bar{1} \end{pmatrix}$$
or alternatively~$\bt$ and
\vspace{-0.5mm}
$$\br := \bt^{-1} \bv^{-1} \bt^{-1} = 
\begin{pmatrix} \bar{1} & \bar{0} \\ -\bar{1} & \bar{1} \end{pmatrix}$$
However, the defining relations for these generators are not easy to obtain. To write them down, we introduce the abbreviation
$\gd(q):= \bv \bt^{q'} \bv^{-1} \bt^q \bv \bt^{q'}$
for $q,q' \in \Z$ such that $qq' \equiv 1 \pmod{N}$.
Although we want to understand this expression here as an abbreviation
for a word in the generators, it is of course also possible to compute
the corresponding matrix in~$\SL(2,\Z_N)$:
\begin{align*}
\gd(q) &= 
\begin{pmatrix} \bar{0} & -\bar{1} \\ \bar{1} & \bar{0} \end{pmatrix} 
\begin{pmatrix} \bar{1} & \bar{q}' \\ \bar{0} & \bar{1} \end{pmatrix}
\begin{pmatrix} \bar{0} & \bar{1} \\ -\bar{1} & \bar{0} \end{pmatrix} 
\begin{pmatrix} \bar{1} & \bar{q} \\ \bar{0} & \bar{1} \end{pmatrix}
\begin{pmatrix} \bar{0} & -\bar{1} \\ \bar{1} & \bar{0} \end{pmatrix} 
\begin{pmatrix} \bar{1} & \bar{q}' \\ \bar{0} & \bar{1} \end{pmatrix} \\
&= 
\begin{pmatrix} \bar{0} & -\bar{1} \\  \bar{1} & \bar{q}' \end{pmatrix}
\begin{pmatrix} \bar{0} & \bar{1} \\ -\bar{1} & -\bar{q}\end{pmatrix}
\begin{pmatrix} \bar{0} & -\bar{1}  \\ \bar{1} & \bar{q}'\end{pmatrix} 
= 
\begin{pmatrix} \bar{1} & \bar{q} \\ -\bar{q}' & \bar{0} \end{pmatrix}
\begin{pmatrix} 0 & -\bar{1}  \\ \bar{1} & \bar{q}'\end{pmatrix}
= \begin{pmatrix} \bar{q} & \bar{0} \\ \bar{0} & \bar{q}' \end{pmatrix}
\end{align*}
Note also that we have intentionally suppressed the dependence of~$\gd(q)$ on~$q'$, on which it, as an abbreviation for a word in the generators, in principle depends. 

The following proposition, which is adapted from~\cite{CostGann}, lists one possible set of defining relations:
\vspace{-2mm}
\begin{prop}
Write $N=2^e m$, where~$m$ is odd. Then the relations
\vspace{-4mm}
\begin{enumerate}
\item $\bv^4 = 1 \qquad (\bt\bv)^3 = \bv^2 \qquad \bt^N = 1$

\item $\bt^{2^e} (\bv \bt^m \bv^{-1}) = (\bv \bt^m \bv^{-1}) \bt^{2^e}$

\item $\gd(q) \bv = \bv \gd(q)^{-1}$

\item $\gd(q) \bt = \bt^{q^2} \gd(q)$
\end{enumerate}
\vspace{-3.5mm}
for all $q \in \Z$ that are relatively prime to~$N$, are defining relations for~$\SL(2,\Z_N)$.
\end{prop}
\begin{pf}
This is proved in~\cite{CostGann}, \S~2.2, Lem.~1.c, p.~5, where also further references are given. Note that the generator~$\gv$ is defined differently in~\cite{CostGann}, namely as our~$\gv^{-1}$. It is also shown there that the relations~3 and~4 are not necessary for all~$q$ that are relatively prime to~$N$, but only for $q=1-2d$, $q = 2-d$, and $q=2d+1$, where~$d$ is an integer that satisfies
$d \equiv 1 \pmod{2^e}$ and $d \equiv 0 \pmod{m}$.
\qed
\end{pf}

\newpage
\section{Quasitriangular Hopf algebras} \label{Sec:QuasitriHopf}
\subsection[Quasitriangular Hopf Algebras]{} \label{QuasitriHopf}
Recall that a Hopf algebra~$A$ is called 
quasitriangular\endnote{\cite{Kas}, Def.~VIII.2.2, p.~173; \cite{M}, Def.~10.1.5, p.~180; \cite{Tur}, Sec.~XI.2.1, p.~496.} if its antipode is invertible and it possesses a so-called R-matrix, which is an invertible element 
$R=\sum_{i=1}^m a_i \o b_i \in A \o A$ that satisfies
$\da^{\cop}(a) = R \da(a) R^{-1}$ as well as
$$(\da \o \id)(R) = \sum_{i,j=1}^m a_i \o a_j \o b_i b_j \qquad
(\id \o \da)(R) = \sum_{i,j=1}^m a_i a_j \o b_j \o b_i$$
Associated with the R-matrix is the Drinfel'd element
$\ua:=\sum_{i=1}^m \sa(b_i) a_i$. This is an invertible element that satisfies\endnote{\cite{Kas}, Sec.~VIII.4, p.~179; \cite{M}, Thm.~10.1.13, p.~181; \cite{Tur}, Sec.~XI.2.2, p.~498.}
$$\da(\ua)= (\ua \o \ua)(R'R)^{-1} = (R'R)^{-1} (\ua \o \ua)
\qquad \sa^2(a) = \ua a \ua^{-1}$$
where $R':=\sum_{i=1}^m b_i \o a_i$ arises from the R-matrix by interchanging the tensorands. The inverse Drinfel'd element is given by
$\ua^{-1}=\sum_{i=1}^m \sa^{-2}(b_i) a_i$. In this context, it should be noted that the element~$R'^{-1}$ always also is an \mbox{R-matrix} for~$A$. The Hopf algebra is called triangular if these two choices for the \mbox{R-matrices} coincide.

\subsection[The Drinfel'd double construction]{}  \label{DrinfDouble}
An important source of quasitriangular Hopf algebras is the Drinfel'd double construction.\endnote{\cite{Kas}, Chap.~IX, p.~199; \cite{M}, \S~10.3, p.~187; \cite{Tur}, Sec.~XI.2.4, p.~499.} For an arbitrary finite-dimensional Hopf algebra~$H$, the Drinfel'd double~$D:=D(H)$ is a Hopf algebra whose underlying vector space is $H^* \o H$. The coalgebra structure is the tensor product coalgebra structure $H^{* \cop} \o H$,
so that coproduct and counit are given by the formulas
$$\dd(\varphi \o h) = (\varphi_\2 \o h_\1) \o (\varphi_\1 \o h_\2)
\qquad \ed(\varphi \o h) = \varphi(1) \varepsilon(h)$$
The formula for the product is a little more involved; it reads
$$(\varphi \o h)(\varphi' \o h') = 
\varphi'_\1(\sh^{-1}(h_\3)) \varphi'_\3(h_\1) \;
\varphi \varphi'_\2 \o h_\2 h'$$
Finally, the antipode is given by the formula
$\sd(\varphi \o h) = (\varepsilon \o \sh(h))(\sh^{-1*}(\varphi) \o 1)$.

To establish the assertion that the Drinfel'd double is quasitriangular, we have to endow it with an R-matrix, which is explicitly given as follows: If $b_1,\ldots,b_n$ is a basis of~$H$ with dual basis $b^*_1,\ldots,b^*_n$,
then the R-matrix\index{$R$} is 
$$R = \sum_{i=1}^n (\varepsilon \o b_i) \o (b^*_i \o 1)$$
The associated Drinfel'd element~$\ud$\index{$\ud$} and its inverse are therefore
$$\ud = \sum_{i=1}^n \sh^{-1*}(b^*_i) \o b_i \qquad 
\ud^{-1} = \sum_{i=1}^n \sh^{2*}(b^*_i) \o b_i$$
The Drinfel'd element has an analogue in the dual~$D^*$, namely the evaluation form
$$e: D \rightarrow K,~\varphi \o h \mapsto \varphi(h)$$
The evaluation form is a symmetric Frobenius homomorphism.\endnote{\cite{Kerl}, Prop.~7, p.~366; \cite{YYY1}, Par.~2, p.~89.} It is invertible with inverse
$e^{-1}(\varphi \o h)= \varphi(\sh^{-1}(h))$.

\subsection[Integrals of the Drinfel'd double]{} \label{IntDrinfDouble}
The integrals of the Drinfel'd double can be described in terms of the integrals of the original Hopf algebra~$H$. If we choose left integrals $\Lmh \in H$ and $\lmh \in H^*$ as well as right integrals $\Gh \in H$ and
$\rhh \in H^*$, then
$$\Lmd = \lmh \o \Gh$$
is a two-sided integral of the Drinfel'd double, which in particular tells that the Drinfel'd double is unimodular.\endnote{\cite{M}, Thm.~10.3.12, p.~192; \cite{R2}, Thm.~4, p.~303.} Similarly, the functions~$\lmd$ and~$\rhd$ in~$D^*$ defined by 
$$\lmd(\varphi \o h) = \varphi(\Gh) \lmh(h) \qquad
\rhd(\varphi \o h) = \varphi(\Lmh) \rhh(h)$$
are left resp.~right integrals on~$D$. Using the forms of the Drinfel'd element and its inverse given in Paragraph~\ref{DrinfDouble}, we see that
\begin{align*}
&\lmd(\ud) = \lmh(\sh^{-1}(\Gh)) \qquad 
\lmd(\ud^{-1}) = \lmh(\sh^{2}(\Gh)) \\
&\rhd(\ud) = \rhh(\sh^{-1}(\Lmh)) \qquad 
\rhd(\ud^{-1}) = \rhh(\sh^{2}(\Lmh))
\end{align*}

Using these integrals, it is possible to relate the Drinfel'd element~$\ud$ and the evaluation form~$e$:
\begin{lemb}
\begin{enumerate}
\item 
$\Lmd_\1 e(\Lmd_\2) = e(\Lmd) \ud$

\item
$e(\Lmd_\1) \Lmd_\2 = e(\Lmd) \sd(\ud)$

\item
$e^{-1}(\Lmd_\1) \Lmd_\2 = e^{-1}(\Lmd) \ud^{-1}$

\item
$\Lmd_\1 e^{-1}(\Lmd_\2) = e^{-1}(\Lmd) \sd(\ud^{-1})$
\end{enumerate}
\end{lemb}
\begin{pf}
For every $h \in H$, we have 
\begin{align*}
\lmh_\1(\Gh_\2) \lmh_\2(h) \; \Gh_\1 
&= \lmh(\Gh_\2 h) \; \Gh_\1 =
\lmh(\Gh_\2 h_\3) \; \Gh_\1 h_\2 \sh^{-1}(h_\1) \\
&= \lmh(\Gh_\2) \; \Gh_\1 \sh^{-1}(h)
= \lmh(\Gh) \; \sh^{-1}(h)
= e(\Lmd) \; \sh^{-1}(h)
\end{align*}
Therefore, if $b_1,\ldots,b_n$ is a basis of~$H$ with dual basis $b^*_1,\ldots,b^*_n$, we have
\begin{align*}
&\Lmd_\1 e(\Lmd_\2) = \lmh_\2 \o \Gh_\1 \; e(\lmh_\1 \o \Gh_\2)
= \lmh_\2 \o \Gh_\1 \; \lmh_\1(\Gh_\2) \\
&= \sum_{i=1}^n b^*_i \o \lmh_\2(b_i) \Gh_\1 \; \lmh_\1(\Gh_\2)
= e(\Lmd)\sum_{i=1}^n b^*_i \o \sh^{-1}(b_i)
= e(\Lmd) \ud
\end{align*}
This proves the first relation. The second relation follows from this by applying the antipode~$\sd$, because~$\Lmd$ is invariant under the antipode,\endnote{\cite{R4}, Prop.~3, p.~590.} and we have~$\sd^*(e)=e$, since
\begin{align*}
\sd^*(e)(\varphi \o h) &= 
e((\eh \o \sh(h))(\sh^{-1*}(\varphi) \o \H)) \\
&= e((\sh^{-1*}(\varphi) \o \H)(\eh \o \sh(h))) 
= \varphi(h) 
\end{align*}
For the third relation, note that
\begin{align*}
&\lmh_\1(h) \lmh_\2(\sh^{-1}(\Gh_\1)) \sh^{-2}(\Gh_\2)
= \lmh(h \sh^{-1}(\Gh_\1)) \sh^{-2}(\Gh_\2) \\
&= \lmh(h_\2 \sh^{-1}(\Gh_\1)) h_\1 \sh^{-1}(\Gh_\2) \sh^{-2}(\Gh_\3)
= \lmh(h_\2 \sh^{-1}(\Gh)) h_\1 \\
&= \lmh(\sh^{-1}(\Gh)) h 
= e^{-1}(\Lmd) h 
\end{align*}
which implies
\begin{align*}
e^{-1}(\Lmd_\1) \Lmd_\2 &= 
e^{-1}(\lmh_\2 \o \Gh_\1) \; \lmh_\1 \o \Gh_\2 =
\lmh_\2(\sh^{-1}(\Gh_\1) \; \lmh_\1 \o \Gh_\2 \\
&= 
\sum_{i=1}^n \lmh_\1(b_i) \lmh_\2(\sh^{-1}(\Gh_\1)  \;  b_i^* \o \Gh_\2 \\
&= e^{-1}(\Lmd)\sum_{i=1}^n b_i^* \o \sh^2(b_i)
= e^{-1}(\Lmd) \ud^{-1}
\end{align*}
The fourth relation follows as before from the third by applying the antipode.~\qed
\end{pf}

We will also need the corresponding result that expresses the evaluation form in terms of the Drinfel'd element:\endnote{\cite{Kerl}, Prop.~7, p.~366.}
\begin{propb}
\begin{enumerate}
\item 
$\rhd(\ud x) = \rhd(x\ud) = \rhd(\ud) e(x)$

\item 
$\rhd(\sd(\ud^{-1})x) = \rhd(x\sd(\ud^{-1})) = \rhd(\sd(\ud^{-1})) e^{-1}(x)$
\end{enumerate}
\end{propb}
\begin{pf}
We can assume that $x = \varphi \o h$.
For the first relation, we then have
\begin{align*}
\rhd(\ud x) &= \rhd((\sh^{-2*}(\varphi) \o \H) \ud (\eh \o h)) =
\sum_{i=1}^n \rhd(\sh^{-2*}(\varphi)\sh^{-1*}(b_i^*) \o b_i h) \\
&= \sum_{i=1}^n \sh^{-2*}(\varphi)(\Lmh_\1) \sh^{-1*}(b_i^*)(\Lmh_\2) 
\rhh_\1(b_i) \rhh_\2(h)\\
&= \varphi(\sh^{-2}(\Lmh_\1)) \rhh_\1(\sh^{-1}(\Lmh_\2)) \rhh_\2(h)\\
&= \varphi(\sh^{-2}(\Lmh_\1)) \rhh(\sh^{-1}(\Lmh_\2) h) \\
&= \varphi(\sh^{-2}(\Lmh_\1) \sh^{-1}(\Lmh_\2) h_\2) 
\rhh(\sh^{-1}(\Lmh_\3) h_\1)\\
&= \varphi(h_\2) \rhh(\sh^{-1}(\Lmh) h_\1)
= \varphi(h) \rhh(\sh^{-1}(\Lmh)) = e(x) \rhd(\ud)
\end{align*}
This shows that $\rhd(\ud x)= \rhd(\ud) e(x)$; if we replace~$x$ by~$\ud^{-1}x\ud$ 
and use that~$e$ is a symmetric Frobenius homomorphism, we get that also
$\rhd(x\ud) = \rhd(\ud) e(x)$.

For the second relation, we have
\begin{align*}
\rhd(\sd(\ud^{-1})x) &= 
\sum_{i=1}^n \rhd(\sd(\eh \o b_i) \sd(\sh^{2*}(b_i^*) \o \H) x) \\
&= \sum_{i=1}^n \rhd((\eh \o b_i) (\sh^{2*}(b_i^*) \o \H) x)\\
&= \sum_{i=1}^n \rhd((\sh^{2*}(b_i^*) \o \H) x(\eh \o \sh^{-2}(b_i)))\\
&= \sum_{i=1}^n \rhd((b_i^* \o \H) x(\eh \o b_i))
= \sum_{i=1}^n \rhd(b_i^* \varphi \o h b_i) \\
&= \sum_{i=1}^n b_i^*(\Lmh_\1) \varphi(\Lmh_\2) \rhh_\1(h) \rhh_\2(b_i)
= \varphi(\Lmh_\2) \rhh_\1(h) \rhh_\2(\Lmh_\1)\\
&= \varphi(\Lmh_\2) \rhh(h \Lmh_\1)
= \varphi(\sh^{-1}(h_\3) h_\2 \Lmh_\2) \rhh(h_\1 \Lmh_\1)\\
&= \varphi(\sh^{-1}(h_\2)) \rhh(h_\1 \Lmh)
= \varphi(\sh^{-1}(h)) \rhh(\Lmh)
= e^{-1}(x) \rhh(\Lmh)
\end{align*}
For $x=\D$, this yields $\rhd(\sd(\ud^{-1})) = \rhh(\Lmh)$, which we can resubstitute in order to establish one of the claimed identities. The other one follows as before by substituting $\sd(\ud) x \sd(\ud^{-1})$
for~$x$ and using the symmetry of~$e$.
\qed
\end{pf}

\newpage
\section{Factorizable Hopf algebras} \label{Sec:FactorizeHopf}
\subsection[Doubles of quasitriangular Hopf algebras]{} \label{DoubleQuasitri}
If $A$ is already a quasitriangular Hopf algebra, it is of course also possible to form its double~$D=D(A)$. In this case, there exists a Hopf algebra retraction
$$\pi: D(A) \rightarrow A,~\varphi \o a \mapsto 
\bigl((\varphi \o \id)(R)\bigr) a$$
from the double of~$A$ to~$A$ itself.\endnote{\cite{DrinfAlmCocom}, Prop.~6.2, p.~337; \cite{Kas}, Prop.~VIII.2.5, p.~177; \cite{SchneiderFact}, Sec.~4, p.~1896; \cite{TsangZhu}, Thm.~1, p.~2.}
Using the alternative R-matrix~$R'^{-1}$ mentioned in Paragraph~\ref{QuasitriHopf} instead of~$R$, we get another Hopf algebra retraction~$\pi'$. As we have\endnote{\cite{Kas}, Thm.~VIII.2.4, p.~175; \cite{M}, Prop.~10.1.8, p.~180; \cite{Tur}, Lem.~XI.2.1.1, p.~497.} $R^{-1} = (\sa \o \id)(R) = (\id \o \sa^{-1})(R)$, this map is explicitly given as
$$\pi': D(A) \rightarrow A,~\varphi \o a \mapsto 
\bigl((\sa \o \varphi)(R)\bigr) a =
\bigl((\id \o \varphi)(\id \o \sa^{-1})(R)\bigr) a$$
From these two homomorphisms, we derive the algebra homomorphism
$$\Psi: D(A) \rightarrow A \o A,~x \mapsto (\pi \o \pi')(\dd(x))$$
where $A \o A$ carries the canonical algebra structure. Note that this map is in general not a coalgebra homomorphism with respect to the canonical coalgebra structure on~$A \o A$. However, it becomes a Hopf algebra homomorphism if we twist the comultiplication\endnote{\cite{ResSem}, Thm.~2.9, p.~546; \cite{SchneiderFact}, Thm.~4.3, p.~1897; \cite{TsangZhu}, Thm.~2, p.~2.} by the cocycle $F:=\A \o R^{-1} \o \A \in A^{\o 4}$. In other words, $\Psi$ is a Hopf algebra homomorphism if considered as a map to the Hopf algebra $(A \o A)_F$, which has the canonical tensor product algebra structure, but the twisted coproduct
$$\df(a \o b) = F ((a_\1 \o b_\1) \o (a_\2 \o b_\2)) F^{-1}$$
The Hopf algebra $(A \o A)_F$ even becomes quasitriangular by using the twisted R-matrix\endnote{\cite{Kas}, Prop.~XV.3.6, p.~376.} 
$R_F := F_t R_{A \o A} F^{-1}$, where
$$R_{A \o A} := \sum_{i,j=1}^m a_i \o b_j \o b_i \o \sa(a_j)$$
is an R-matrix for the tensor product Hopf algebra~$A \o A$ and $F_t$ arises from~$F$ by interchanging the first and the third as well as the second and the fourth tensor factor. By using this specific R-matrix, $\Psi$ becomes a morphism of quasitriangular Hopf algebras\endnote{\cite{M}, Def.~10.1.15, p.~183.} in the sense that it maps the R-matrix of the Drinfel'd double to~$R_F$. This implies that the image of the Drinfel'd element is given as follows: 
\begin{lemma}
$\Psi(\ud) = \ua \o \ua^{-1}$
\end{lemma}
\begin{pf}
In general, if the coproduct of a Hopf algebra is modified by a cocycle, then the resulting Hopf algebra has the antipode 
$S_F(a) = w \sa(a) w^{-1}$, where $w$ arises from $(\id \o \sa)(F)$ by multiplication of the tensorands.\endnote{\cite{Kas}, Exerc.~XV.6.1, p.~381; \cite{SchneiderFact}, p.~1897.} From this, we see that
$$S_F^2(a) = x \sa^2(a) x^{-1}$$
where $x=w \sa(w^{-1})$, and it can be shown that the Drinfel'd element~$u_F$ that arises from the R-matrix $R_F=F_t R F^{-1}$ is related to the original one via the formula $u_F = x \ua$.

In our case, we find 
$$w = \sum_{i=1}^m (\A \o \sa(a_i)) S_{A \o A}(b_i \o \A)
= \sum_{i=1}^m \sa(b_i) \o \sa(a_i) = R'$$
This implies in this case that~$x=1$, which means that the Drinfel'd elements of~$(A \o A)_{F}$ and $A \o A$ coincide. Because~$\Psi$ is a morphism of quasitriangular Hopf algebras, this element is equal to~$\Psi(\ud)$. But for the Drinfel'd element of~$A\o A$, we find
$$u_{A \o A} = 
\sum_{i,j=1}^m S_{A \o A}(b_i \o \sa(a_j))(a_i \o b_j)
= 
\sum_{i,j=1}^m \sa(b_i)a_i \o \sa^2(a_j) b_j
= \ua \o \sa^2(\ua^{-1})$$
Because the Drinfel'd element is invariant under the square of the antipode, this implies the assertion.
\qed
\end{pf}

Note that by replacing the R-matrix with the alternative R-matrix~$R'^{-1}$, we get a second Hopf algebra homomorphism
$$\Psi': D(A) \rightarrow (A \o A)_{F'},~x \mapsto (\pi' \o \pi)(\dd(x))$$
where this time we have to use the cocycle 
$F':=\A \o R' \o \A \in A^{\o 4}$ to twist the comultiplication on the right-hand side. As before, the Hopf algebra $(A \o A)_{F'}$ is quasitriangular with respect to the R-matrix 
$R_{F'} := F'_t R''_{A \o A} F'^{-1}$, where
$$R''_{A \o A} := \sum_{i,j=1}^m b_i \o a_j \o \sa(a_i) \o b_j$$
is an R-matrix for the tensor product Hopf algebra~$A \o A$ and $F'_t$ arises from~$F'$ by interchanging the first and the third as well as the second and the fourth tensor factor. By using this specific R-matrix, $\Psi'$ becomes a morphism of quasitriangular Hopf algebras in the sense that it maps the R-matrix of the Drinfel'd double to~$R_{F'}$. Furthermore, the Drinfel'd element arising from the alternative R-matrix~$R'^{-1}$ is exactly the inverse~$\ua^{-1}$ of the original one, so that the preceding lemma yields that $\Psi'(\ud) = \ua^{-1} \o \ua$.

\subsection[Factorizable Hopf algebras]{} \label{FactHopf}
In the situation of Paragraph~\ref{DoubleQuasitri}, we have a left action of $A \o A$ on~$A$ by requiring that the element~$a \o a' \in A \o A$
acts on~$b \in A$ by mapping it to~$a' b \sa^{-1}(a)$. Pulling this action back along~$\Psi'$, we get a left action of~$D(A)$ on~$A$ given by
$$x \twoheadrightarrow b = \pi(x_\2) b \sa^{-1}(\pi'(x_\1))$$
If $x = \varphi \o a$ and $R= \sum_{i=1}^m a_i \o b_i$, this action is explicitly given as
\begin{align*}
(\varphi \o a) \twoheadrightarrow b &= 
\pi(\varphi_\1 \o a_\2) b \sa^{-1}(\pi'(\varphi_\2 \o a_\1)) \\
&= \sum_{i,j=1}^m \varphi_\1(a_i) b_i a_\2 b 
\sa^{-1}(\varphi_\2(b_j) \sa(a_j) a_\1) \\
&= \sum_{i,j=1}^m \varphi(a_i b_j) b_i a_\2 b \sa^{-1}(a_\1) a_j
\end{align*}
We will use the same notation for the restrictions of this action to~$A$ and~$A^*$, i.e., we define
\begin{align*}
a \twoheadrightarrow b := 
(\ea \o a) \twoheadrightarrow b &= a_\2 b \sa^{-1}(a_\1) \\
\varphi \twoheadrightarrow b := (\varphi \o \A) \twoheadrightarrow b 
&= \sum_{i,j=1}^m \varphi(a_i b_j) b_i b a_j
\end{align*}
Note that the restriction of the action to~$A$ is just the left adjoint action of~$A^{\cop}$ on itself. The space of invariants for this restricted action therefore is exactly the center~$Z(A)$ of~$A$.

We also introduce the map
$$\Phi: A^* \rightarrow A,~\varphi \mapsto 
(\varphi \twoheadrightarrow \A) = 
\sum_{i,j=1}^m \varphi(a_i b_j) b_i a_j = (\id \o \varphi)(R'R)$$
If~$C(A)$ denotes the subalgebra
$$C(A) := \{\chi \in A^* \mid \chi(ab) = \chi(b\sa^2(a)) \; 
\text{for all} \; a, b \in A\}$$
of~$A^*$, it is known\endnote{\cite{DrinfAlmCocom}, Prop.~3.3, p.~327; see also \cite{EG1}, Lem.~1.1; p.~192; \cite{SchneiderFact}, Thm.~2.1, p.~1892.} that~$\Phi$ has the following property:
\begin{prop}
We have 
$$\Phi(\varphi \chi) = \Phi(\varphi) \Phi(\chi)$$
for all $\varphi \in A^*$ and all $\chi \in C(A)$. 
Furthermore, we have
$$a\twoheadrightarrow \Phi(\varphi)
= \varphi_\1(\sa^{-1}(a_\2)) \varphi_\3(a_\1) \; \Phi(\varphi_\2)$$
for all $\varphi \in A^*$ and all $a \in A$. 
Consequently, $\Phi$ restricts to an algebra homomorphism from~$C(A)$ to~$Z(A)$.
\end{prop}
\begin{pf}
It is possible to verify these properties by direct computation; however, it is interesting to derive them from our construction of~$\Phi$. The second equation holds since
\begin{align*}
a\twoheadrightarrow \Phi(\varphi) &=
(\ea \o a) (\varphi \o \A) \twoheadrightarrow \A
= \varphi_\1(\sa^{-1}(a_\3)) \varphi_\3(a_\1) \; (\varphi_\2 \o a_\2) \twoheadrightarrow \A \\
&=\varphi_\1(\sa^{-1}(a_\2)) \varphi_\3(a_\1) \; \Phi(\varphi_\2)
\end{align*}
This amounts to saying that~$\Phi$ is an $A$-linear map from~$A^*$ to~$A$, where~$A^*$ is considered as an $A$-module via 
$a \o \varphi \mapsto 
\varphi_\1(\sa^{-1}(a_\2)) \varphi_\3(a_\1) \varphi_\2$; this action is the left coadjoint action, built with the inverse antipode, and is actually used in the construction of the Drinfel'd double as a double crossproduct.\endnote{\cite{M}, Def.~10.3.1, p.~188.} The space~$C(A)$ is exactly the space of invariants for this action. As an $A$-linear map takes invariants to invariants, $\Phi$ maps~$C(A)$ to the center~$Z(A)$. 

For the first assertion, note that we clearly have
$\varphi \twoheadrightarrow z = (\varphi \twoheadrightarrow \A) z$
if $z \in Z(A)$, so that 
$$\Phi(\varphi \chi) = (\varphi \chi) \twoheadrightarrow \A
= \varphi \twoheadrightarrow (\chi \twoheadrightarrow \A)
= \varphi \twoheadrightarrow \Phi(\chi)
= (\varphi \twoheadrightarrow \A) \Phi(\chi)
= \Phi(\varphi) \Phi(\chi)$$
if $\chi \in C(A)$. Finally, note that it follows from the elementary properties of R-matrices\endnote{\cite{Kas}, Thm.~VIII.2.4, p.~175; \cite{M}, Prop.~10.1.8, p.~180; \cite{Tur}, Lem.~XI.2.1.1, p.~497.} that~$\Phi$ preserves the unit element. 
\qed
\end{pf}

In a very similar way, we have a right action of $A \o A$ on~$A$ by requiring that the element~$a \o a' \in A \o A$
acts on~$b \in A$ by mapping it to~$\sa^{-1}(a') b a$. Pulling this action back along~$\Psi$, we get a right action of~$D(A)$ on~$A$ given by
$$b \twoheadleftarrow x =  \sa^{-1}(\pi'(x_\2)) b \pi(x_\1)$$
The two actions are related via the formulas
$$\sa(x \twoheadrightarrow b) = \sa(b) \twoheadleftarrow \sd(x) 
\quad \text{and} \quad 
\sa(b \twoheadleftarrow x) = \sd(x) \twoheadrightarrow \sa(b)$$
If $x = \varphi \o a$ and $R= \sum_{i=1}^m a_i \o b_i$, this action is explicitly given as
\begin{align*}
b \twoheadleftarrow (\varphi \o a) &= 
\sa^{-1}(\pi'(\varphi_\1 \o a_\2)) b \pi(\varphi_\2 \o a_\1) \\
&= \sum_{i,j=1}^m \sa^{-1}(\varphi_\1(b_i) \sa(a_i) a_\2) b 
\varphi_\2(a_j) b_j a_\1 \\
&= \sum_{i,j=1}^m \varphi(b_i a_j) \sa^{-1}(a_\2) a_i b b_j a_\1
\end{align*}
As before, we use the same notation for the restrictions of this action to~$A$ and~$A^*$, so that
\begin{align*}
b \twoheadleftarrow a:= 
b \twoheadleftarrow (\ea \o a) &= \sa^{-1}(a_\2) b a_\1 \\
b \twoheadleftarrow \varphi := b \twoheadleftarrow (\varphi \o \A)
&= \sum_{i,j=1}^m \varphi(b_i a_j) a_i b b_j
\end{align*}
Note that the restriction of the action to~$A$ is just the right adjoint action of~$A^{\cop}$ on itself, and as for the left adjoint action considered before, the space of invariants is the center~$Z(A)$.

We also introduce the map
$$\bar{\Phi}: A^* \rightarrow A,~\varphi \mapsto 
(\A \twoheadleftarrow \varphi) = 
\sum_{i,j=1}^m \varphi(b_i a_j) a_i b_j = (\varphi \o \id)(R'R)$$
This map has similar properties as the map~$\Phi$:
If~$\bar{C}(A)$ denotes the subalgebra
$$\bar{C}(A) := \{\chi \in A^* \mid \chi(ab) = \chi(b\sa^{-2}(a)) \; 
\text{for all} \; a, b \in A\}$$
of~$A^*$, we have 
$$\bar{\Phi}(\chi \varphi) = \bar{\Phi}(\chi) \bar{\Phi}(\varphi)$$
for all $\chi \in \bar{C}(A)$ and all $\varphi \in A^*$. 
Furthermore, $\bar{\Phi}$ satisfies
$$\bar{\Phi}(\varphi) \twoheadleftarrow a 
= \varphi_\1(a_\2) \varphi_\3(\sa^{-1}(a_\1)) \; \bar{\Phi}(\varphi_\2)$$
i.e., it is an $A$-linear map from~$A^*$ with the right coadjoint action, built with the inverse antipode, to~$A$ with right adjoint action of~$A^{\cop}$. The space of invariants of the right coadjoint action is~$\bar{C}(A)$, whereas the space of invariant of the right adjoint action is the center~$Z(A)$, so that $\bar{\Phi}$ restricts to an algebra homomorphism from~$\bar{C}(A)$ to~$Z(A)$. These properties can be verified directly\endnote{\cite{DrinfAlmCocom}, Prop.~3.3, p.~327; \cite{Kerl}; Lem.~2, p.~362.} or derived from our construction of~$\bar{\Phi}$ in a way similar to the proof of the preceding proposition. It is also possible to derive them from the corresponding properties of~$\Phi$, as we have
$$\sa(\Phi(\varphi)) = \sa(\varphi \twoheadrightarrow \A) = 
\sa(\A) \twoheadleftarrow \sa^{-1*}(\varphi) = \bar{\Phi}(\sa^{-1*}(\varphi))$$
and similarly $\sa(\bar{\Phi}(\varphi))  = \Phi(\sa^{-1*}(\varphi))$.

The mappings~$\Phi$ and~$\bar{\Phi}$ are related in various ways to the Drinfel'd element~$\ua$. Besides the equations
$$\Phi(\varphi) = (\id \o \varphi)((\ua \o \ua)\da(\ua^{-1})) 
= (\id \o \varphi)(\da(\ua^{-1})(\ua \o \ua))$$
and  
$$\bar{\Phi}(\varphi) = (\varphi \o \id)((\ua \o \ua)\da(\ua^{-1})) 
= (\varphi \o \id)(\da(\ua^{-1})(\ua \o \ua))$$
which are direct consequences of the identity 
for~$\da(\ua)$ stated in Paragraph~\ref{QuasitriHopf}, we also have the following relation: The element $g:=\ua\sa(\ua^{-1})$ is a grouplike element.\endnote{\cite{DrinfAlmCocom}, Prop.~3.2, p.~327; \cite{M}, Thm.~10.1.13, p.~181.} If $\chi \in C(A)$, define~$\chi' \in A^*$ by
$\chi'(a) := \chi(ag^{-1})$. Then~$\chi' \in \bar{C}(A)$, and\endnote{\cite{DrinfAlmCocom}, Prop.~3.4, p.~328.}
$\bar{\Phi}(\chi')=\Phi(\chi)$.

As it turns out,\endnote{\cite{ResSem}, Thm.~2.9, p.~546; \cite{SchneiderFact}, Thm.~4.3, p.~1897.} the four conditions that~$\Phi$ is bijective, that~$\bar{\Phi}$ is bijective, that~$\Psi$ is bijective, and that~$\Psi'$ is bijective, are all equivalent. If these conditions are satisfied, the Hopf algebra~$A$ is called factorizable.\endnote{\cite{ResSem}, Def.~2.1, p.~543; see also \cite{EG1}, Lem.~1.1, p.~192; \cite{SchneiderFact}, p.~1892.}

\subsection[The coproduct of the evaluation form]{} \label{CoprodEvalForm}
If $A$ is the Drinfel'd double of a finite-dimensional Hopf algebra~$H$,
the mapping~$\bar{\Phi}$ takes a very simple form. To make this explicit, we decompose the dual of the double in the form $D(H)^* \cong H \o H^*$, where we use the isomorphism
$$H \o H^* \rightarrow D(H)^*,~h \o \varphi \mapsto 
(\varphi' \o h' \mapsto \varphi'(h) \varphi(h'))$$
From the form of the R-matrix of the Drinfel'd double described in Paragraph~\ref{DrinfDouble}, we see that
$$R'R = \sum_{i,j=1}^n (b_j^* \o b_i) \o (\eh \o b_j)(b_i^*\o \H)$$
so that
$\bar{\Phi}(h \o \varphi) = (h \o \varphi \o \id)(R'R)
= (\eh \o h)(\varphi \o \H)$.
A similar formula holds for~$\Phi$ under additional restrictions:\endnote{\cite{YYY2}, Par.~6.2, p.~44.}
\begin{lemma}
Suppose that $\psi = \sum_j h_j \o \varphi_j \in D(H)^*$.
\begin{enumerate}
\item 
If $\psi \in C(D(H))$, then 
$\Phi(\psi) = \sum_j \sh^{2*}(\varphi_j) \o h_j$.

\item 
If $\psi \in \bar{C}(D(H))$, then
$\Phi(\psi) = \sum_j \sh^{-2*}(\varphi_j) \o h_j$.
\end{enumerate}
\end{lemma}
\begin{pf}
For the first assertion, we have
\begin{align*}
\Phi(\psi) &= 
\sum_{i,j=1}^n (b_j^* \o b_i) \o \psi((\eh \o b_j)(b_i^*\o \H)) \\
&= 
\sum_{i,j=1}^n (b_j^* \o b_i) \o \psi((b_i^*\o \H)\sd^2(\eh \o b_j))
\\
&= 
\sum_{i,j=1}^n (b_j^* \o b_i) \o \psi(b_i^*\o \sh^2(b_j))
 = \sum_j \sh^{2*}(\varphi_j) \o h_j
\end{align*}
The proof of the second assertion is very similar.
\qed
\end{pf}

At the end of Paragraph~\ref{FactHopf}, we have expressed the mappings~$\Phi$ and~$\bar{\Phi}$ in terms of the Drinfel'd element.
In the case of a Drinfel'd double, we can replace the Drinfel'd element by the evaluation form to get very similar formulas for their inverses:
\begin{prop}
Suppose that $D=D(H)$ is the Drinfel'd double of a finite-dimensional Hopf algebra~$H$. Then we have for $x \in D(H)$ that
\begin{enumerate}
\item 
$\bar{\Phi}^{-1}(x) = 
((\id_{H^*} \o \sh^2)^* \circ \sd^{-1*})(e_\1^{-1} e) \; (e_\2^{-1} e)(x)$

\item
$\Phi^{-1}(x) = 
(\sd^* \circ(\sh^{2*} \o \id)^*)(e e_\2^{-1}) \; (e e_\1^{-1})(x)$
\end{enumerate}
\end{prop}
\begin{pf}
We can assume that $x = \varphi \o h$ for $\varphi \in H^*$ and~$h \in H$. Suppose that $b_1,\ldots,b_n$ is a basis of~$H$ with dual basis $b^*_1,\ldots,b^*_n$. Using the form of the R-matrix given in Paragraph~\ref{DrinfDouble}, we find
\begin{align*}
&\bar{\Phi}(((\id_{H^*} \o \sh^2)^* \circ \sd^{-1*})(e_\1^{-1} e)) \; (e_\2^{-1} e)(x) \\
&= \sum_{i,j=1}^n \sd^{-1*}(e_\1^{-1} e)(b^*_i \o \sh^2(b_j))
(\eh \o b_i)(b^*_j \o \H) 
e_\2^{-1}(\varphi_\2 \o h_\1) e(\varphi_\1 \o h_\2) \\
&= \sum_{i,j=1}^n e(b^*_{i\2} \o \sh^2(b_{j\1}))
e_\1^{-1}(\sd^{-1}(b^*_{i\1} \o \sh^2(b_{j\2})))
(\eh \o b_i)(b^*_j \o \H) \\
&\mspace{400mu}
e_\2^{-1}(\varphi_\2 \o h_\1) \varphi_\1(h_\2)
\end{align*}
where we have used the fact that~$e$ is invariant under the antipode observed in the proof of Lemma~\ref{IntDrinfDouble}. Using the definition of~$e$, this becomes
\begin{align*}
&\bar{\Phi}(((\id_{H^*} \o \sh^2)^* \circ \sd^{-1*})(e_\1^{-1} e)) \; (e_\2^{-1} e)(x) =\\
& \sum_{i,j=1}^n b^*_{i\2}(\sh^2(b_{j\1}))
e^{-1}(\sd^{-1}(b^*_{i\1} \o \sh^2(b_{j\2}))(\varphi_\2 \o h_\1))
\varphi_\1(h_\2)\\
&\mspace{450mu}
(\eh \o b_i)(b^*_j \o \H)=\\
& \sum_{i,j=1}^n b^*_{i\2}(\sh^2(b_{j\1}))
e^{-1}(\sh^*(b^*_{i\1}) \varphi_\2 \o h_\1 \sh(b_{j\2}))
\varphi_\1(h_\2) (\eh \o b_i)(b^*_j \o \H)
\end{align*}
where, in the last step, we have used that the antipode is antimultiplicative and that~$e^{-1}$ is cocommutative. Using the explicit form of~$e^{-1}$, this becomes
\begin{align*}
&\bar{\Phi}(((\id_{H^*} \o \sh^2)^* \circ \sd^{-1*})(e_\1^{-1} e)) \; (e_\2^{-1} e)(x) = \\
& \sum_{i,j=1}^n b^*_{i\2}(\sh^2(b_{j\1}))
\bigl(\sh^*(b^*_{i\1}) \varphi_\2 \bigr) (b_{j\2} \sh^{-1}(h_\1))
\varphi_\1(h_\2) (\eh \o b_i)(b^*_j \o \H) = \\
& \sum_{i,j=1}^n b^*_{i\2}(\sh^2(b_{j\1})) \varphi_\1(h_\3) \\
&\mspace{100mu}
\sh^*(b^*_{i\1})(b_{j\2} \sh^{-1}(h_\2))
\varphi_\2(b_{j\3} \sh^{-1}(h_\1)) (\eh \o b_i)(b^*_j \o \H) = \\
& \sum_{i,j=1}^n b^*_{i}(h_\2\sh(b_{j\2}) \sh^2(b_{j\1})) 
\varphi(h_\3 b_{j\3} \sh^{-1}(h_\1)) (\eh \o b_i)(b^*_j \o \H)
\end{align*}
Using the antipode equation, we can cancel two terms to get
\begin{align*}
&\bar{\Phi}(((\id_{H^*} \o \sh^2)^* \circ \sd^{-1*})(e_\1^{-1} e)) \; (e_\2^{-1} e)(x) = \\
& \sum_{i,j=1}^n b^*_{i}(h_\2) 
\varphi(h_\3 b_{j} \sh^{-1}(h_\1)) (\eh \o b_i)(b^*_j \o \H)= \\
&\sum_{j=1}^n \varphi(h_\3 b_{j} \sh^{-1}(h_\1)) (\eh \o h_\2)(b^*_j \o \H)= \\
&\varphi_\1(h_\3) \varphi_\3(\sh^{-1}(h_\1))  \;
(\eh \o h_\2)(\varphi_\2 \o \H)
= \varphi \o h = x
\end{align*}
This proves the first formula. For the second formula, recall from 
Paragraph~\ref{FactHopf} that
$\Phi^{-1} = \sd^* \circ \bar{\Phi}^{-1}  \circ \sd$, so that
we get from the first formula
\begin{align*}\Phi^{-1}(x) &= 
(\sd^* \circ \bar{\Phi}^{-1}  \circ \sd)(x) = 
(\sd^* \circ (\id \o \sh^2)^* \circ \sd^{-1*})(e_\1^{-1} e) \;  
(e_\2^{-1} e)(\sd(x)) \\
&= 
(\sd^* \circ (\sh^{2*} \o \id)^* \circ \sd^*)(e_\1^{-1} e) \;
(e_\2^{-1} e)(\sd(x)) \\
&= (\sd^* \circ(\sh^{2*} \o \id)^*)(e e_\2^{-1}) \; (e e_\1^{-1})(x)
\end{align*}
as asserted.
\qed
\end{pf}
Note that, in the case where the antipode of~$H$ is an involution, these formulas reduce to
$$\bar{\Phi}^{-1}(x) = 
\sd^*(e_\1^{-1} e) \; (e_\2^{-1} e)(x) \qquad \qquad
\Phi^{-1}(x) = \sd^*(e e_\2^{-1}) \; (e e_\1^{-1})(x)$$
which should be compared with the ones given in Paragraph~\ref{FactHopf} using the Drinfel'd element.

\subsection[The double and the tensor product]{} \label{DoubleTens}
If $f: A \rightarrow B$ is a morphism of quasitriangular Hopf algebras, then it follows directly from the formulas for~$\Phi$ and~$\bar{\Phi}$
given in Paragraph~\ref{FactHopf} that the diagrams
$$
\square<1`1`1`-1;1000`700>[B^*`A^*`B`A;f^*`\Phi_B`\Phi_A`f]
\mspace{100mu}
\square<1`1`1`-1;1000`700>[B^*`A^*`B`A;f^*`\bar{\Phi}_B`\bar{\Phi}_A`f]$$
commute, where the index indicates to which Hopf algebra the mapping belongs. Applying this to~$\Psi$, we see that the diagrams
$$
\square<1`1`1`-1;1000`700>[(A \o A)_{F}^*`D(A)^*`(A \o A)_{F}`D(A);\Psi^*`\Phi`\Phi`\Psi]
\mspace{100mu}
\square<1`1`1`-1;1000`700>[(A \o A)_{F}^*`D(A)^*`(A \o A)_{F}`D(A);\Psi^*`\bar{\Phi}`\bar{\Phi}`\Psi]$$
commute. Now we have by definition that $(A \o A)_{F} = A \o A$ as algebras. But these two Hopf algebras have even more things in common:
\begin{lemma}
We have $C((A \o A)_{F}) = C(A \o A)$ and 
$\bar{C}((A \o A)_{F}) = \bar{C}(A \o A)$ as algebras.
\end{lemma}
\begin{pf}
We have already seen in the proof of Lemma~\ref{DoubleQuasitri} that the Drinfel'd elements of $(A \o A)_{F}$ and $A \o A$ coincide. Therefore, the squares of the antipodes coincide, too. This implies the asserted equalities as vector spaces. That the products also agree follows from the fact that $(S^2 \o S^2)(F) = F$ in our case.
\qed
\end{pf}

From this lemma, we can extract the following information about the restrictions of~$\Phi$ and~$\bar{\Phi}$:
\begin{prop}
The diagrams
$$
\square<1`1`1`-1;1400`700>[C(A) \o C(A)`C(D(A))`A \o A`Z(D(A));\Psi^*`\Phi \o (\sa \circ\Phi)`\Phi`\Psi]$$
and
$$\square<1`1`1`-1;1400`700>[\bar{C}(A) \o \bar{C}(A)`\bar{C}(D(A))`A \o A`Z(D(A));\Psi^*`\bar{\Phi} \o (\sa \circ \bar{\Phi})`\bar{\Phi}`\Psi]$$
commute.
\end{prop}
\begin{pf}
From the formula $R_F = F_t R_{A \o A} F^{-1}$ for the twisted R-matrix, we get immediately that $R'_F R_F = F R'_{A \o A} R_{A \o A} F^{-1}$. Now we have
\begin{align*}
R'_{A \o A} R_{A \o A}  = \sum_{i,j,k,l=1}^m 
b_k a_i \o  \sa(a_l) b_j  \o a_k b_i  \o b_l \sa(a_j)
\end{align*}
and therefore
\begin{align*}
R'_F R_F &= (\sum_{p=1}^m \A \o \sa(a_p) \o b_p \o \A)
R'_{A \o A} R_{A \o A}
(\sum_{q=1}^m \A \o a_q \o b_q \o \A)\\
&= 
\sum_{i,j,k,l,p,q=1}^m 
b_k a_i \o \sa(a_p) \sa(a_l) b_j a_q \o b_p a_k b_i b_q \o b_l \sa(a_j)
\end{align*}
For two elements $\varphi, \psi \in C(A)$, we therefore find
\begin{align*}
\Phi(\varphi \o \psi) = \sum_{i,j,k,l,p,q=1}^m 
b_k a_i \o \sa(a_p) \sa(a_l) b_j a_q  \varphi(b_p a_k b_i b_q)  
\psi(b_l \sa(a_j))
\end{align*}
Now note that
\begin{align*}
\sum_{j,l=1}^m \sa(a_l) b_j \psi(b_l \sa(a_j))
&= \sum_{j,l=1}^m \sa(a_l) b_j \psi(\sa(a_j) \sa^2(b_l))
= \sum_{j,l=1}^m a_l b_j \psi(\sa(a_j) \sa(b_l)) \\
&= \sum_{j,l=1}^m \sa^{-1}(\sa(b_j) \sa(a_l)) \psi(\sa(a_j) \sa(b_l)) \\
&= \sum_{j,l=1}^m \sa^{-1}(b_j a_l) \psi(a_j b_l)
= \sa^{-1}(\Phi(\psi))
\end{align*}
which is a central element by Proposition~\ref{FactHopf}. We can therefore rewrite the above expression in the form
\begin{align*}
\Phi(\varphi \o \psi) &= \sum_{i,k,p,q=1}^m 
b_k a_i \o \sa(a_p) a_q \sa^{-1}(\Phi(\psi)) \varphi(b_p a_k b_i b_q) 
\end{align*}
But in this expression, we can cancel the summation over~$p$ and~$q$, as we have
\begin{align*}
&\sum_{p,q=1}^m \sa(a_p) a_q \varphi(b_p a_k b_i b_q) 
=\sum_{p,q=1}^m \sa(a_p) a_q \varphi(a_k b_i b_q \sa^2(b_p))  \\
&=\sum_{p,q=1}^m a_p a_q \varphi(a_k b_i b_q \sa(b_p)) 
=\sum_{p,q=1}^m a_p a_q \varphi(a_k b_i \sa(b_p \sa^{-1}(b_q)))
= \A \varphi(a_k b_i)
\end{align*}
After this cancellation, we get 
\begin{align*}
\Phi(\varphi \o \psi) &= 
\sum_{i,k=1}^m b_k a_i \o \sa^{-1}(\Phi(\psi)) \varphi(a_k b_i) \\
&= \Phi(\varphi) \o \sa^{-1}(\Phi(\psi)) 
= \Phi(\varphi) \o \sa(\Phi(\psi))
\end{align*}
where the last step follows from the fact that 
$\sa^2(\Phi(\psi)) = \ua \Phi(\psi) \ua^{-1}= \Phi(\psi)$ since
$\Phi(\psi)$ is central.
This shows the commutativity of the first diagram. The commutativity of the second diagram follows from similar computations.
\qed
\end{pf}
Note that, in contrast to the factorizable case, it is not true in general that $\Psi$ maps the center of~$D(A)$
to the center of~$A \o A$, as the example of a group ring with an R-matrix equal to the unit shows.

In the whole discussion, it is possible to replace the original R-matrix by~$R'^{-1}$. We have already pointed out above that this interchanges~$\Psi$ and~$\Psi'$ as well as~$F$ and~$F'$. The analogue of~$\Phi$ is the map that assigns to every~$\varphi \in A^*$ the element $(\id \o \varphi)(R^{-1} R'^{-1})$. If $\varphi \in C(A)$, this element can be expressed in terms of the original map~$\Phi$ as follows:
\begin{align*}
(\id \o \varphi)(R^{-1} R'^{-1}) &= 
\sum_{i,j=1}^m \sa(a_i) b_j \varphi(b_i \sa(a_j)) = 
\sum_{i,j=1}^m \sa(a_i) b_j \varphi(\sa(a_j) \sa^2(b_i)) \\
&= \sa^{-1}(\sum_{i,j=1}^m \sa(b_j) \sa^2(a_i) 
\varphi(\sa(a_j) \sa^2(b_i))) \\
&= \sa^{-1}(\sum_{i,j=1}^m b_j a_i \varphi(a_j b_i))
= \sa^{-1}(\Phi(\varphi))= \sa(\Phi(\varphi))
\end{align*}
where the last step uses the argument from the end of the proof of the preceding proposition. Something similar holds for the analogue of~$\bar{\Phi}$: If
$\varphi \in \bar{C}(A)$, we have 
$$(\varphi \o \id)(R^{-1} R'^{-1}) = \sa(\bar{\Phi}(\varphi))$$
as can be seen by a similar computation. Therefore, the corresponding commutative diagrams for~$\Psi'$ are
$$
\square<1`1`1`-1;1400`700>[C(A) \o C(A)`C(D(A))`A \o A`Z(D(A));\Psi'^*`(\sa \circ\Phi) \o \Phi`\Phi`\Psi']$$
and
$$\square<1`1`1`-1;1400`700>[\bar{C}(A) \o \bar{C}(A)`\bar{C}(D(A))`A \o A`Z(D(A));\Psi'^*`(\sa \circ \bar{\Phi}) \o \bar{\Phi}`\bar{\Phi}`\Psi']$$

\subsection[Integrals of factorizable Hopf algebras]{} \label{IntFact}
Let us now assume in addition that our quasitriangular Hopf algebra~$A$ is also factorizable. It is then finite-dimensional and unimodular,\endnote{\cite{RadfKnotInv}, Prop.~3, p.~224; \cite{SchneiderFact}, Rem.~4.4, p.~1898.} which implies\endnote{\cite{R4}, Thm.~3, p.~594.} that a right integral~$\rha \in A^*$ must be contained in~$\bar{C}(A)$. The image of this integral under~$\bar{\Phi}$ is again an integral:
\begin{lemma}
$\bar{\Phi}(\rha)$ is a two-sided integral of~$A$. 
\end{lemma}
\begin{pf}
It follows from the discussion in Paragraph~\ref{FactHopf} that 
$$\bar{\Phi}(\rha) \bar{\Phi}(\varphi) = \bar{\Phi}(\rha \varphi) = \varphi(\A) \bar{\Phi}(\rha ) = 
\ea(\bar{\Phi}(\varphi)) \bar{\Phi}(\rha) $$
for all $\varphi \in A^*$. Because $A$ is factorizable, every $a \in A$ can be written in the form $a = \bar{\Phi}(\varphi)$, showing that $\bar{\Phi}(\rha)$ is a right integral. It is also a left integral since~$A$ is unimodular, or alternatively because $\bar{\Phi}(\rha)$
is central, as we saw in Paragraph~\ref{FactHopf}.
\qed
\end{pf}
From this lemma, it follows in particular that
$\rha(\bar{\Phi}(\rha)) = (\rha \o \rha)(R'R) \neq 0$ \linebreak
if $\rha \neq 0$, because nonzero integrals do not vanish on nonzero integrals.\endnote{\cite{LR2}, Prop.~1, p.~269.} It also shows that a factorizable Hopf algebra is semisimple if and only if it is cosemisimple, because $\rha(\A)=\ea(\bar{\Phi}(\rha))$, and, by Maschke's theorem,
$A$~is cosemisimple if and only if $\rha(\A) \neq 0$, and semisimple if and only if $\ea(\bar{\Phi}(\rha)) \neq 0$.\endnote{\cite{Tak}, Prop.~4.4, p.~639.} Finally, it also shows that, for a left integral~$\lma \in A^*$,
the element $\Phi(\lma)$ is a two-sided integral of~$A$, since we get from Paragraph~\ref{FactHopf} that
$\sa(\bar{\Phi}(\rha)) = \Phi(\sa^{-1*}(\rha))$.

We now discuss how the integrals behave under~$\Psi$. It is obvious that $\rha \o \rha$ is a right integral on the tensor product 
Hopf algebra~$A \o A$. Therefore, a general result on the behavior of integrals under twisting,\endnote{\cite{AEGN}, Thm.~3.4, p.~488.}
together with the discussion in Lemma~\ref{DoubleQuasitri}, yields 
that it is also a right integral on~$(A \o A)_F$. As~$\Psi$ is a Hopf algebra isomorphism, $\Psi^*(\rha \o \rha)$ must be a right integral on~$D(A)^*$. To make this more precise, we decompose the dual of the double in the form $D(A)^* \cong A \o A^*$ as described in Paragraph~\ref{CoprodEvalForm}.
We then get the following formulas for the integrals:
\begin{prop}
Suppose that $\lma \in A^*$ is a left integral, that $\rha \in A^*$ is a right integral, and that $\Ga \in A$ is a two-sided integral. Then we have
\begin{enumerate}
\item $\Psi(\lma \o \Ga) = \Phi(\lma) \o \Ga$
\item $\Psi^*(\lma \o \lma) = \Phi(\lma) \o \lma$
\item $\Psi^*(\rha \o \rha) = \bar{\Phi}(\rha) \o \rha$
\end{enumerate}
\end{prop}
\begin{pf}
\begin{list}{(\arabic{num})}{\usecounter{num} \leftmargin0cm \itemindent5pt}
\item
To establish the first assertion, we compute as follows:
\begin{align*}
&\Psi(\lma \o \Ga) = \pi(\lma_\2 \o \Ga_\1) \o \pi'(\lma_\1 \o \Ga_\2) \\
&= \sum_{i,j=1}^m 
\lma_\2(a_i) b_i \Ga_\1 \o \lma_\1(b_j) \sa(a_j) \Ga_\2 
= \sum_{i,j=1}^m \lma(b_j a_i) b_i \Ga_\1 \o \sa(a_j) \Ga_\2 \\
&= \sum_{i,j=1}^m \lma(b_j a_i) b_i \sa^2(a_j) \Ga_\1 \o \Ga_\2
= \sum_{i,j=1}^m \lma(a_i \sa^2(b_j)) b_i \sa^2(a_j) \Ga_\1 \o \Ga_\2\\
&= \sum_{i,j=1}^m \lma(a_i b_j) b_i a_j \Ga_\1 \o \Ga_\2
= \Phi(\lma) \Ga_\1 \o \Ga_\2
= \Phi(\lma) \o \Ga
\end{align*}
where unimodularity is used in particular for the fifth equation.

\item
For the second assertion, we get from Proposition~\ref{DoubleTens} that
\begin{align*}
\Phi(\Psi^*(\lma \o \lma)) = \Psi^{-1}(\Phi(\lma) \o \sa(\Phi(\lma)))
= \Psi^{-1}(\Phi(\lma) \o \Phi(\lma))
\end{align*}
since $\lma \in C(A)$. But this gives
$\Phi(\Psi^*(\lma \o \lma)) = \lma \o \Phi(\lma)$
by the assertion just established. Since~$C(D(A))$ and~$Z(D(A))$
have the same dimension,
it now follows from Lemma~\ref{CoprodEvalForm} that\endnote{\cite{SchneiderFact}, Lem.~2.2, p.~1893; \cite{R4}, Prop.~3, p.~590.}
$$\Psi^*(\lma \o \lma) = \Phi^{-1}(\lma \o \Phi(\lma)) = 
\Phi(\lma) \o \sa^{-2*}(\lma) = 
\Phi(\lma) \o \lma$$

\item
For the third assertion, we substitute~$\bar{\Phi}(\rha)$ for~$\Ga$
and~$\sa^*(\rha)$ for~$\lma$ into the first assertion to get,
using another result from Paragraph~\ref{FactHopf}, that
\begin{align*}
\Psi(\sa^*(\rha) \o \bar{\Phi}(\rha)) = 
\Phi(\sa^*(\rha)) \o \bar{\Phi}(\rha) = 
\sa^{-1}(\bar{\Phi}(\rha)) \o \bar{\Phi}(\rha) = 
\bar{\Phi}(\rha) \o \bar{\Phi}(\rha)
\end{align*}
Using the second part of Proposition~\ref{DoubleTens} on the right-hand side, this becomes
$\Psi(\sa^*(\rha) \o \bar{\Phi}(\rha)) = 
(\Psi \circ \bar{\Phi} \circ \Psi^*)(\rha \o \rha)$,
so that
$$\sa^*(\rha) \o \bar{\Phi}(\rha) = 
(\bar{\Phi} \circ \Psi^*)(\rha \o \rha)$$
Since the left-hand side is a two-sided integral of~$D(A)$, it is invariant under the antipode, so that we can rewrite this equation as
\begin{align*}
(\bar{\Phi} \circ \Psi^*)(\rha \o \rha) &=
\sd(\sa^*(\rha) \o \bar{\Phi}(\rha)) =
\sd(\eh \o \bar{\Phi}(\rha))\sd(\sa^*(\rha) \o \A) \\
&= (\eh \o \bar{\Phi}(\rha))(\rha \o \A)
= \bar{\Phi}(\bar{\Phi}(\rha) \o \rha)
\end{align*}
by the discussion in Paragraph~\ref{CoprodEvalForm}.
Now cancelling~$\bar{\Phi}$ gives the assertion.
\qed
\end{list}
\end{pf}

\newpage
\section{The action of the modular group} \label{Sec:ActModGroup}
\subsection[The role of the integral]{} \label{RoleInt}
In our factorizable Hopf algebra~$A$, we now fix a nonzero right integral~$\rha$, and introduce the map
$$\iota: A \rightarrow A^*,~a \mapsto \iota(a) := \rha_\1(a) \; \rha_\2$$
The fact that $\rha$ is a Frobenius homomorphism\endnote{\cite{M}, Thm.~2.1.3, p.~18.} implies that~$\iota$ is bijective. The fact 
that~$\rha \in \bar{C}(A)$ implies that~$\iota$ is an
$A$-linear map from~$A$ with the right adjoint action of~$A$
to~$A^*$ with the right coadjoint action, built with the inverse antipode, which we have considered in Paragraph~\ref{FactHopf}. In particular, $\iota$ induces an isomorphism between the spaces of invariants of these actions; in other words, it restricts to a bijection between the center~$Z(A)$ and the algebra~$\bar{C}(A)$.\endnote{\cite{SchneiderFact}, Lem.~2.2, p.~1893.}
Following\endnote{See also \cite{Ly2}, Def.~6.2, p.~320; \cite{LyuMaj}, Thm.~1.1, p.~507.}  \cite{Kerl}, Eq.~(2.55), p.~369, we use~$\iota$ to introduce the maps $\v \in \End(A)$ and~$\v_* \in \End(A^*)$ as
$$\v := \sa \circ \bar{\Phi} \circ \iota \qquad \qquad
\v_* := \sa^{-1*} \circ \iota \circ \Phi$$
It is important to distinguish~$\v_*$ from the transpose~$\v^*$ of~$\v$. 
Using the form of the R-matrix given in Paragraph~\ref{QuasitriHopf}, we have explicitly
$$\v(a) = \sum_{i,j=1}^m \rha(a b_i a_j) \sa(a_i b_j)$$
The relationship of these maps is clarified in the following proposition, which also lists some of their basic properties:
\begin{prop}
$\v$ is an $A$-linear map from the right adjoint representation of~$A$ to itself. In particular, $\v$ preserves the center~$Z(A)$ of~$A$. Furthermore, the diagrams
$$
\square<1`1`1`1;1000`700>[A^*`A^*`A`A;\v_*`\Phi`\Phi`\v]
\qquad \qquad
\square<1`-1`-1`1;1000`700>[A^*`A^*`A`A;\v^*`\iota`\iota`\v]$$
are commutative, and we have 
$\v^* \circ \sa^* = \sa^* \circ \v_*$.
\end{prop}
\begin{pf}
The composition of $A$-linear maps is $A$-linear. From the linearity properties of~$\bar{\Phi}$ and~$\iota$ discussed so far, we get that 
$\bar{\Phi} \circ \iota$ is $A$-linear from the right adjoint representation of~$A$ to the right adjoint representation of~$A^{\cop}$; i.e., satisfies
$$(\bar{\Phi} \circ \iota)(\sa(a_\1) b a_\2) = 
\sa^{-1}(a_\2)(\bar{\Phi} \circ \iota)(b) a_\1$$
Applying the antipode to this equation gives
$\v(\sa(a_\1) b a_\2) = \sa(a_\1)\v(b) a_\2$,
which is the first assertion. The preservation of the center is a direct consequence. The commutativity of the first diagram follows from the equation
$$\Phi \circ \v_* = \Phi \circ \sa^{-1*} \circ \iota \circ \Phi
= \sa \circ \bar{\Phi} \circ \iota \circ \Phi = \v \circ \Phi$$
where we have used the fact that 
$\sa \circ \bar{\Phi} = \Phi \circ \sa^{-1*}$ established in Paragraph~\ref{FactHopf}.
To establish the commutativity of the second diagram, we first derive the formula~$\v^* = \iota \circ \sa \circ \bar{\Phi}$. For~$\varphi \in A^*$
and $a \in A$, we have
\begin{align*}
\v^*(\varphi)(a) &= \varphi(\v(a)) = \varphi(\sa(\bar{\Phi}(\iota(a))))
= (\iota(a) \o \sa^*(\varphi))(R'R) = \iota(a)(\Phi(\sa^*(\varphi))) \\
&= \iota(a)(\sa^{-1}(\bar{\Phi}(\varphi))) 
= \rha(a\sa^{-1}(\bar{\Phi}(\varphi))) = \rha(\sa(\bar{\Phi}(\varphi))a)
= \iota(\sa(\bar{\Phi}(\varphi)))(a)
\end{align*}
where we have used for the second last equality that $\rha \in \bar{C}(A)$.
From this, we immediately get the commutativity of the second diagram, as we now have  
$\v^* \circ \iota = \iota \circ \sa \circ \bar{\Phi} \circ \iota 
= \iota \circ \v$. Furthermore, using the results from Paragraph~\ref{FactHopf}, we can rewrite the formula for~$\v^*$ above in the form $\v^* = \iota \circ \Phi \circ \sa^{-1*}$, which yields
$\v^* \circ \sa^* = \iota \circ \Phi = \sa^* \circ \v_*$. \qed
\end{pf}

It may be noted that the commutativity of the second diagram is equivalent to the condition $\rha(\v(a)b) = \rha(a\v(b))$, which is an adjunction property of~$\v$ with respect to the associative bilinear form determined by the Frobenius homomorphism~$\rha$.

\subsection[The inverse of~$\v$]{} \label{InvSigma}
As we have explained in Paragraph~\ref{QuasitriHopf}, $R'^{-1}$ is always an alternative choice for the R-matrix of a quasitriangular Hopf algebra.
This raises the question how~$\v$ is modified if one replaces~$R$ by~$R'^{-1}$. The answer to this question is that, up to a scalar multiple,
$\v$ turns into its inverse:
\begin{prop}
$\displaystyle
\v^{-1}(a) = \frac{1}{(\rha \o \rha)(R'R)} 
\sum_{i,j=1}^m \rha(a a_i b_j) \sa^2(a_j) b_i$
\end{prop}
\begin{pf}
Because $A$ is finite-dimensional, it suffices to prove that
$$\sum_{i,j=1}^m \rha(a a_i b_j) \v(\sa^2(a_j) b_i)
= (\rha \o \rha)(R'R) \; a$$
To see this, we use the identity
$\rha(a b_\1) \sa(b_\2) = \rha(a_\1 b) a_\2$ to compute 
\begin{align*}
&\sum_{i,j=1}^m \rha(a a_i b_j) \v(\sa^2(a_j) b_i)
= \sum_{i,j,k,l=1}^m \rha(a a_i b_j) \; \rha(\sa^2(a_j) b_i b_k a_l) \;
\sa(a_k b_l)  \\
&= \sum_{i,j,k,l=1}^m \rha(a a_i b_j) \; 
\rha(b_i b_k a_l a_j) \; \sa(a_k b_l)  
= \sum_{i,j=1}^m  \rha(a a_i{}_\1 b_j{}_\1)\; 
\rha(b_i a_j) \sa(a_i{}_\2 b_j{}_\2)   \\
&= \sum_{i,j=1}^m \rha(a_\1 a_i b_j) \; \rha(b_i a_j) \; a_\2 
= \rha(a_\1 \bar{\Phi}(\rha))  \; a_\2 = \rha(\bar{\Phi}(\rha)) \; a
\end{align*}
where the last step follows from Lemma~\ref{IntFact}.
\qed
\end{pf}

This formula has an interesting consequence for the restriction of~$\v$ to the center:
\begin{corollary}
For all $a \in Z(A)$, we have $\v^2(a) = (\rha \o \rha)(R'R) \; \sa(a)$.
\end{corollary}
\begin{pf}
In this case, we get from the formula in the preceding proposition that
\begin{align*}
(\rha \o \rha)(R'R) \; \v^{-1}(a) &= \sum_{i,j=1}^m \rha(a a_i b_j) \sa^2(a_j) b_i
= \sum_{i,j=1}^m \rha(a a_i \sa^{-2}(b_j)) a_j b_i \\
&= \sum_{i,j=1}^m \rha(b_j a a_i) a_j b_i 
= \sum_{i,j=1}^m \rha(a b_j a_i) a_j b_i = \sa^{-1}(\v(a))
\end{align*}
which becomes the assertion if we insert~$\v(a)$ for~$a$.
\qed
\end{pf}
By rescaling~$\rha$ if necessary, we can achieve that $(\rha \o \rha)(R'R) = 1$, as we saw in Paragraph~\ref{IntFact} that this expression is nonzero, and in our algebraically closed field every element has a square root. In this case, the formula in the preceding corollary asserts that
$\v^2(a) = \sa(a)$ for all $a \in Z(A)$.

\subsection[Ribbon elements]{} \label{RibEl}
Recall\endnote{\cite{Tur}, Sec.~XI.3.1, p.~500; note the difference to \cite{Kas}, Def.~XIV.6.1, p.~361.} that a ribbon element is a nonzero central element $v \in A$ that satisfies
$$\da(v) = (R'R) (v \o v) \quad \text{and} \quad \sa(v) = v$$
It follows\endnote{\cite{Kas}, Cor.~XIV.6.3, p.~362.} that~$v$ is an invertible element that satisfies~$\ea(v)=1$ as well as 
$v^{-2} = \ua \sa(\ua)$. We use it to define the endomorphism
$$\t: A \rightarrow A,~a \mapsto v a$$
which is just multiplication by the central element~$v$. The fact that~$v$ is central directly yields the equation
$\rha(\t(a)b) = \rha(a\t(b))$, which can, as for~$\v$ in Paragraph~\ref{RoleInt}, be expressed by saying that the diagram
$$\square<1`-1`-1`1;1000`700>[A^*`A^*`A`A;\t^*`\iota`\iota`\t]$$
is commutative.

The decisive relation between~$\v$ and~$\t$ is the following:\endnote{\cite{Kerl}, Prop.~13, p.~372; \cite{Ly2}, Thm.~6.5, p.~321; \cite{LyuMaj}, Eq.~(1), p.~507.}
\begin{prop}
$\v \circ \t \circ \v
= \rha(v) \; \t^{-1} \circ \v \circ \t^{-1}$
\end{prop}
\begin{pf}
From Proposition~\ref{InvSigma}, we have
\begin{align*}
\rha(\bar{\Phi}(\rha)) \;& (\v^{-1} \circ \t \circ \v)(a) 
= \rha(\bar{\Phi}(\rha)) \;
\sum_{i,j=1}^m \rha(a b_i a_j) \v^{-1}(v \sa(a_i b_j))
\\ &= 
\sum_{i,j,k,l=1}^m \rha(a b_i a_j) \rha(v \sa(a_i b_j) a_k b_l)  \sa^2(a_l) b_k \\ &= 
\sum_{i,j,k,l=1}^m 
\rha(a b_i a_j) \rha(\sa^2(b_l) v \sa(b_j) \sa(a_i) a_k) \sa^2(a_l) b_k \\ &= \sum_{i,j,k,l=1}^m 
\rha(a \sa^{-1}(b_i) \sa^{-1}(a_j)) \rha(v b_l b_j a_i a_k) a_l b_k \\ 
&= \sum_{i,j=1}^m 
\rha(a \sa^{-1}(b_{i\2}) \sa^{-1}(a_{j\2})) \rha(v b_j a_i) a_{j\1} b_{i\1}
\end{align*}
Using that $(\A \o v^{-1}) \da(v) = \sum_{i,j=1}^m v b_j a_i \o a_j b_i$,
we can write this in the form
\begin{align*}
\rha(\bar{\Phi}(\rha)) \; (\v^{-1} \circ \t \circ \v)(a) 
&= \rha(a \sa^{-1}(v^{-1}_{\2} v_{\3})) \rha(v_\1) v^{-1}_{\1} v_{\2} \\
&= \rha(a \sa^{-1}(v^{-1}_{\2})) \rha(v) v^{-1}_{\1} 
\end{align*}
On the other hand, we have that
$$\da(v^{-1}) = (v^{-1} \o v^{-1}) (R'R)^{-1} 
= \sum_{i,j=1}^m v^{-1}  a_i b_j\o  v^{-1} \sa^{-1}(b_i) \sa(a_j) $$
so that Proposition~\ref{InvSigma} also implies that
\begin{align*}
\rha(\bar{\Phi}(\rha)) \; (\t^{-1} \! \circ \v^{-1} \!\circ \t^{-1})(a) 
= \sum_{i,j=1}^m \rha(v^{-1} a a_i b_j) v^{-1} \sa^2(a_j) b_i
= \rha(a v_\1^{-1}) \; \sa(v_\2^{-1})  
\end{align*}
Comparing both expressions and using $\sa(v^{-1}) = v^{-1}$, we get that
$$(\v^{-1} \circ \t \circ \v)(a) 
= \rha(v) \; (\t^{-1} \circ \v^{-1} \circ \t^{-1})(a)$$
which is equivalent to the assertion.
\qed
\end{pf}

The restrictions of~$\v$ and~$\t$ to the center of~$A$
induce of course also automorphisms of the corresponding projective space~$P(Z(A))$ of one-dimensional subspaces of~$Z(A)$, which are even independent of the choice of the integral~$\rha$. It is a consequence of the results above that these automorphisms yield a projective representation of the modular group:
\begin{corollary}
There is a unique homomorphism from~$\SL(2,\Z)$ to~$\PGL(Z(A))$ 
that maps~$\gv$ to the equivalence class of~$\v$ and~$\gt$ to the equivalence class of~$\t$.
\end{corollary}
\begin{pf}
The homomorphism is unique because~$\gv$ and~$\gt$ generate the modular group, as discussed in Paragraph~\ref{GenRel}. For the existence question, recall the defining relations $\gv^4 = 1$ and $(\gt\gv)^3 = \gv^2$. Because the square of the antipode is given by conjugation with the Drinfel'd element,\endnote{\cite{M}, Prop.~10.1.4, p.~179; \cite{Tur}, Sec.~XI.2.2, Eq.~(2.2.c), p.~498; \cite{Kas}, Prop.~VIII.4.1, p.~180.} it restricts to the identity on the center, and therefore Corollary~\ref{InvSigma} implies the first relation needed. The second defining relation $\gt\gv\gt\gv\gt\gv = \gv^2$ can also be written in the form
$\gv\gt\gv = \gt^{-1} \gv \gt^{-1}$, and therefore it follows from the preceding proposition that this relation is satisfied, too.
\qed
\end{pf}

The proof shows that if $(\rha \o \rha)(R'R) = 1$ and $\rha(v)=1$, we even get a linear representation $\SL(2,\Z) \rightarrow \GL(Z(A))$ by assigning~$\v$ to~$\gv$ and~$\t$ to~$\gt$. However, this happens if and only if~$\rha(v) = \rha(v^{-1})$, as we see from the following lemma:\endnote{\cite{RadfKnotInv}, Cor.~2, p.~226.}
\begin{lemma}
$(\rha \o \rha)(R'R) = \rha(v) \rha(v^{-1})$
\end{lemma}
\begin{pf}
We have seen in Lemma~\ref{IntFact} 
that~$\bar{\Phi}(\rha) = \rha(v^{-1} v_\1) v^{-1} v_\2 \in A$ 
is a two-sided integral. We therefore have
$$\rha(v^{-1} v_\1) v_\2 = v \bar{\Phi}(\rha) 
= \ea(v) \bar{\Phi}(\rha) = \bar{\Phi}(\rha)$$
Now there is a grouplike element~$g \in A$, called the right modular element, that satisfies $a_\1 \rha(a_\2) = g \rha(a)$ for all~$a \in A$, and furthermore $\rha(ag) = \rha(\sa^{-1}(a))$.\endnote{\cite{R4}, Prop.~3, p.~590.}
The preceding computation therefore yields that
\begin{align*}
(\rha \o \rha)(R'R) &= \rha(\bar{\Phi}(\rha)) 
= \rha(v^{-1} v_\1) \rha(v_\2) \\
&= \rha(v^{-1} g) \rha(v) 
= \rha(\sa^{-1}(v^{-1})) \rha(v)= \rha(v^{-1}) \rha(v)
\end{align*}
as asserted.
\qed
\end{pf}
It may be noted that we have discussed after Lemma~\ref{IntFact} that this quantity is nonzero if~$\rha$ is nonzero. Furthermore, since~$A$ is unimodular, the right modular element~$g$ that appears in the preceding proof is exactly the grouplike element, also denoted by~$g$, that appeared in the discussion at the end of Paragraph~\ref{FactHopf}.\endnote{\cite{M}, Prop.~10.1.14, p.~183.}
It should also be noted that we do not claim that it is impossible to modify the representation so that it becomes linear.\endnote{\cite{BakKir}, Rem.~3.1.9, p.~52; \cite{Beyl}, Sec.~4, p.~34.}

\subsection[Integrals, ribbon elements, and the double]{} \label{IntRibDouble}
We have seen in Proposition~\ref{IntFact} that
$\rhd := \Psi^*(\rha \o \rha) = \bar{\Phi}(\rha) \o \rha$
is a right integral in~$D(A)^*$. As in the case of~$A$ itself, we therefore get an isomorphism
$$\iota_D: D(A) \rightarrow D(A)^*,~x \mapsto \iota_D(x) := 
\rha_{D\1}(x) \; \rha_{D\2}$$
Since $\rhd$ is defined as the image of~$\rha \o \rha$ under the Hopf algebra homomorphism~$\Psi^*$, it is obvious that the diagram
$$\square<-1`-1`-1`1;1400`700>[D(A)^*`A^* \o A^*`D(A)`A \o A;
\Psi^*`\iota_D`\iota \o \iota`\Psi]$$
commutes, as we have already pointed out in Paragraph~\ref{IntFact}
that~$\rho \o \rho$ is also a right integral in~$(A \o A)_F^*$.
Combining this with Proposition~\ref{DoubleTens}, we get the following commutative diagram:
$$\square<1`-1`-1`1;1400`700>[Z(D(A))`Z(A) \o Z(A)`Z(D(A))`Z(A) \o Z(A);\Psi`\bar{\Phi} \circ \iota_D`(\bar{\Phi} \circ \iota)\o 
(\sa \circ \bar{\Phi} \circ \iota)`\Psi]$$
From the general formula for the antipode of a twist mentioned in the proof of Lemma~\ref{DoubleQuasitri}, it is immediate that the antipode of
$(A \o A)_F$ coincides with the antipode of~$A \o A$ on the center. This implies that the following diagram is also commutative:
$$\square<1`-1`-1`1;1400`700>[Z(D(A))`Z(A) \o Z(A)`Z(D(A))`Z(A) \o Z(A);\Psi`\v`\v \o (\sa \circ \v)`\Psi]$$

The ribbon element can be treated in a similar way. It is immediate from the definition that a ribbon element $v \in A$ satisfies
$$\da(v^{-1}) = (R^{-1}R'^{-1}) (v^{-1} \o v^{-1})$$
which means that~$v^{-1}$ is a ribbon element for~$A$ endowed with the alternative R-matrix~$R'^{-1}$.
This implies that $v \o v^{-1}$ is a ribbon element for~$A \o A$, endowed with the R-matrix~$R_{A \o A}$ considered in Paragraph~\ref{DoubleQuasitri}. It is not difficult to see that a ribbon element stays a ribbon element if the coproduct of the Hopf algebra is twisted, and therefore~$v \o v^{-1}$ is also a ribbon element 
for~$(A \o A)_F$. This enables us to define a ribbon element~$v_D$ of the Drinfel'd double~$D(A)$ by setting~$v_D:=\Psi^{-1}(v \o v^{-1})$. 
So, if we define~$\t \in \End(D(A))$ to be the multiplication by~$v_D$, as in Paragraph~\ref{RibEl}, it is obvious that the following diagram is commutative:
$$\square<1`-1`-1`1;1400`700>[D(A)`A \o A`D(A)`A \o A;\Psi`\t`\t \o \t^{-1}`\Psi]$$
Note that it follows from Lemma~\ref{DoubleQuasitri} that in the case that the ribbon element is the inverse Drinfel'd element, so that~$v=\ua^{-1}$, the
arising ribbon element of~$D(A)$ is also the inverse Drinfel'd element.

\subsection[The modular group and the double]{} \label{ModGrDoub}
As we have discussed there, the projective representation of the modular group on the center of~$A$ described in Corollary~\ref{RibEl} is not induced by a linear representation in general. However, the situation is better for a certain tensor product:
\begin{lemma}
There is a unique homomorphism from~$\SL(2,\Z)$ to~$\GL(Z(A) \o Z(A))$ 
that maps~$\gv$ to~$\v \o \v^{-1}$ and~$\gt$ to~$\t \o \t^{-1}$.
\end{lemma}
\begin{pf}
As in the proof of Corollary~\ref{RibEl}, we have to check the defining relations $\gv^4 = 1$ and $(\gt\gv)^3 = \gv^2$. For the first relation, we have by Corollary~\ref{InvSigma} that
$$(\v \o \v^{-1})^2 = 
\frac{(\rha \o \rha)(R'R)}{(\rha \o \rha)(R'R)} \; \sa \o \sa^{-1}$$
which implies the assertion, since $\sa^2=\id$ on the center.
It should be noted that, in contrast to~$\v$, the endomorphism~$\v \o \v^{-1}$ is independent of the choice of an integral.
The second defining relation can be rewritten in the form $\gv\gt\gv = \gt^{-1} \gv \gt^{-1}$, 
as in the proof of Corollary~\ref{RibEl}. This now follows from Proposition~\ref{RibEl}, too, as we have
$$
(\v \o \v^{-1}) \circ (\t \o \t^{-1}) \circ (\v \o \v^{-1}) =
\frac{\rha(v)}{\rha(v)}
(\t^{-1} \o \t) \circ (\v \o \v^{-1}) \circ (\t^{-1} \o \t)
$$
and the factors~$\rha(v)$ involved now cancel.
\qed
\end{pf}
The associated projective representation on the projective space
$P(Z(A) \o Z(A))$ is the tensor product of two projective representations:
The first is the one constructed in Corollary~\ref{RibEl}, and the second is the first one twisted by the conjugation with the matrix~$\ga$ described in Paragraph~\ref{GenRel}. This implies that when
$\rha(v)=\rha(v^{-1})=1$, in which case these two projective representations lift to ordinary linear representations, we can write
$$g.(z \o z') = g.z \o \tilde{g}.z'$$
This equation holds because it suffices to check it on generators,
and for the generators we observed in Paragraph~\ref{GenRel} that
$\tilde{g} = g^{-1}$.

From Corollary~\ref{RibEl}, we also get a projective representation of the modular group on the center of the Drinfel'd double~$D(A)$, using the integral~$\rhd = \Psi^*(\rha \o \rha)$ and the ribbon element~$v_D=\Psi^{-1}(v \o v^{-1})$ introduced in Paragraph~\ref{IntRibDouble}. Suppose now that~$\rha$ is normalized so that
$(\rha \o \rha)(R'R) = 1$. By Lemma~\ref{RibEl}, we then have
$$\rhd(v_D) = \rha(v) \rha(v^{-1}) = (\rha \o \rha)(R'R) = 1$$
and similarly also that $\rhd(v_D^{-1}) = 1$. Therefore again by Lemma~\ref{RibEl}, we can conclude for the R-matrix of the Drinfel'd double that~$(\rhd \o \rhd)(R'R)=1$. By the discussion in Paragraph~\ref{RibEl}, this means that the projective representation on~$Z(D(A))$ is induced from an ordinary linear representation.
Clearly, $\Psi$ is equivariant with respect to this action
and the one considered in the preceding lemma:
\begin{prop}
For all $g \in \Gamma$ and all $z \in Z(D(A))$, we have
$\Psi(g.z) = g.\Psi(z)$.
\end{prop}
\begin{pf}
It suffices to check this on generators, i.e., in the case~$g=\gv$ and~$g=\gt$. But in these cases the assertion is exactly what we have established in Paragraph~\ref{IntRibDouble}, because
$\sa \circ \v  = \v^{-1}$ by Corollary~\ref{InvSigma}.
\qed
\end{pf}

It should be noted that in the case $\rha(v)=1$, in which the projective representation on~$Z(A)$ lifts to a linear representation, the formula in the proposition can be written as
$$\Psi(g.z) = (g \o \tilde{g}).\Psi(z)$$

\newpage
\section{The semisimple case} \label{Sec:FaktSemisim}
\subsection[The character ring]{} \label{CharRing}
Let us now assume that our quasitriangular Hopf algebra~$A$ is factorizable, semisimple, and that the base field~$K$ is algebraically closed of characteristic zero. In this case, $A$ is also cosemisimple and the antipode is an involution.\endnote{\cite{LR2}, Thm.~3.3, p.~276; \cite{LR1}, Thm.~4, p.~195.} By Wedderburn's theorem,\endnote{\cite{FarbDennis}, Thm.~1.11, p.~40.} we can decompose~$A$ into a direct sum of simple two-sided ideals:
$$A = \bigoplus_{i=1}^k \; I_i$$
For every $i=1,\ldots,k$, we can then find a simple module such that the corresponding representation maps~$I_i$ isomorphically to~$\End(V_i)$ and vanishes on the other two-sided ideals~$I_j$ if~$j \neq i$. We denote the dimension of~$V_i$ by~$n_i$. We can assume that $V_1=K$, the base field, considered as a trivial module via the counit. We denote the character of~$V_i$ by~$\chi_i$, so that, for $a \in A$, 
$$\chi_i(a) := \Tr(a \mid_{V_i})$$
is the trace of the action of~$a$ on~$V_i$. We then have that the character~$\chi_R$ of the regular representation, i.e., the representation given by left multiplication on~$A$ itself, has the form
$$\chi_R(a) = \sum_{i=1}^k n_i \chi_i(a)$$
This character is a two-sided integral in~$A^*$.\endnote{\cite{LR2}, Prop.~2.4, p.~273.}

The subspace of $A^*$ spanned by the characters $\chi_1,\ldots,\chi_k$ is called the character ring of $A$, and is denoted by $\Ch(A)$. It is easy to see that it really is a subalgebra of $A^*$, which consists precisely of the cocommutative elements. Because the antipode is an involution, this means that the character ring~$\Ch(A)$ coincides with both of the algebras~$C(A)$ and~$\bar{C}(A)$ introduced in Paragraph~\ref{FactHopf}, and Proposition~\ref{FactHopf} therefore asserts that~$\Phi$ induces an isomorphism between the character ring and the center~$Z(A)$, which is spanned by the centrally primitive idempotents~$e_i \in I_i$. The first idempotent~$e_1$ is then a two-sided integral normalized such 
that~$\ea(e_1)= 1$. Note that it follows from the discussion in 
Paragraph~\ref{FactHopf} that the restrictions of~$\Phi$ and~$\bar{\Phi}$
to~$\Ch(A)$ are equal, because the grouplike element $g:=\ua\sa(\ua^{-1})$
is the unit element in this case.\endnote{\cite{DrinfAlmCocom}, Prop.~6.2, p.~337; \cite{M}, Prop.~10.1.14, p.~183; see also \cite{EG1}, Eq.~(3), p.~192.} The center~$Z(A)$ is a commutative semisimple algebra that admits exactly~$k$ distinct algebra homomorphisms~$\omega_1,\ldots,\omega_k$ to the base field, which are explicitly given as
$$\omega_i: Z(A) \rightarrow K,~z \mapsto \frac{1}{n_i} \chi_i(z)$$
These mappings are called the central characters; they satisfy $\omega_i(e_j)=\delta_{ij}$. Because~$\Phi$ is an algebra isomorphism between~$\Ch(A)$ and~$Z(A)$, $\Ch(A)$ is also a commutative semisimple algebra,\endnote{\cite{Z2}, Lem.~2, p.~55; \cite{YYY2}, Prop.~6.2, p.~44.}
whose~$k$ distinct algebra homomorphisms~$\xi_1,\ldots,\xi_k$ to the base field are given as~$\xi_i := \omega_i \circ \Phi$. The primitive idempotents~$p_1,\ldots,p_k$ of the character ring are accordingly 
given as $p_j:=\Phi^{-1}(e_j)$ and satisfy~$\xi_i(p_j)=\delta_{ij}$.
The first primitive idempotent is then proportional to the character of the regular representation; more precisely, we have
$\chi_R = n p_1$.\endnote{\cite{So4}, Prop.~3.3, p.~208.}

Because the pairing between the character ring and the center is nondegenerate, a linear functional on the character ring can be uniquely represented by an element in the center. Therefore, there exist
elements\endnote{\cite{With}, Eq.~(4.1), p.~888.} 
$z_1,\ldots,z_k \in Z(A)$ such that 
$\xi_i(\chi) = \chi(z_i)$ for all $\chi \in \Ch(A)$. They are explicitly
given as
$$z_i = \sum_{j=1}^k \frac{\xi_i(\chi_j)}{n_j} e_j$$
We will call $z_1,\ldots,z_k$ the class sums, as they are related to the normalized conjugacy class sums in the group ring of a finite group. Note that~$z_1=\A$.

For every simple module~$V_i$, its dual~$V_i^*$ is again simple. Therefore, there is a unique index~$i^* \in \{1,\ldots,k\}$ such 
that $V_i^* \cong V_{i^*}$. The character of this module will also be denoted by~$\chi_i^* := \chi_{i^*}$. The map $i \mapsto i^*$ is an involution on the index set
$\{1,\ldots,k\}$. Because the use of the antipode in the definition of the dual module, characters, centrally primitive idempotents, and central characters behave as follows with respect to dualization:
$$\chi_{i^*} = \sa^*(\chi_i) \qquad e_{i^*} = \sa(e_i) \qquad 
\omega_{i^*} = \omega_i \circ \sa$$
We have derived in Paragraph~\ref{FactHopf} that
$\sa(\Phi(\varphi)) = \bar{\Phi}(\sa^{-1*}(\varphi))$. 
Since~$A$ is involutory and $\Phi$ and~$\bar{\Phi}$ agree on the character ring, we get furthermore that
$$\xi_{i^*} = \xi_i \circ \sa^* \qquad p_{i^*} = \sa^*(p_i)$$
which implies the formula $z_{i^*} = \sa(z_i)$ for the class sums.

Using duals, we can express the character~$\chi_A$ of the left adjoint representation in the form\endnote{\cite{So4}, Subsec.~3.3, p.~208.}
$$\chi_A = \sum_{i=1}^k \chi_i \chi_i^*$$
This in turn enables to invert the expansion 
$\chi_j = \sum_{i=1}^k \xi_i(\chi_j) p_i$
of the characters in terms of the idempotents:
\begin{prop}
For $i=1,\ldots,k$, we have
$$p_i = \frac{1}{\xi_i(\chi_A)} \sum_{j=1}^k \xi_i(\chi_j) \chi_{j^*}$$
\end{prop}
\begin{pf}
The element $\sum_{j=1}^k \chi_j \o \chi_{j^*}$ is a Casimir element;\endnote{\cite{So4}, Prop.~3.5, p.~211.} i.e., it satisfies
$$\sum_{j=1}^k \chi \chi_j \o \chi_{j^*} 
= \sum_{j=1}^k \chi_j \o \chi_{j^*} \chi$$
for all $\chi \in \Ch(A)$. Applying $\xi_i$ to the first tensor factor, we get
$$\xi_i(\chi) \sum_{j=1}^k  \xi_i(\chi_j) \chi_{j^*} 
= \bigl(\sum_{j=1}^k \xi_i(\chi_j) \chi_{j^*}\bigr) \chi$$
This shows that $\sum_{j=1}^k \xi_i(\chi_j) \chi_{j^*}$ is proportional to~$p_i$. Since $\xi_i(p_i)=1$, we find that the proportionality factor is
$\xi_i(\sum_{j=1}^k \xi_i(\chi_j) \chi_{j^*}) = \xi_i(\chi_A)$. This proportionality factor cannot be zero, since the element itself is not zero, so the assertion follows. 
\qed
\end{pf}

\subsection[The Verlinde matrix]{} \label{VerlMat}
As $A$ is involutory, $\ua$ is central.\endnote{\cite{M}, Prop.~10.1.4, p.~179; \cite{Tur}, Sec.~XI.2.2, Eq.~(2.2.c), p.~498; \cite{Kas}, Prop.~VIII.4.1, p.~180.} As explained in Paragraph~\ref{CharRing}, $\ua$ is also invariant under the antipode, and therefore we can in this situation use $\ua^{-1}$ as a ribbon element.
For this ribbon element, the quantum trace coincides with the usual trace,\endnote{\cite{Tur}, Lem.~XI.3.3, p.~501; \cite{Kas}, Prop.~XIV.6.4, p.~363.} so that the categorical dimensions coincide with the ordinary dimensions~$n_i$ introduced above. Clearly, we can expand the Drinfel'd element and its inverse in terms of the centrally primitive idempotents:
$$\ua = \sum_{i=1}^k u_i e_i \qquad \ua^{-1} = \sum_{i=1}^k \frac{1}{u_i} e_i$$
for numbers $u_i \in K$, which are nonzero because the Drinfel'd element is invertible. Using these numbers, we define the diagonal matrix\endnote{\cite{Tur}, Sec.~II.3.9, p.~98.}
$$\T := (\frac{1}{u_i} \delta_{ij})_{i,j=1,\ldots,k}$$
Because for this ribbon element~$\t$ is the multiplication by~$\ua^{-1}$, this matrix represents the restriction of~$\t$ to the center with respect to the basis consisting of the centrally primitive idempotents.
Furthermore, we will need an auxiliary matrix, the so-called charge conjugation matrix~$\C:=(\delta_{i,j^*})_{i,j=1,\ldots,k}$, which is the matrix representation of the action of the antipode on the center of~$A$ with respect to the basis consisting of the centrally primitive idempotents.
It is also the matrix representation of the action of the dual antipode on the character ring with respect to the basis consisting of the irreducible characters.

We define still another matrix, the so-called Verlinde matrix~$\V$, which should not be confused with the antipode:\endnote{\cite{Tur}, Sec.~II.1.4, p.~74; note the difference to \cite{EG1}, p.~192 and \cite{SchneiderFact}, Rem.~3.4, p.~1895.} 
\begin{defn}
The Verlinde matrix is the matrix $\V=(s_{ij})_{i,j=1,\ldots,k}$ with entries
$$s_{ij} := (\chi_i \o \chi_j)(R'R) =\chi_i(\Phi(\chi_j))$$
\end{defn}

We list some well-known properties of the Verlinde matrix:\endnote{\cite{BakKir}, Eq.~(3.1.3), p.~48; \cite{SchneiderFact}, Rem.~3.4, p.~1895; \cite{Tur}, p.~74f, p.~90.}

\pagebreak

\begin{lemma}
The Verlinde matrix is invertible. Its entries satisfy
\begin{enumerate}
\item 
$s_{ij} = s_{ji}$

\item 
$s_{ij} = s_{i^* j^*}$

\item 
$s_{ij} = n_i \xi_i(\chi_j)$
\end{enumerate}
\end{lemma}
\begin{pf}
The first property follows from the trace property of the characters. 
The second property can be deduced from the fact that
$(\sa \o \sa)(R)=R$. The third property follows from the definitions:
$$\xi_i(\chi_j) = \omega_i(\Phi(\chi_j)) 
= \frac{1}{n_i} \chi_i(\Phi(\chi_j))= \frac{1}{n_i} s_{ij}$$
This also shows that the Verlinde matrix is invertible: Expanding the characters in terms of the idempotents, we have
$\chi_j = \sum_{i=1}^k \xi_i(\chi_j) p_i$,
so that the matrix $(\xi_i(\chi_j))$ is invertible as a base change matrix, and the Verlinde matrix is, by the third property, the product of this matrix and an invertible diagonal matrix.
\qed
\end{pf}

In contrast to the matrix~$\T$, the Verlinde matrix is not exactly the matrix representation of~$\v$ with respect to the centrally primitive idempotents, although these two matrices are closely related.
To understand this relation, recall that~$\v$ depends, via~$\iota$, on the choice of an integral~$\rha \in A^*$. Because the space of integrals is one-dimensional, $\rha$ has to be proportional to the character~$\chi_R$
of the regular representation, so that $\rha = \kappa \chi_R$ for a nonzero number~$\kappa \in K$. Although it is in principle possible to fix~$\kappa$ by normalizing the integral in some way, we will see that it is convenient not to do that at the moment and to keep~$\kappa$ as a free parameter. With this parameter introduced, let us see how the maps $\Phi$ and~$\iota$ behave with respect to the new bases introduced in Paragraph~\ref{CharRing}:
\begin{prop}
For all $i=1,\ldots,k$, we have
\begin{enumerate}
\item 
$\displaystyle
\Phi(\chi_i) = n_i z_i$

\item
$\displaystyle
\iota(e_i) = \kappa n_i \chi_i$

\item
$\displaystyle
\iota(z_i) = \kappa \xi_i(\chi_A) p_{i^*}$
\end{enumerate}
\end{prop}
\begin{pf}
For the first assertion, we note that by the above lemma 
$$n_i \xi_i(\chi_j) = s_{ij} = s_{ji} =  n_j \xi_j(\chi_i)$$
so that the definition of the class sums from Paragraph~\ref{CharRing}
can be rewritten in the form
$$z_i = \sum_{j=1}^k \frac{\xi_i(\chi_j)}{n_j} e_j 
= \sum_{j=1}^k \frac{\xi_j(\chi_i)}{n_i} e_j$$
Expanding $\Phi(\chi_i)$ in terms of centrally primitive idempotents, we therefore have
$$\Phi(\chi_i) = \sum_{j=1}^k \omega_j(\Phi(\chi_i)) e_j
= \sum_{j=1}^k \xi_j(\chi_i) e_j = n_i z_i$$
For the second assertion, we have by definition of~$\iota$ that
\begin{align*}
\iota(e_i)(a) = \kappa \chi_R(e_i a)
= \kappa \sum_{j=1}^k n_j \chi_j(e_i a)
= \kappa n_i \chi_i(a)
\end{align*}
The third assertion follows from the second assertion, together with Proposition~\ref{CharRing} and the formula for the class sums given in that paragraph. We then find
\begin{align*}
\iota(z_i) = \sum_{j=1}^k \frac{\xi_i(\chi_j)}{n_j} \iota(e_j)
&= 
\sum_{j=1}^k \frac{\xi_i(\chi_j)}{n_j} \kappa n_j \chi_j \\
&= \kappa \sum_{j=1}^k \xi_i(\chi_j) \chi_j
= \kappa \xi_i(\chi_A) p_{i^*}
\end{align*}
where we have used that $\xi_{i^*}(\chi_{j^*}) = \xi_i(\chi_j)$
and $\xi_{i^*}(\chi_A) = \xi_i(\chi_A)$.
\qed
\end{pf}

From this proposition, we can deduce the precise relation of the Verlinde matrix and the matrix representation of~$\v$ resp.~$\v_*$. Up to scalar multiples, $\v$ maps idempotents to class sums, and vice versa. Similarly, $\v_*$ maps idempotents to multiples of characters and characters to multiples of idempotents:
\begin{corb}
\begin{enumerate}
\item 
$\displaystyle 
\v(z_j) = \kappa \xi_j(\chi_A) e_j
= \kappa \sum_{i=1}^k \frac{n_i}{n_j} s_{ji^*} z_{i}$

\item 
$\displaystyle 
\v(e_j) = \kappa n_j^2 z_{j^*}
= \kappa \sum_{i=1}^k \frac{n_j}{n_i} s_{j^*i} e_i$

\item 
$\displaystyle 
\v_*(\chi_j) = \kappa n_j \xi_j(\chi_A) p_j
= \kappa \sum_{i=1}^k s_{ji^*} \chi_i$

\item 
$\displaystyle 
\v_*(p_j) = \kappa n_j \chi_{j^*} 
= \kappa \sum_{i=1}^k \frac{n_j}{n_i} s_{ij^*} p_i$
\end{enumerate}
\end{corb}
\begin{pf}
If we apply~$\Phi$ to the formula in Proposition~\ref{CharRing} and use the preceding proposition, then we get
\begin{align*}
e_j = \frac{1}{\xi_j(\chi_A)} \sum_{i=1}^k \xi_j(\chi_i) n_i z_{i^*}
= \frac{1}{\xi_j(\chi_A)} \sum_{i=1}^k \frac{n_i}{n_j} s_{ji} z_{i^*}
\end{align*}
Since~$\Phi$ and~$\bar{\Phi}$ coincide on the character ring, we get from the definition of~$\v$ in Paragraph~\ref{RoleInt} that
\begin{align*}
\v(z_j) = \kappa \xi_j(\chi_A) \sa(\bar{\Phi}(p_{j^*}))
= \kappa \xi_j(\chi_A) e_j
\end{align*}
Combining these two formulas, we get the first statement.
For the second statement, we get from the preceding proposition that
$$\v(e_j) = \kappa n_j \sa(\bar{\Phi}(\chi_j))
= \kappa n_j^2 z_{j^*}
= \kappa n_j^2 
\sum_{i=1}^k \frac{\xi_{j^*}(\chi_i)}{n_i} e_i$$
For the third statement, recall that by its definition in Paragraph~\ref{RoleInt} we have
$\v_* = \sa^* \circ \iota \circ \Phi$, so that the preceding proposition gives
\begin{align*}
\v_*(\chi_j) &= n_j \sa^*(\iota(z_j))
= \kappa n_j \xi_j(\chi_A) p_j 
= \kappa \sum_{i=1}^k n_j \xi_j(\chi_{i^*}) \chi_i
= \kappa \sum_{i=1}^k s_{ji^*} \chi_i
\end{align*}
where the third equation follows from Proposition~\ref{CharRing}. For the fourth statement, we have
\begin{align*}
\v_*(p_j) &= 
\sa^*(\iota(e_j)) = \kappa n_j \chi_{j^*} 
= \kappa \sum_{i=1}^k n_j \xi_i(\chi_{j^*}) p_i
= \kappa \sum_{i=1}^k \frac{n_j}{n_i} s_{ij^*} p_i
\end{align*}
as asserted.
\qed
\end{pf}

\subsection[Matrix identities]{} \label{MatIdent}
The fact that the matrices~$\V$ and~$\T$ are essentially the matrix representations of~$\v$ and~$\t$ implies that they essentially satisfy the defining relations of the modular group. More precisely, 
they satisfy the following relations:\endnote{\cite{Tur}, Sec.~II.3.9,  p.~98.}
\begin{prop}
$$\V^2 = \dim(A) \; \C \qquad \qquad 
\V \T \V = \chi_R(\ua^{-1})\; \T^{-1} \V \C \T^{-1}$$
\end{prop}
\begin{pf}
By Corollary~\ref{VerlMat}, we have
\begin{align*}
\v^2(e_j) = 
\kappa \sum_{l=1}^k \frac{n_j}{n_l} s_{j^*l} \v(e_l) = 
\kappa^2 \sum_{i,l=1}^k \frac{n_j}{n_i} s_{j^*l} s_{i^*l} e_i
\end{align*}
On the other hand, it follows from Corollary~\ref{InvSigma} and Lemma~\ref{RibEl} that
\begin{align*}
\v^2(a) = \rha(\ua) \rha(\ua^{-1}) \sa(a) 
\end{align*}
Inserting $a=e_j$ into this equation and comparing it with the preceding one, we find 
\begin{align*}
\rha(\ua) \rha(\ua^{-1}) e_{j^*} = 
\kappa^2 \sum_{i,l=1}^k \frac{n_j}{n_i} s_{j^*l} s_{i^*l} e_i 
\end{align*}
which implies
$\rha(\ua) \rha(\ua^{-1}) \delta_{ij^*} = 
\kappa^2 \frac{n_j}{n_i} \sum_{l=1}^k  s_{j^*l} s_{i^*l} $
by comparing coefficients. Now note that by Lemma~\ref{RibEl}
$$\frac{1}{\kappa^2} \rha(\ua) \rha(\ua^{-1}) = 
\chi_R(\ua) \chi_R(\ua^{-1}) = \chi_R(\Phi(\chi_R)) = \dim(A)$$
because $\Phi(\chi_R)$ is an integral satisfying 
$\ea(\Phi(\chi_R)) = \dim(A)$ by Lemma~\ref{IntFact}.
This shows that\endnote{\cite{EG1}, Lem.~1.2, p.~193; \cite{SchneiderFact}, Rem.~3.4, p.~1895; \cite{Tur}, Sec.~II.3.8, Eq.~(3.8.a), p.~97.}
\begin{align*}
\sum_{l=1}^k s_{il} s_{lj} = \dim(A) \delta_{ij^*}
\end{align*}
which is the first assertion.

For the second assertion, recall that
$\v \circ \t \circ \v = \rha(u^{-1}) \; \t^{-1} \circ \v \circ \t^{-1}$ 
by Proposition~\ref{RibEl}. By Corollary~\ref{VerlMat}, we have
\begin{align*}
(\v \circ \t \circ \v)(e_j) 
= \kappa \sum_{l=1}^k \frac{n_j}{n_l} s_{j^*l} \frac{1}{u_l} \v(e_l)
= \kappa^2 
\sum_{i,l=1}^k \frac{n_j}{n_i} \frac{s_{j^*l} s_{l^*i}}{u_l} e_i
\end{align*}
as well as
\begin{align*}
(\t^{-1} \circ \v \circ \t^{-1})(e_j) 
= \kappa \sum_{i=1}^k \frac{n_j}{n_i} u_i u_j s_{j^*i} e_i
\end{align*}
Comparing coefficients, we find that\endnote{\cite{Tur}, Sec.~II.3.8, Eq.~(3.8.c), p.~97.}
\begin{align*}
\kappa \sum_{l=1}^k \frac{s_{j^*l} s_{l^*i}}{u_l} 
= \rha(u^{-1}) \; u_i u_j s_{j^*i} 
\end{align*}
or alternatively that
$\sum_{l=1}^k \frac{s_{jl} s_{li}}{u_l} 
= \chi_R(\ua^{-1}) u_i u_j s_{ji^*}$, which gives the second relation.
\qed
\end{pf}

In the proof of the first matrix identity above, we have used one of the two formulas for~$\v(e_j)$ given in Corollary~\ref{VerlMat}. Using the other form reveals another interesting identity:
\begin{corollary}
$\displaystyle
\xi_i(\chi_A) = \frac{\dim(A)}{n_i^2}$
\end{corollary}
\begin{pf}
From Corollary~\ref{VerlMat}, we have 
\begin{align*}
\v^2(e_j) = \kappa n_j^2 \v(z_{j^*}) 
= \kappa^2 n_j^2 \xi_j(\chi_A) e_{j^*}
\end{align*}
But in the proof of the preceding proposition, we have derived that
$$\v^2(e_j) = \rha(\ua) \rha(\ua^{-1}) e_{j^*}$$
and also that $\rha(\ua) \rha(\ua^{-1}) = \kappa^2 \dim(A)$. Therefore, the assertion follows by comparing coefficients.
\qed
\end{pf}
It should be noted that $\xi_i(\chi_A)$ is an eigenvalue of the multiplication with the character~$\chi_A$ corresponding to the eigenvector~$p_i$, and therefore an algebraic integer. As the corollary shows, it is also a rational number, and therefore an integer. This
yields the well-known\endnote{\cite{EG1}, Thm.~1.4, p.~193; \cite{SchneiderFact}, Thm.~3.2, p.~1894; \cite{Tak}, Thm.~5.7, p.~641; \cite{TsangZhu}, Thm.~3, p.~5.} result that~$n_i^2$ divides~$\dim(A)$.

It should furthermore be noted that the preceding corollary can also be
used to give a different proof of the equation $\V^2 = \dim(A) \; \C$, as we have
\begin{align*}
\sum_{l=1}^k s_{il} s_{lj} = 
\sum_{l=1}^k n_i \xi_i(\chi_l) n_j \xi_j(\chi_l) =
n_i n_j \xi_j(\chi_A) \xi_i(p_{j^*}) =
n_j^2 \xi_j(\chi_A) \delta_{ij^*}
\end{align*}
where the first equation follows from Lemma~\ref{VerlMat} and the second equation from Proposition~\ref{CharRing}.

\subsection[A comparison]{} \label{CompTuraev}
We proceed to carry out a more precise comparison of our setup with the setup in \cite{Tur}. For this, we need the following lemma:
\begin{lemma}
$\displaystyle \sum_{i=1}^k n_i u_i^{-1} u_j^{-1} s_{ij} = n_j \chi_R(\ua^{-1})$
\end{lemma}
\begin{pf}
Since $\da(\ua) = (R' R)^{-1} (\ua \o \ua) = (\ua \o \ua)(R' R)^{-1}$, we have
\begin{align*}
\sum_{i=1}^k n_i u_i^{-1} u_j^{-1} s_{ij} &=
\sum_{i=1}^k n_i u_i^{-1} u_j^{-1} (\chi_i \o \chi_j)(R' R) 
= \sum_{i=1}^k n_i (\chi_i \o \chi_j)(\da(\ua^{-1})) \\
&= (\chi_R \o \chi_j)(\da(\ua^{-1})) =
(\chi_R \chi_j)(\ua^{-1}) = n_j \chi_R(\ua^{-1})
\end{align*}
where the last equality follows from the fact that the character of the regular representation is an integral.
\qed
\end{pf}

As we have already pointed out in Paragraph~\ref{VerlMat}, categorical dimensions and ordinary dimensions coincide for our choice of a ribbon element, so that the numbers $\dim(i)$ introduced in \cite{Tur}, Sec.~II.1.4, p.~74 are equal to~$n_i$. Also, since our ribbon element is
$\ua^{-1}$, it is clear that the numbers~$v_i$ and~$\Delta$ introduced in \cite{Tur}, Sec.~II.1.6, p.~76 are equal to~$1/u_i$ resp.~$\chi_R(\ua)$. It therefore follows from the preceding lemma that the parameters~$d_i$ introduced in  
\cite{Tur}, Sec.~II.3.2, p.~87 are in our case equal to~$d_i = n_i/\chi_R(\ua^{-1})$, which is in accordance with \cite{Tur}, Lem.~II.3.2.3, p.~89. For the rank $\cal D$, we have the two choices 
${\cal D}=\pm\sqrt{\dim(A)}$. The equation $\dim(A)=\chi_R(\ua)\chi_R(\ua^{-1})$ observed in Paragraph~\ref{MatIdent}
then becomes the equation~$\Delta = d_0 {\cal D}^{2}$ in 
\cite{Tur}, Sec.~II.3.2, Eq.~(3.2.j), p.~89.

\subsection[Radford's example]{} \label{ExamRadf}
We now illustrate the preceding considerations by inspecting an example given by D.~E.~Radford.\endnote{\cite{RadfAntiQuasi}, Sec.~3, p.~10; \cite{RadfKnotInv}, Sec.~2.1, p.~219.} Consider a cyclic group~$G$ of order~$n$. Denote the group ring by~$A=K[G]$, and fix a generator~$g$ of~$G$. As~$A$ is cocommutative, $A$ is certainly quasitriangular with respect to the R-matrix~$\A \o \A$. However, with respect to this R-matrix, it is not factorizable. Radford has determined all possible R-matrices for~$A$, and shown that~$A$ can only be factorizable if~$n$ is odd, what we will assume for the rest of this paragraph, and that in this case the R-matrix necessarily has the form
$$R = \frac{1}{n} \sum_{i,j=0}^{n-1} \zeta^{-ij} g^i \o g^{j}$$
where~$\zeta$ is a primitive $n$-th root of unity.\endnote{\cite{RadfKnotInv}, Sec.~2.3, p.~227.}
To follow his convention, we will deviate in this paragraph from the enumeration introduced in Paragraph~\ref{CharRing} and enumerate the 
(centrally) primitive idempotents in the form $e_0,\ldots,e_{n-1}$ instead of $e_1,\ldots,e_{k}$; note that~$n=k$ in the present situation. They are then given by the formula\endnote{\cite{RadfKnotInv}, Sec.~1.1, p.~210.}
$$e_j = \frac{1}{n} \sum_{i=0}^{n-1} \zeta^{-ij} g^i$$
The irreducible characters are determined by~$\chi_i(e_j) = \delta_{ij}$
and therefore satisfy $\chi_i(g) = \zeta^i$. Radford also gives the following formulas for the Drinfel'd element and its inverse:\endnote{\cite{RadfKnotInv}, Sec.~1.1, p.~211; Sec.~2.1, p.~219; Sec.~2.3, p.~227.}
$$u = \sum_{i=0}^{n-1} \zeta^{-i^2} e_i \qquad \qquad
u^{-1} = \sum_{i=0}^{n-1} \zeta^{i^2} e_i$$
He also gives the formula
$\da(u^{-1}) = \sum_{i,j=0}^{n-1} \zeta^{(i+j)^2} e_i \o e_j$
for the coproduct of the inverse Drinfel'd element, from which we get that
$$R' R = (\ua \o \ua) \da(u^{-1}) =
\sum_{i,j=0}^{n-1} \zeta^{2ij} e_i \o e_j$$
This means that the entries of the Verlinde matrix are given as
$$s_{ij} = (\chi_i \o \chi_j)(R' R) = \zeta^{2ij}$$
so that, using Corollary~\ref{VerlMat}, we find the expressions 
$$\v(e_j) =  \kappa
\sum_{i=0}^{n-1} \zeta^{-2ij} e_i \qquad \qquad
\t(e_j) = \zeta^{j^2} e_j$$
for the mappings~$\v$ and~$\t$, since~$\t$ is the multiplication by~$u^{-1}$.

The reason for mentioning this example is its following feature:
\begin{prop}
We have
$$\chi_R(u^{-1}) =
\begin{cases}
\chi_R(u) &\text{if} \mspace{20mu} n \equiv 1 \pmod{4} \\
-\chi_R(u) &\text{if} \mspace{20mu} n \equiv 3 \pmod{4}
\end{cases}
$$
\end{prop}
\begin{pf}
From the form of the inverse Drinfel'd element given above, we see that
$\chi_R(u^{-1}) = \sum_{i=0}^{n-1} \zeta^{i^2}$ is the quadratic Gaussian sum. As the quadratic Gaussian sum transforms with the Jacobi symbol under change of the root of unity,\endnote{\cite{Nag}, Chap.~V, Exerc.~122, p.~187; \cite{RadElemNum}, Chap.~11, Thm.~38, p.~93.} we have
$$\chi_R(u^{-1}) = \jac{-1}{n} \chi_R(u)$$
The assertion now follows from the first supplement to Jacobi's reciprocity law.\endnote{\cite{Lan}, Chap.~I.VI, Thm.~91, p.~66; \cite{Nag}, Chap.~IV, Sec.~42, p.~146; \cite{RadElemNum}, Chap.~11, p.~91.}~\qed
\end{pf}

For the discussion in Paragraph~\ref{RibEl}, this result means that the two conditions $(\rha \o \rha)(R'R) = 1$ and $\rha(v)=1$ can not always be simultaneously satisfied. It also means that Lemma~\ref{RibEl} can, in a sense, be considered as a generalization of the formula for the absolute value of the quadratic Gaussian sum.\endnote{\cite{Lan}, Chap.~IV.VI.2, p.~208; \cite{RadElemNum}, Chap.~11, Eq.~(11.7), p.~88.} We will further elaborate on this analogy in Paragraph~\ref{HopfSymb}.

\newpage
\section{The case of the Drinfel'd double} \label{Sec:CaseDrinfDoub}
\subsection[The role of the evaluation form]{} \label{RoleEval}
In the case where~$A$ is the Drinfel'd double~$D=D(H)$ of a semisimple Hopf algebra, it is possible to give another description of the action of the modular group that will play an important role in the sequel. We therefore suppose now that $H$ is a semisimple Hopf algebra over an algebraically closed field of characteristic zero, and set $A=D(H)$, its Drinfel'd double. First, recall from Paragraph~\ref{IntDrinfDouble} that the two-sided integral of the Drinfel'd double has the form
$\Lmd = \lmh \o \Gh$
for an integral $\lmh \in H^*$ and an integral $\Gh \in H$. 
We can choose these integrals in such a way that~$\eh(\Gh)=1$ and $\lmh(\Gh)=1$. They are then uniquely determined and satisfy $\lmh(\sh(\Gh))=1$ as well as
$\lmh(\H)=\dim(H)$.\endnote{\cite{R4}, Prop.~1.e, p.~587; Prop.~2.c, p.~589.} From Paragraph~\ref{IntDrinfDouble}, we then know 
that the right integral~$\rhd$ on~$D$ given by
$$\rhd(\varphi \o h) = \varphi(\Gh) \lmh(h)$$
satisfies $\rhd(\ud^{-1}) = \rhd(\ud) = 1$, which implies that 
$(\rhd \o \rhd)(R'R)=1$ by Lemma~\ref{RibEl} and $\rhd(\Lmd) = \lmh(\Gh)^2 = 1$. By comparing normalizations, we see that the character of the regular representation 
is $\chi_R = \dim(H) \; \rhd$, so that we have
$$\chi_R(\ud^{-1}) = \chi_R(\ud) = \dim(H)$$
This means on the one hand that the parameter~$\kappa = \kappa_D$, introduced in Paragraph~\ref{VerlMat}, is in the case of the Drinfel'd double with these normalizations given by
$\kappa_D = \frac{1}{\dim(H)}$, and on the other hand means that,
as discussed in Paragraph~\ref{RibEl}, the representation of~$\SL(2,\Z)$ on the center is not only a projective representation, but rather is linear.

Recall from Lemma~\ref{CoprodEvalForm} that, under the correspondence
$H \o H^* \cong D^*$ described there, the restrictions of~$\Phi$ and~$\bar{\Phi}$ to the character ring are just the interchange of the tensorands. From this, and the fact that the antipode is an involution, it is clear that the evaluation form~$e$ introduced in Paragraph~\ref{DrinfDouble}, which is contained in the character ring, is mapped under $\Phi$ and~$\bar{\Phi}$ to the inverse Drinfel'd element~$\ud^{-1}$, which, as discussed in Paragraph~\ref{CharRing}, can be used as a ribbon element.
Another consequence of these considerations is the following fact:
\begin{lemma}
Suppose that $z = \sum_j \varphi_j \o h_j$ is a central element.
Then we have also 
$$z = \sum_j (\eh \o h_j)(\varphi_j \o \H)$$
\end{lemma}
\begin{pf}
Put $\chi := \Phi^{-1}(z)$. By Lemma~\ref{CoprodEvalForm}, we then have
$\chi = \sum_j h_j \o \varphi_j$. But we have also seen in Paragraph~\ref{CoprodEvalForm} that
$\bar{\Phi}(\chi) = \sum_j (\eh \o h_j)(\varphi_j \o \H)$.
Since~$\bar{\Phi}$ agrees with~$\Phi$ on the character ring, the assertion follows.
\qed
\end{pf}

For the right and the left multiplication with the evaluation form, we now introduce the notation~$\t_*$ and~$\bar{\t}_*$. In other words, we define the endomorphisms~$\t_*$ and~$\bar{\t}_*$ of~$D^*$ by
$$\t_*(\psi) = \psi e  \qquad \bar{\t}_*(\psi) = e \psi$$
This notation is justified by the following proposition, which should be compared with Proposition~\ref{RoleInt}:
\begin{prop}
The diagrams
$$\square<1`1`1`1;1000`700>[D^*`D^*`D`D;\t_*`\Phi`\Phi`\t]
\qquad \qquad
\square<1`1`1`1;1000`700>[D^*`D^*`D`D;\bar{\t}_*`\bar{\Phi}`\bar{\Phi}`\t]
$$
are commutative.
\end{prop}
\begin{pf}
The commutativity of the first diagram follows directly from Proposition~\ref{FactHopf}, as we have
$$
\Phi(\t_*(\psi)) = \Phi(\psi e) = \Phi(\psi) \Phi(e) = \Phi(\psi) \ud^{-1}
= \t(\Phi(\psi))$$
Also in Paragraph~\ref{FactHopf}, we saw that
$\bar{\Phi}(e\psi) = \bar{\Phi}(e) \bar{\Phi}(\psi) = \ud^{-1}\bar{\Phi}(\psi)$,
which yields the commutativity of the second diagram.
\qed
\end{pf}

This proposition implies that we also get a representation of the modular group on the character ring of the Drinfel'd double by mapping the generator~$\gv$ to the restriction of~$\v_*$ and the generator~$\gt$ to the restriction of~$\t_*$. This action is, via~$\Phi$, just conjugate
to the action on the center constructed in Corollary~\ref{RibEl}.
Note that~$\t_*$ and~$\bar{\t}_*$ really preserve the character ring, as they are left resp.~right multiplication with the character~$e$. As the character ring is commutative, these two endomorphisms in fact coincide on the character ring.

\subsection[The new maps]{} \label{NewMaps}
The second construction of the modular group action alluded to above is based on the following maps~$\r$ and~$\bar{\r}$, which should not be confused with the R-matrix:
\begin{defn}
We define the endomorphisms~$\r$ and~$\bar{\r}$ of~$D$ by setting
$$\r(a) := e(a_\1) a_\2 \qquad \bar{\r}(a) := e(a_\2) a_\1$$
Furthermore, we define the endomorphism~$\r_*$ of~$D^*$ as
$$\r_*(\psi)(a) = \psi(\ud^{-1} a)$$
\end{defn}
In other words, we set $\r_*=\t^*$, the transpose of~$\t$. Note that we also have $\r^*=\bar{\t}_*$ and $\bar{\r}^*=\t_*$.

These maps are related in a similar way as the ones considered earlier:
\begin{prop}
The diagrams
$$\square<1`1`1`1;1000`700>[D^*`D^*`D`D;\r_*`\Phi`\Phi`\r]
\qquad \qquad
\square<1`1`1`1;1000`700>[D^*`D^*`D`D;\r_*`\bar{\Phi}`\bar{\Phi}`\bar{\r}]$$
are commutative.
\end{prop}
\begin{pf}
\begin{list}{(\arabic{num})}{\usecounter{num} \leftmargin0cm \itemindent5pt}
\item
Recall from Paragraph~\ref{DrinfDouble} the formula 
$\ud^{-1} = \sum_{i=1}^n b^*_i \o b_i$, where $b_1,\ldots,b_n$ is a basis of~$H$ with dual basis $b^*_1,\ldots,b^*_n$. Recall further that we have set up in Paragraph~\ref{RoleEval} a correspondence between~$H \o H^*$ and~$D^*$, so that we can associate with every~$h \in H$ and~$\varphi \in H^*$ an element~$\psi \in D^*$. For this element, we find that 
\begin{align*}
\bar{\Phi}(\psi) &= 
\sum_{i,j=1}^n \psi(b_i^* \o b_j) \; (\eh \o b_i)(b_j^* \o \H) \\
&= \sum_{i,j=1}^n b_i^*(h) \varphi(b_j) \; (\eh \o b_i)(b_j^* \o \H)
= (\eh \o h)(\varphi \o \H)
\end{align*}
which implies that 
\begin{align*}
\bar{\r}(\bar{\Phi}(\psi)) &= 
e((\eh \o h_\2)(\varphi_\1 \o \H)) \; (\eh \o h_\1)(\varphi_\2 \o \H) \\
&= e((\varphi_\1 \o \H)(\eh \o h_\2)) \; (\eh \o h_\1)(\varphi_\2 \o \H) \\
&= \varphi_\1(h_\2) \; (\eh \o h_\1)(\varphi_\2 \o \H)
\end{align*}

\item
On the other hand, for~$a = \varphi' \o h' \in D$, we have
by the centrality of the Drinfel'd element that
\begin{align*}
\r_*(\psi)(a) &= \psi(\ud^{-1}a) = \psi((\varphi' \o \H)\ud^{-1}(\eh \o h')) 
= \sum_{i=1}^n \psi(\varphi'b^*_i \o b_i h') \\
&= \sum_{i=1}^n (\varphi'b^*_i)(h) \varphi(b_i h')
= \sum_{i=1}^n \varphi'(h_\1) b^*_i(h_\2) \varphi_\1(b_i) \varphi_\2(h') \\
&= \varphi'(h_\1) \varphi_\1(h_\2) \varphi_\2(h')
\end{align*}
This means that $\r_*(\psi) \in D^*$ corresponds to
$\varphi_\1(h_\2) h_\1 \o \varphi_\2 \in H \o H^*$.
By the preceding computation of~$\bar{\Phi}$, we therefore have
$$\bar{\Phi}(\r_*(\psi)) 
= \varphi_\1(h_\2) \; (\eh \o h_\1)(\varphi_\2 \o \H)
= \bar{\r}(\bar{\Phi}(\psi))$$
This establishes the commutativity of the second diagram.

\item
The commutativity of the first diagram is a consequence of 
the commutativity of the second. As discussed in Paragraph~\ref{CharRing},
we have $\sd(\ud) = \ud$, and therefore
$\sd^* \circ \r_* = \r_* \circ \sd^*$. From the proof of Lemma~\ref{IntDrinfDouble}, we have $\sd^*(e)=e$, which implies that $\sd \circ \r = \bar{\r} \circ \sd$. As we saw in Paragraph~\ref{FactHopf}, we also have 
$\sd \circ \Phi = \bar{\Phi} \circ \sd^{*}$ as a consequence of involutivity, so that
\begin{align*}
\sd \circ \Phi \circ \r_* &= \bar{\Phi} \circ \sd^* \circ \r_*
= \bar{\Phi} \circ \r_* \circ \sd^* \\
&= \bar{\r} \circ \bar{\Phi} \circ \sd^*
= \bar{\r} \circ \sd \circ \Phi
= \sd \circ \r \circ \Phi
\end{align*}
After cancelling the antipode, this is the commutativity of the first diagram.
\qed
\end{list}
\end{pf}

It is interesting to look at the linearity properties of our new maps:
\begin{lemma}
For all $a,b \in D$, we have
\begin{enumerate}
\item 
$\r(\sd^{-1}(b_\1) a b_\2) = \sd^{-1}(b_\1) \r(a) b_\2$

\item
$\bar{\r}(b_\1 a \sd^{-1}(b_\2)) = b_\1 \bar{\r}(a) \sd^{-1}(b_\2)$
\end{enumerate}
In particular, $\r$ and~$\bar{\r}$ both preserve the center of~$D$.
\end{lemma}
\begin{pf}
By the symmetry property of the evaluation form~$e$ recorded in 
Paragraph~\ref{DrinfDouble}, we have
\begin{align*}
\r(\sd^{-1}(b_\1) a b_\2) &= 
e(\sd^{-1}(b_\2) a_\1 b_\3) \sd^{-1}(b_\1) a_\2 b_\4 \\
&= e(a_\1 b_\3 \sd^{-1}(b_\2)) \sd^{-1}(b_\1) a_\2 b_\4 \\
&= e(a_\1) \sd^{-1}(b_\1) a_\2 b_\2
= \sd^{-1}(b_\1) \r(a) b_\2
\end{align*}
For the second assertion, we have similarly that
\begin{align*}
\bar{\r}(b_\1 a \sd^{-1}(b_\2)) = 
b_\1 a_\1 \sd^{-1}(b_\4) e(b_\2 a_\2 \sd^{-1}(b_\3)) =
b_\1 \bar{\r}(a) \sd^{-1}(b_\2)
\end{align*}
These computations do not use that the antipode is an involution; however, this is necessary for the statement about the center, because~$a$ is central if and only if $\sd(b_\1) a b_\2 = \ed(b) \; a$ for all~$b \in D$, a relation that is then preserved by~$\r$. A similar argument shows that~$\bar{\r}$ preserves the center.
\qed
\end{pf}

\subsection[The first relation]{} \label{FirstRel}
One important step in the second approach toward the action of the modular group is the following relations between our maps:
\begin{prop}
$$\bar{\t}_* \circ \r_* \circ \bar{\t}_* 
= \r_* \circ \bar{\t}_* \circ \r_*
\qquad \qquad
\t_* \circ \r_* \circ \t_* 
= \r_* \circ \t_* \circ \r_*$$
\end{prop}
\begin{pf}
As discussed in Paragraph~\ref{CharRing}, we have $\sd(\ud)=\ud$ in our case. It then follows from Proposition~\ref{IntDrinfDouble}
that $e^{-1} = \rhd_\2(\ud^{-1}) \rhd_\1$. We therefore get
\begin{align*}
\r_*^{-1}(\bar{\t}_*^{-1}(\psi)) = \r_*^{-1}(e^{-1}\psi)
= \rhd_\2(\ud^{-1}) \r_*^{-1}(\rhd_\1 \psi)
\end{align*}
Since $\rhd$ is an integral, we can rewrite this as\endnote{\cite{LR2}, Lem.~1.2, p.~270.}
\begin{align*}
\r_*^{-1}(\bar{\t}_*^{-1}(\psi)) 
&= (\rhd_\2 \sd^*(\psi))(\ud^{-1}) \r_*^{-1}(\rhd_\1) \\
&= \rhd_\3(\ude^{-1}) \sd^*(\psi)(\udz^{-1}) \rhd_\2(\ud) \rhd_\1\\
&= \rhd_\2(\ud\ude^{-1}) \r_*(\sd^*(\psi))(\ud\udz^{-1})  \rhd_\1
\end{align*}
Using the formula that expresses~$\Phi$ in terms of the Drinfel'd element given at the end of Paragraph~\ref{FactHopf}, this can be written as
\begin{align*}
\r_*^{-1}(\bar{\t}_*^{-1}(\psi)) 
&= \rhd_\2(\Phi(\r_*(\sd^*(\psi)))) \rhd_\1 \\
&= \rhd_\2(\r(\Phi(\sd^*(\psi)))) \rhd_\1
= (e \rhd_\2)(\Phi(\sd^*(\psi))) \rhd_\1
\end{align*}
where the second equality follows from Proposition~\ref{NewMaps}.
Using the relation between~$\Phi$ and the Drinfel'd element backwards, this becomes
\begin{align*}
\r_*^{-1}(\bar{\t}_*^{-1}(\psi)) 
= (e \rhd_\2)(\ud \ude^{-1}) \; \sd^*(\psi)(\ud \udz^{-1}) \; \rhd_\1
\end{align*}
Now $\rhd$ is in our case also a left integral, and furthermore we saw in Paragraph~\ref{IntDrinfDouble} that~$e$ is invariant under the antipode, so that we can rewrite the preceding equation in the form
\begin{align*}
\r_*^{-1}(\bar{\t}_*^{-1}(\psi)) 
&= \rhd_\2(\ud \ude^{-1}) \; \sd^*(\psi)(\ud \udz^{-1}) \; e\rhd_\1 \\
&= \rhd_\2(\ud) \rhd_\3(\ude^{-1})\; \r_*^{-1}(\sd^*(\psi))(\udz^{-1}) \; e\rhd_\1 \\
&= (\rhd_\2\; \r_*^{-1}(\sd^*(\psi)))(\ud^{-1}) \; e  \r_*^{-1}(\rhd_\1)
\end{align*}
As discussed in Paragraph~\ref{NewMaps}, we have
$\sd^* \circ \r_* = \r_* \circ \sd^*$, and by using this and the properties of the integral again, we can rewrite this expression as
\begin{align*}
\r_*^{-1}(\bar{\t}_*^{-1}(\psi)) 
&= (\rhd_\2\; \sd^*(\r_*^{-1}(\psi)))(\ud^{-1}) \; e  \r_*^{-1}(\rhd_\1) \\
&= \rhd_\2 (\ud^{-1}) \; e  \r_*^{-1}(\rhd_\1\r_*^{-1}(\psi)) \\
&= e  \r_*^{-1}(e^{-1} \r_*^{-1}(\psi)) 
= \bar{\t}_*(\r_*^{-1}(\bar{\t}_*^{-1}(\r_*^{-1}(\psi))))
\end{align*}
where the third equation uses Proposition~\ref{IntDrinfDouble} again.
This proves $\r_*^{-1} \circ \bar{\t}_*^{-1}
= \bar{\t}_* \circ \r_*^{-1} \circ \bar{\t}_*^{-1} \circ \r_*^{-1}$, which is equivalent to the first assertion. The second assertion follows from the first by conjugating with the antipode~$\sd^*$ of~$D^*$, as we have already noted that~$\sd^*$ commutes with~$\r_*$, and 
$\sd^* \circ \t_* = \bar{\t}_* \circ \sd^*$ holds since~$\sd^*$ is antimultiplicative and preserves~$e$.
\qed
\end{pf}

Instead of using endomorphisms of~$D^*$, we can use endomorphisms of~$D$. The corresponding relation then has the following form:
\begin{corollary}
$$\t \circ \bar{\r} \circ \t = \bar{\r} \circ \t \circ \bar{\r}
\qquad \qquad
\t \circ \r \circ \t = \r \circ \t \circ \r$$
\end{corollary}
\begin{pf}
Using Proposition~\ref{RoleEval} and Proposition~\ref{NewMaps}, this follows by conjugating the first formula in the preceding proposition by~$\bar{\Phi}$ and the second formula in the preceding proposition by~$\Phi$.
\qed
\end{pf}

\subsection[The second approach to the action of the modular group]{} \label{SecApprMod}
Corollary~\ref{FirstRel} suggests that we might get a representation of the modular group by assigning~$\t$ to the generator~$\gt$ and~$\r$ to the generator~$\gr$, the two alternative generators described in Paragraph~\ref{GenRel}. The first defining relation 
$\gt \gr \gt = \gr \gt \gr$ then follows from this corollary. However, we still need the second defining relation~$(\gr \gt)^6 = 1$. This relation only holds for the restrictions of~$\t$ and~$\r$ to the center of~$D$. We now proceed not only to verify this relation, but also to check that the representation of the modular group that we construct in this way agrees with the one constructed earlier. To do this, we introduce the following analogue of~$\iota$:
$$\iota_*: D^* \rightarrow D,~\psi \mapsto \Lmd_\1 \psi(\Lmd_\2)$$
These maps together satisfy the following relations:\endnote{\cite{R4}, p.~588.}
$$(\iota \circ \iota_*)(\psi) = \rhd(\Lmd) \; \sd^*(\psi) 
\qquad 
(\iota_* \circ \iota)(x) = \rhd(\Lmd) \; \sd(x)$$
With the help of this map, we can deduce the following fact, which relates the two approaches to the representation of the modular group:
\begin{prop}
$$\v_* = \t_*^{-1} \circ \r_*^{-1} \circ \t_*^{-1}
= \r_*^{-1} \circ \t_*^{-1} \circ \r_*^{-1}$$
\end{prop}
\begin{pf}
The second equality is just the inversion of the second identity in Proposition~\ref{FirstRel}; it is therefore sufficient to show the first equality. We have the commutation relation
$\t \circ \iota_* = \iota_* \circ \r_*$ because, as discussed in Paragraph~\ref{CharRing}, $\ud^{-1}$ is invariant under the antipode, and therefore we have\endnote{\cite{LR2}, Lem.~1.2, p.~270.}
\begin{align*}
\t(\iota_*(\psi)) &= \ud^{-1} \Lmd_\1 \psi(\Lmd_\2) = 
\Lmd_\1 \psi(\sd^{-1}(\ud^{-1}) \Lmd_\2) \\
&= \Lmd_\1 \psi(\ud^{-1} \Lmd_\2) = \Lmd_\1 \r_*(\psi)(\Lmd_\2) 
= \iota_*(\r_*(\psi))
\end{align*}
It follows from Lemma~\ref{IntDrinfDouble} that
$\iota_*(e^{-1}) = \sd(\ud^{-1}) = \ud^{-1}$,
because $e^{-1}(\Lmd) = \lmh(\sh^{-1}(\Gh)) = 1$ in our case.
From the expression for~$\Phi$ in terms of the Drinfel'd element given in
Paragraph~\ref{FactHopf}, we therefore get
\begin{align*}
\Phi(\psi) &= \ud \ude^{-1} \psi(\ud \udz^{-1}) 
= \t^{-1}(\ude^{-1}) \r_*^{-1}(\psi)(\udz^{-1}) \\
&= \t^{-1}(\Lmd_\1) \; \bigl(\r_*^{-1}(\psi)e^{-1}\bigr)(\Lmd_\2)\\
&= \t^{-1}(\Lmd_\1) \; \bigl(\t_*^{-1}(\r_*^{-1}(\psi))\bigr)(\Lmd_\2)
\\
&= (\t^{-1} \circ \iota_* \circ \t_*^{-1} \circ \r_*^{-1})(\psi)
= (\iota_* \circ \r_*^{-1} \circ \t_*^{-1} \circ \r_*^{-1})(\psi)
\end{align*}
so that, by the definition of~$\v_*$ in Paragraph~\ref{RoleInt} and the properties of~$\iota$ and~$\iota_*$ mentioned above, we have
\begin{align*}
\v_* = \sa^{-1*} \circ \iota \circ \Phi 
= \sa^{-1*} \circ \iota \circ \iota_* \circ \r_*^{-1} \circ \t_*^{-1} \circ \r_*^{-1}
= \r_*^{-1} \circ \t_*^{-1} \circ \r_*^{-1}
\end{align*}
as asserted.
\qed
\end{pf}

It is easy to convert the preceding proposition from a statement about endomorphism of~$D^*$ into a statement about endomorphism of~$D$: If we conjugate the identity by~$\Phi$ and use Proposition~\ref{RoleInt}, 
Proposition~\ref{RoleEval}, and Proposition~\ref{NewMaps}, we get 
$$\v = \t^{-1} \circ \r^{-1} \circ \t^{-1}
= \r^{-1} \circ \t^{-1} \circ \r^{-1}$$

We have proved in Proposition~\ref{RoleInt} that~$\v$ preserves the center, and this is also true for~$\t$ by the centrality of the Drinfel'd element. In Lemma~\ref{NewMaps}, we have seen that~$\r$ preserves the center. We use the same symbols for the restrictions of these maps to the center. We then have
\begin{corollary}
There is a unique homomorphism from~$\SL(2,\Z)$ to~$\GL(Z(D))$ that
maps~$\gr$ to~$\r$ and~$\gt$ to~$\t$.
\end{corollary}
\begin{pf}
The homomorphism is unique because~$\gr$ and~$\gt$ generate the modular group, as discussed in Paragraph~\ref{GenRel}. For the existence question, recall the defining relations $\gt \gr \gt = \gr \gt \gr$ and $(\gr \gt)^6 = 1$. The first relation holds by Corollary~\ref{FirstRel}. We have 
$(\rhd \o \rhd)(R'R) = 1$, and 
therefore Corollary~\ref{InvSigma} yields that the restriction of~$\v^2$ to the center coincides with the antipode. Together with the 
above considerations, this shows that
$$(\r \circ \t)^6 = (\r \circ \t \circ \r \circ \t \circ \r \circ \t)^2 = \v^{-4} = \sd^{-2} = \id$$
on the center, which is the second relation needed.
\qed
\end{pf}

Because $\gv = \gt^{-1} \gr^{-1} \gt^{-1}$, it is clear that this representation of the modular group agrees with the one constructed in Corollary~\ref{RibEl}.

\subsection[Matrix representations of the new maps]{} \label{MatNewMap}
We have discussed the matrix representations of~$\t$ and~$\v$ in Paragraph~\ref{VerlMat}. It is possible to give a similar discussion of the matrix representations of~$\r$, $\bar{\r}$, and~$\r_*$:
\begin{propb}
\begin{enumerate}
\item 
$\displaystyle
\r_*(\chi_i) = \frac{1}{u_i} \chi_i= \xi_i(e) \chi_i$

\item 
$\displaystyle
\r(z_i) = \bar{\r}(z_i) =\frac{1}{u_i} z_i= \xi_i(e) z_i$
\end{enumerate}
\end{propb}
\begin{pf}
For the first assertion, note that 
$$\r_*(\chi_i)(a) = \chi_i(\ud^{-1}a) = \frac{1}{u_i} \chi_i(a)$$
which gives the first equation. The second equation holds since
$$\xi_i(e) = \omega_i(\Phi(e))= \omega_i(\ud^{-1}) = \frac{1}{u_i}$$
The second assertion follows by applying~$\Phi$ to the first assertion and using Proposition~\ref{VerlMat} and Proposition~\ref{NewMaps}; note that we discussed in Paragraph~\ref{RoleEval} that~$\Phi$ and~$\bar{\Phi}$ agree on the character ring.
\qed
\end{pf}

Using this proposition, we can expand the evaluation form explicitly in
terms of the irreducible characters:
\begin{corollary}
\begin{align*}
e = \frac{1}{\dim(H)} \sum_{i=1}^k n_i \xi_i(e^{-1}) \chi_i
\qquad
e^{-1} = \frac{1}{\dim(H)} \sum_{i=1}^k n_i \xi_i(e) \chi_i
\end{align*}
\end{corollary}
\begin{pf}
Since $\chi_R$ is an integral, we can deduce from Proposition~\ref{IntDrinfDouble} that
$$\r_*(\chi_R) = \chi_R(\ud^{-1}) e^{-1}= \dim(H) e^{-1}$$
The second assertion therefore follows from the preceding proposition by applying~$\r_*$ to the equation $\chi_R = \sum_{i=1}^k n_i \chi_i$. The first assertion follows in a very similar way by applying~$\r_*^{-1}$, as we have
$\r_*^{-1}(\chi_R) = \dim(H) e$ by Proposition~\ref{IntDrinfDouble} and
$\r_*^{-1}(\chi_i) = \xi_i(e^{-1}) \chi_i$ by the preceding proposition.
\qed
\end{pf}

It should be pointed out in this context that these two elements are interchanged by~$\v_*$:
\begin{lemma}
$$\v_*(e) = e^{-1} \qquad \qquad \v_*(e^{-1}) = e$$
\end{lemma}
\begin{pf}
It follows from the definition that $\r_*(\ed) = \ed$.
Therefore, we get by Proposition~\ref{SecApprMod} that
$$\v_*(e) = (\t_*^{-1} \circ \r_*^{-1} \circ \t_*^{-1})(e)
= \t_*^{-1} (\r_*^{-1}(\ed)) = \t_*^{-1}(\ed)= e^{-1}$$
This proves the first assertion. The second assertion follows from the first by applying~$\v_*$, because we have $\v_*^2(\chi)=\sd^*(\chi)$ for all $\chi \in \Ch(D)$ by Proposition~\ref{RoleInt} and Corollary~\ref{InvSigma}, and we have seen in the proof of Lemma~\ref{IntDrinfDouble} that~$\sd^*(e)=e$.~\qed
\end{pf}

\newpage
\section{Induced modules} \label{Sec:IndMod}
\subsection[Induction]{} \label{Induct}
Suppose that $H$ is a semisimple Hopf algebra over an algebraically closed field~$K$ of characteristic zero, and consider its Drinfel'd double~$D=D(H)$. For an $H$-module~$V$, we can form the induced $D$-module:
$$D \o_H V = (H^* \o H) \o_H V \cong H^* \o V$$
where the last isomorphism maps~$\varphi \o h \o_H v$ 
to~$\varphi \o h.v$. This isomorphism is $D$-linear if we consider
$H^* \o V$ as a $D$-module via the module structure\endnote{\cite{YYY2}, Par.~6.4, p.~47.}
$$(\varphi \o h).(\varphi' \o v) := 
\varphi'_\1(\sh(h_\3)) \varphi'_\3(h_\1) \; \varphi \varphi'_\2 \o h_\2.v$$
We will view the induced module from this latter viewpoint in the sequel and therefore write $\Ind(V) := H^* \o V$, considered as a $D$-module with this module structure.

Suppose now that~$W$ is another $H$-module. We introduce the following map:
\begin{defn}
Suppose that $b_1,\ldots,b_n$ is a basis of~$H$ with dual basis $b^*_1,\ldots,b^*_n$. We define
\begin{align*}
\beta_{V,W}: \Ind(V \o W) \rightarrow \Ind(W \o V),~\varphi \o v \o w \mapsto \sum_{i=1}^n  \varphi b_i^* \o w \o b_i.v
\end{align*}
\end{defn}

Let us record some first properties of this map:
\begin{lemma}
$\beta_{V,W}$ is a $D$-linear isomorphism. The inverse is given by
$$\beta_{V,W}^{-1}(\varphi \o w \o v) = 
\sum_{i=1}^n  \varphi b_i^* \o \sh(b_i).v \o w$$
Furthermore, we have
$$(\beta_{W,V} \circ \beta_{V,W})(x) = \ud^{-1}.x$$
for all $x \in \Ind(V \o W)$.
\end{lemma}
\begin{pf}
To establish $D$-linearity, $\beta_{V,W}$ has to commute with elements of the form~$\varphi \o \H$ and elements of the form~$\eh \o h$. As it clearly commutes with elements of the first form, we can concentrate on elements of the second form. We have
\begin{align*}
&(\eh \o h).\beta_{V,W}(\varphi \o v \o w) = \\ 
&\sum_{i=1}^n  (\varphi_\1 b_{i\1}^*)(\sh(h_\4)) (\varphi_\3 b_{i\3}^*)(h_\1) \varphi_\2 b_{i\2}^* \o h_\2.w \o h_\3 b_i. v = \\ 
&\sum_{i=1}^n  
\varphi_\1(\sh(h_\6)) b_{i\1}^*(\sh(h_\5))
\varphi_\3(h_\1) b_{i\3}^*(h_\2) 
\varphi_\2 b_{i\2}^* \o h_\3.w \o h_\4 b_i. v
\end{align*}
Using the dual basis formulas stated in the introduction, this becomes
\begin{align*}
&(\eh \o h).\beta_{V,W}(\varphi \o v \o w) = \\ 
&\sum_{i_1,i_2,i_3=1}^n  
\varphi_\1(\sh(h_\6)) b_{i_1}^*(\sh(h_\5))
\varphi_\3(h_\1) b_{i_3}^*(h_\2) \\
&\mspace{300mu}
\varphi_\2 b_{i_2}^* \o h_\3.w \o h_\4 b_{i_1}b_{i_2}b_{i_3}. v = \\ 
&\sum_{i_2=1}^n \varphi_\1(\sh(h_\6)) \varphi_\3(h_\1) \;
\varphi_\2 b_{i_2}^* \o h_\3.w \o h_\4 \sh(h_\5) b_{i_2}h_\2. v= \\ 
&\sum_{i=1}^n \varphi_\1(\sh(h_\4)) \varphi_\3(h_\1) \;
\varphi_\2 b_{i}^* \o h_\3.w \o b_{i}h_\2. v
\end{align*}
But this is exactly $\beta_{V,W}((\eh \o h).(\varphi \o v \o w))$,
which establishes the $D$-linearity. To establish the form of the inverse, we note that 
\begin{align*}
\beta_{V,W}(&\sum_{i=1}^n \varphi b_i^* \o \sh(b_i).v \o w)
= \sum_{i,j=1}^n  \varphi b_i^* b_j^* \o w \o b_j\sh(b_i).v \\
&= \sum_{i,j=1}^n  \varphi b_i^*  \o w \o b_{i\2} \sh(b_{i\1}).v 
= \varphi \o w \o v 
\end{align*}
by the dual basis formulas, which establishes that the map stated is a
right inverse of~$\beta_{V,W}$. It can be shown similarly that it is also a left inverse.

To establish the last property, we can assume that
$x = \varphi \o v \o w$ is decomposable. As discussed in Paragraph~\ref{DrinfDouble}, 
$\ud^{-1} = \sum_{i=1}^n b_i^* \o b_i$ is central, since $H$ is involutory.
We therefore have
$$\ud^{-1}.x = \sum_{i=1}^n \varphi b_i^* \o b_{i\1}.v \o b_{i\2}.w
= \sum_{i,j=1}^n \varphi b_i^* b_j^* \o b_{i}.v \o b_{j}.w$$
by the dual basis formula. But this is 
exactly $(\beta_{W,V} \circ \beta_{V,W})(x)$.
\qed
\end{pf}

From the point of view of category theory, $\beta$ is a natural transformation between the functors $(V,W) \mapsto \Ind(V \o W)$ and $(V,W) \mapsto \Ind(W \o V)$. The natural transformation~$\beta$ also satisfies the following coherence properties:
\begin{prop}
If $U$, $V$, and $W$ are $H$-modules, the following diagram commutes:
$$\Vtriangle<1`1`-1;900>[\Ind(U \o V \o W)`\Ind(W \o U \o V)`
\Ind(V \o W \o U);
\beta_{U \o V,W}`\beta_{U, V \o W}`\beta_{V, W \o U}]$$

In addition, the following diagrams also commute:
$$\square<1`1`1`1;1000`700>[\Ind(V \o K)`
\Ind(K \o V)`\Ind(V)`\Ind(V);\beta_{V,K}`\cong`\cong`w \mapsto \ud^{-1}.w]
\qquad
\square<1`1`1`1;1000`700>[\Ind(K \o V)`
\Ind(V \o K)`\Ind(V)`\Ind(V);\beta_{K,V}`\cong`\cong`\id]$$
Here, the vertical maps are induced from the canonical isomorphisms.
\end{prop}
\begin{pf}
If $\varphi \in H^*$, $u \in U$, $v \in V$, and $w \in W$, we have 
\begin{align*}
(\beta_{V, W \o U} \circ \beta_{U, V \o W})&(\varphi \o u \o v \o w) =
\beta_{V, W \o U}(\sum_{i=1}^n \varphi b_i^* \o v \o w  \o b_i.u) \\
&= \sum_{i,j=1}^n \varphi b_i^* b_j^* \o w  \o b_i.u \o b_j.v
\end{align*}
on the one hand and 
\begin{align*}
\beta_{U \o V,W}(\varphi \o u \o v \o w)
= \sum_{i=1}^n \varphi b_i^* \o w  \o b_{i\1}.u \o b_{i\2}.v
\end{align*}
on the other hand. By the dual basis formula, both expressions agree,
proving $\beta_{U \o V,W} = \beta_{V, W \o U} \circ \beta_{U, V \o W}$, which establishes the commutativity of the first diagram. The commutativity of the two remaining diagrams follows directly from the definitions.
\qed
\end{pf}

\subsection[Induction and duality]{} \label{IndDual}
For a finite-dimensional module~$V$, it turns out that the induced module of the dual is isomorphic to the dual of the induced module. More generally, suppose that~$V$ and~$V'$ are two finite-dimensional $H$-modules endowed with a nondegenerate pairing
$\langle \cdot, \cdot \rangle: V \times V' \rightarrow K$ that satisfies
$$\langle h.v, v' \rangle = \langle v, \sh(h).v' \rangle$$
for all $v \in V$, $v' \in V'$, and $h \in H$.
If we then choose a nonzero integral~$\Lmh \in H$ and define a pairing 
$\langle \cdot, \cdot \rangle_\Lmh: \Ind(V) \times \Ind(V') 
\rightarrow K$ as
$$\langle \varphi \o v, \psi \o v' \rangle_\Lmh :=
(\sh^*(\varphi) \psi)(\Lmh) \; \langle v, v' \rangle$$
for $\varphi, \psi \in H^*$, $v \in V$, and $v' \in V'$,
this pairing has the following properties:
\begin{lemma}
$\langle \cdot, \cdot \rangle_\Lmh$ is nondegenerate. For $x \in D$, we have
$$\langle x.(\varphi \o v), \psi \o v' \rangle_\Lmh =
\langle \varphi \o v, \sd(x).(\psi \o v') \rangle_\Lmh$$
\end{lemma}
\begin{pf}
The nondegeneracy follows from the nondegeneracy of the pairing
$(\varphi,\psi) \mapsto (\sh^*(\varphi) \psi)(\Lmh)$.\endnote{\cite{M}, Thm.~2.1.3, p.~18.} To prove the second assertion, it suffices to show this in the cases $x= \varphi' \o \H$ and $x= \eh \o h$. In the first case,
this amounts to the identity
$$(\sh^*(\varphi' \varphi) \psi)(\Lmh) \; \langle v, v' \rangle
= (\sh^*(\varphi) \sh^*(\varphi') \psi)(\Lmh) \; \langle v, v' \rangle$$
In the second case, this amounts to the identity
\begin{align*}
&\varphi_\1(\sh(h_\3)) \; \varphi_\3(h_\1) \; 
(\sh^*(\varphi_\2) \psi)(\Lmh) \; \langle h_\2.v, v' \rangle \\
&= \psi_\1(h_\1) \; \psi_\3(\sh(h_\3)) \; (\sh^*(\varphi) \psi_\2)(\Lmh) \; \langle v, \sh(h_\2).v' \rangle
\end{align*}
which by the property of the original pairing will follow from
\begin{align*}
&\varphi_\1(\sh(h_\3)) \; \varphi_\3(h_\1) \; 
\varphi_\2(\sh(\Lmh_\1)) \psi(\Lmh_\2) \; h_\2 \\
&= \psi_\1(h_\1) \; \psi_\3(\sh(h_\3)) \; \varphi(\sh(\Lmh_\1)) \psi_\2(\Lmh_\2) \; h_\2
\end{align*}
This can be written as
\begin{align*}
&\varphi(\sh(h_\3) \sh(\Lmh_\1) h_\1) \;  \psi(\Lmh_\2) \; h_\2 
= \varphi(\sh(\Lmh_\1)) \; \psi(h_\1 \Lmh_\2 \sh(h_\3)) \;  h_\2
\end{align*}
which is a consequence of the fact that $\sh(\Lmh_\1) \o \Lmh_\2$
is a symmetric Casimir element.\endnote{\cite{LR2}, Lem.~1.2, p.~270.}
\qed
\end{pf}

Suppose now that~$W$ and~$W'$ is another pair of finite-dimensional $H$-modules endowed with another nondegenerate pairing
$\langle \cdot, \cdot \rangle: W \times W' \rightarrow K$ that satisfies
$$\langle h.w, w' \rangle = \langle w, \sh(h).w' \rangle$$
for all $w \in W$, $w' \in W'$, and $h \in H$. We can then form a pairing between the tensor products~$V \o W$ and~$W' \o V'$ that has the form
$$\langle v \o w, w' \o v' \rangle_\o := 
\langle v, v' \rangle \langle w, w' \rangle$$
This pairing is also nondegenerate and satisfies
$$\langle h.(v \o w), w' \o v' \rangle_\o =
\langle v \o w, \sh(h).(w' \o v') \rangle_\o$$

We can therefore invoke the preceding lemma to get a nondegenerate pairing~$\langle \cdot, \cdot \rangle_\Lmh$
between~$\Ind(V \o W)$ and~$\Ind(W' \o V')$ that has the explicit form
$$\langle \varphi \o v \o w, \psi \o w' \o v'\rangle_\Lmh =
(\sh^*(\varphi) \psi)(\Lmh) \; 
\langle v, v' \rangle \; \langle w, w' \rangle$$
Interchanging the roles of~$V$ and~$W$, we also get a pairing between~$\Ind(W \o V)$ and~$\Ind(V' \o W')$, for which we use the same notation and which is explicitly given as
$$\langle \varphi \o w \o v, \psi \o v' \o w'\rangle_\Lmh =
(\sh^*(\varphi) \psi)(\Lmh) \; \langle w, w' \rangle \; 
\langle v, v' \rangle $$
These pairings are compatible with the morphisms introduced in Paragraph~\ref{Induct} in the following way:
\begin{prop}
$$\langle \beta_{V,W}(\varphi \o v \o w), \psi \o v' \o w'\rangle_\Lmh =
\langle \varphi \o v \o w, \beta_{V',W'}(\psi \o v' \o w') \rangle_\Lmh$$
\end{prop}
\begin{pf}
On the one hand, we have
\begin{align*}
\langle \beta_{V,W}(\varphi \o v \o w)&, \psi \o v' \o w'\rangle_\Lmh 
= \sum_{i=1}^n 
\langle \varphi b_i^* \o w \o b_i.v, \psi \o v' \o w'\rangle_\Lmh \\
&= \sum_{i=1}^n (\sh^*(\varphi b_i^*) \psi) (\Lmh)  
\langle w, w' \rangle \langle b_i.v, v' \rangle \\
&= \sum_{i=1}^n (\sh^*(b_i^*) \sh^*(\varphi) \psi) (\Lmh) 
\langle w, w' \rangle \langle v, \sh(b_i).v' \rangle \\
&= \sum_{i=1}^n (b_i^* \sh^*(\varphi) \psi) (\Lmh) 
\langle w, w' \rangle \langle v, b_i.v' \rangle 
\end{align*}
On the other hand, we have
\begin{align*}
\langle \varphi \o v \o w, \beta_{V',W'}(\psi \o v' \o w') \rangle_\Lmh
&= \sum_{i=1}^n 
\langle \varphi \o v \o w, \psi b_i^* \o w' \o b_i.v' \rangle_\Lmh \\
&= \sum_{i=1}^n 
(\sh^*(\varphi) \psi b_i^*) (\Lmh) \langle v, b_i.v' \rangle \langle w, w' \rangle
\end{align*}
Both expressions are equal because~$\Lmh$ is cocommutative.\endnote{\cite{LR2}, Prop.~5.4, p.~282.}
\qed
\end{pf}
It is of course possible to choose the dual~$V^*$ for~$V'$ and the dual~$W^*$ for~$W'$. The above discussion then shows that
$\Ind(V \o W)^* \cong \Ind(W^* \o V^*)$ and also
$\Ind(W \o V)^* \cong \Ind(V^* \o W^*)$.
Using these identifications, it follows from the above proposition that
we have
$$\beta_{V,W}^* = \beta_{V^*,W^*}$$
for the transpose of~$\beta_{V,W}$.

\subsection[The relation with the center construction]{} \label{RelCent}
As pointed out by D.~Nikshych,\endnote{Private communication, Chicago, 2007.}
the natural transformation~$\beta$ introduced in Definition~\ref{Induct} can be related to the categorical center construction. Recall\endnote{\cite{Kas}, Thm.~XIII.5.1, p.~333.}
that the category of modules over the Drinfel'd double~$D=D(H)$ can be considered as the center of the category of $H$-modules. This implies in particular that for every $D$-module~$U$ and every $H$-module~$V$ we have the isomorphism
$$c_{V,U}: V \o U \rightarrow U \o V,~v \o u \mapsto 
\sum_{i=1}^n (b_i^* \o \H).u \o b_i.v$$
which consists in the application of the R-matrix followed by interchanging the tensorands. Its inverse is therefore given by
$$c_{V,U}^{-1}: U \o V \rightarrow V \o U,~u \o v \mapsto 
\sum_{i=1}^n  \sh(b_i).v \o (b_i^* \o \H).u$$
To relate this map to the isomorphism~$\beta_{V,W}$, where~$W$ is another
$H$-module, note that
$$\Hom_D(\Ind(V \o W),U) \cong \Hom_H(V \o W,U)$$
by the Frobenius reciprocity theorem.\endnote{\cite{Lang}, Chap.~XVIII, \S~7, p.~689.} Therefore, there is a unique isomorphism~$\beta'_{V,W;U}$
that makes the diagram 
$$\square<1`1`1`1;1600`700>[\Hom_D(\Ind(W \o V),U)`
\Hom_D(\Ind(V \o W),U)`\Hom_H(W \o V,U)`
\Hom_H(V \o W,U);\circ \beta_{V,W}```\beta'_{V,W;U}]$$
commutative, where~$\circ \beta_{V,W}$ denotes the map coming from~$\beta_{V,W}$ by composition on the right.
The isomorphism~$\beta'_{V,W;U}$ is given explicitly as
$$\beta'_{V,W;U}(f) (v \o w) = \sum_{i=1}^n (b_i^* \o \H).f(w \o b_i.v)$$
for $f \in \Hom_H(W \o V,U)$, as we have
\begin{align*}
(&g \circ \beta_{V,W})(\eh \o v \o w)
= g(\sum_{i=1}^n b_i^* \o w \o b_i.v) 
= \sum_{i=1}^n (b_i^*\o \H). g(\eh \o w \o b_i.v)
\end{align*}
for $g \in \Hom_D(\Ind(W \o V),U)$.

Besides the adjunction between induction and restriction that appears in the Frobenius reciprocity theorem, there are two other pairs of adjoint functors that appear in this setting: The composition\endnote{\cite{Kas}, Prop.~XIV.2.2b, p.~343f.}
$$\Hom_K(W, U \o V^*) \rightarrow \Hom_K(W \o V, U \o V^* \o V) 
\rightarrow \Hom_K(W \o V, U)$$
where the first map takes~$f$ to~$f \o \id_V$ and the second evaluates~$V^*$ on~$V$, defines a homomorphism from $\Hom_K(W, U \o V^*)$ to~$\Hom_K(W \o V, U)$. The image~$g \in \Hom_K(W \o V, U)$ of~$f \in \Hom_K(W, U \o V^*)$ is given explicitly as
$$g(w \o v) = \sum_j \psi_j(v) u_j$$
if $f(w) = \sum_j u_j \o \psi_j$. This composition is bijective if~$V$ is finite-dimensional. Because the evaluation map 
$V^* \o V \rightarrow K,~\psi \o v \mapsto \psi(v)$
is $H$-linear, both mappings that appear in the composition preserve
the subspace of $H$-linear maps, so that we can restrict this composition to a map from $\Hom_H(W, U \o V^*)$ to~$\Hom_H(W \o V, U)$, which is an isomorphism in the finite-dimensional case.

Similarly, the composition
$$\Hom_K(W, V^* \o U) \rightarrow \Hom_K(V \o W, V \o V^* \o U) 
\rightarrow \Hom_K(V \o W, U)$$
obtained by tensoring with~$\id_V$ on the left and then evaluating~$V^*$ on~$V$ leads to a homomorphism that takes a 
linear map~$f \in \Hom_K(W, V^* \o U)$ to the map 
$g \in \Hom_K(V \o W, U)$ that satisfies
$$g(v \o w) = \sum_j \psi_j(v) u_j$$
if $f(w) = \sum_j \psi_j \o u_j$. This time, the second homomorphism in the composition uses the evaluation map
$$V \o V^* \rightarrow K, v \o \psi \mapsto \psi(v)$$
But as the antipode is an involution, this evaluation map
is also $H$-linear, so that we again get a homomorphism 
from~$\Hom_H(W, V^* \o U)$ to~$\Hom_H(V \o W, U)$
by restriction, which is an isomorphism if~$V$ is finite-dimensional.

All these mappings come together in the following proposition:
\begin{prop}
The diagram
$$\square<1`-1`-1`1;1600`700>[\Hom_H(W \o V,U)`\Hom_H(V \o W,U)`\Hom_H(W,U \o V^*)`
\Hom_H(W, V^* \o U);\beta'_{V,W;U}```c_{V^*,U}^{-1} \circ ]$$
is commutative, where $c_{V^*,U}^{-1} \circ$ denotes composition with~$c_{V^*,U}^{-1}$ on the left.
\end{prop}
\begin{pf}
Suppose that $f \in \Hom_H(W,U \o V^*)$. The two possible
paths in the diagram give two elements of~$\Hom_H(V \o W,U)$,
and we have to prove that they are equal. For this, it suffices to show
that they agree for every decomposable tensor $v \o w \in V \o W$.
To see this, write $f(w) = \sum_j u_j \o \psi_j$. Then we have
$$c_{V^*,U}^{-1}(f(w)) = 
\sum_j \sum_{i=1}^n \sh(b_i).\psi_j \o (b_i^* \o \H).u_j$$
This means that the homomorphism that arises from composing the lower
and the right arrow maps our decomposable tensor~$v \o w$ to
$$\sum_j \sum_{i=1}^n (\sh(b_i).\psi_j)(v) \; (b_i^* \o \H).u_j
= \sum_j \sum_{i=1}^n \psi_j(b_i.v) \; (b_i^* \o \H).u_j$$
where we have used that the antipode is an involution.

On the other hand, if $g \in \Hom_H(W \o V, U)$ is the image of~$f$ under
the left arrow, then we have
$g(w \o v) = \sum_j \psi_j(v) u_j$,
so that
$$\beta'_{V,W;U}(g)(v \o w)= \sum_{i=1}^n (b_i^* \o \H).g(w \o b_i.v)
= \sum_j \sum_{i=1}^n \psi_j(b_i.v) \; (b_i^* \o \H).u_j$$
which is exactly the result coming from the other path.
\qed
\end{pf}

If we insert the definition of~$\beta'_{V,W;U}$ in this diagram, we can 
extend the vertical arrows by the isomorphisms coming from the Frobenius reciprocity theorem to get the diagram 
$$\square<1`-1`-1`1;1600`700>[\Hom_D(\Ind(W \o V),U)`\Hom_D(\Ind(V \o W),U)`\Hom_H(W,U \o V^*)`
\Hom_H(W, V^* \o U);\circ \beta_{V,W}```c_{V^*,U}^{-1} \circ ]$$
which exhibits the relation between the natural transformation~$\beta$
and the categorical center construction.

\subsection[The relation of the coherence properties]{} \label{RelCoherProp}
The coherence properties of the natural transformation~$\beta$ stated in  Proposition~\ref{Induct} can also be related to the coherence properties of the natural transformation~$c$ required in the categorical center construction.\endnote{\cite{Kas}, Def.~XIII.4.1, p.~330.} If~$V_1$ and~$V_2$ are $H$-modules and 
$U$ is a $D$-module, then the natural transformation~$c$ makes
the triangle 
$$\Vtriangle<1`1`-1;900>[
U \o V_2^* \o V_1^*`
V_2^* \o V_1^* \o U`
V_2^* \o U \o V_1^*;
c^{-1}_{V_2^* \o V_1^*,U}`
c^{-1}_{V_2^*,U} \o \id_{V_1^*}`
\id_{V_2^*} \o c^{-1}_{V_1^*,U}]$$
commutative. If $W$ is another $H$-module, we therefore have that the diagram
$$\Vtriangle<1`1`-1;900>[\Hom_H(W,U \o V_2^* \o V_1^*)`
\Hom_H(W,V_2^* \o V_1^* \o U)`
\Hom_H(W,V_2^* \o U \o V_1^*);
c^{-1}_{V_2^* \o V_1^*,U} \circ`
(c^{-1}_{V_2^*,U} \o \id_{V_1^*}) \circ`
(\id_{V_2^*} \o c^{-1}_{V_1^*,U}) \circ]$$
also commutes, where we have used the circle notation from Paragraph~\ref{RelCent}. To translate this diagram into a diagram for the natural transformation~$\beta$, we need as an intermediate step on the left side the diagram
$$
\bfig
\putsquare<1`1`1`1;2000`600>(700,600)[\Hom_H(W,U \o V_2^* \o V_1^*)`
\Hom_H(W,V_2^* \o U \o V_1^*)`
\Hom_H(W \o V_1,U \o V_2^*)`
\Hom_H(W \o V_1,V_2^* \o U);
(c^{-1}_{V_2^*,U} \o \id_{V_1^*}) \circ```
c^{-1}_{V_2^*,U} \circ]
\putmorphism(700,600)(0,-1)[`\Hom_H(W \o V_1 \o V_2,U)`]{600}1m
\putmorphism(2700,600)(0,-1)[`\Hom_H(V_2 \o W \o V_1,U)`]{600}1m
\putmorphism(1080,0)(1,0)[\phantom{T'TT'}`
  \phantom{T'T'}`\beta'_{V_2,W \o V_1;U}]{1200}1b
\efig
$$
where the first diagram commutes because of the naturality of the adjunction and the second diagram commutes by Proposition~\ref{RelCent}. Similarly, we need on the right side the diagram
$$
\bfig
\putsquare<1`1`1`1;2000`600>(700,600)[
\Hom_H(W,V_2^* \o U \o V_1^*)`
\Hom_H(W,V_2^* \o V_1^* \o U)`
\Hom_H(V_2 \o W,U \o V_1^*)`
\Hom_H(V_2 \o W,V_1^* \o U);
(\id_{V_2^*} \o c^{-1}_{V_1^*,U}) \circ```
c^{-1}_{V_1^*,U} \circ]
\putmorphism(700,600)(0,-1)[`\Hom_H(V_2 \o W \o V_1,U)`]{600}1m
\putmorphism(2700,600)(0,-1)[`\Hom_H(V_1 \o V_2 \o W,U)`]{600}1m
\putmorphism(1080,0)(1,0)[\phantom{T'TT'}`
\phantom{T'T'}`\beta'_{V_1, V_2 \o W;U}]{1200}1b
\efig
$$
which commutes for exactly the same reasons. Using this, the commuting triangle for~$c$ above translates into the diagram
$$\Vtriangle<1`1`-1;900>[\Hom_H(W \o V_1 \o V_2,U)`
\Hom_H(V_1 \o V_2 \o W,U)`
\Hom_H(V_2 \o W \o V_1,U);
\beta'_{V_1 \o V_2, W;U}`
\beta'_{V_2,W \o V_1;U}`
\beta'_{V_1, V_2 \o W;U}]$$
where the map at the top has been translated by Proposition~\ref{RelCent},
using the fact that $(V_1 \o V_2)^* \cong V_2^* \o V_1^*$ in a way that 
is compatible with the translation.

Using the adjunction between induction and restriction again, we can translate the last triangle further into the triangle
$$\Vtriangle<1`1`-1;900>[\Hom_D(\Ind(W \o V_1 \o V_2),U)`
\Hom_D(\Ind(V_1 \o V_2 \o W),U)`
\Hom_D(\Ind(V_2 \o W \o V_1),U);
\circ \beta_{V_1 \o V_2, W}`
\circ \beta_{V_2,W \o V_1}`
\circ \beta_{V_1, V_2 \o W}]$$
which in turn by the Yoneda lemma implies the first coherence property of~$\beta$ as given in Proposition~\ref{Induct}. Note that, although all the above diagrams commute in any case, this amounts to a new proof of
first coherence property of~$\beta$ only in the case where~$V_1$ and~$V_2$ are finite-dimensional, because otherwise the commutativity of the second triangle above does not logically imply the commutativity of the third triangle.

\subsection[Adjoint functors]{} \label{AdjFunct}
To analyse the relation between~$\beta$ and~$c$ further, we need some preparation. So far in this section, we have basically used two pairs of adjoint functors: The adjunction between induction and restriction and the adjunction between tensoring with a module and tensoring with its dual. These two adjunctions can be related by the following map:
\begin{lemma}
For an $H$-module~$V$ and a $D$-module~$W$, the map
$$r_{V,W}: \Ind(V \o W) \rightarrow \Ind(V) \o W,~x \o_H (v \o w) \mapsto
(x_\1 \o_H v) \o x_\2.w$$
is a $D$-linear isomorphism.
\end{lemma}
\begin{pf}
First, $r_{V,W}$ is well-defined, as we have
\begin{align*}
r_{V,W}(&xh \o_H (v \o w)) = (x_\1 h_\1 \o_H v) \o x_\2 h_\2.w \\
&= (x_\1 \o_H h_\1 v) \o x_\2 h_\2.w
= r_{V,W}(x \o_H (h_\1.v \o h_\2.w))
\end{align*}
The reader should verify at this point that $r_{V,W}$ is $D$-linear.
The bijectivity will follow if we can show that the potential inverse 
$$r_{V,W}^{-1}: \Ind(V) \o W \rightarrow \Ind(V \o W),~(x \o_H v) \o w \mapsto x_\1 \o_H (v \o \sa(x_\2).w)$$
is also well-defined. This holds because
\begin{align*}
&r_{V,W}^{-1}((xh \o_H v) \o w) = 
x_\1 h_\1\o_H (v \o \sa(x_\2 h_\2).w) \\
&= x_\1 \o_H (h_\1.v \o h_\2 \sa(h_\3) \sa(x_\2).w)
= x_\1 \o_H (h.v \o \sa(x_\2).w)
\end{align*}
establishing the assertion.
\qed
\end{pf}
Note that, under the correspondence of both $\Ind(V \o W)$ and
$\Ind(V) \o W$ to $H^* \o V \o W$ introduced in Paragraph~\ref{Induct}, the map $r_{V,W}$ becomes
$$r_{V,W}: H^* \o V \o W \rightarrow H^* \o V \o W,~\varphi \o v \o w
\mapsto \varphi_\2 \o v \o (\varphi_\1 \o \H).w$$

Our claim that~$r_{V,W}$ relates the two adjunctions is justified 
by the following proposition:
\begin{prop}
For an $H$-module~$V$ and $D$-modules~$W$ and~$U$, the diagram
$$
\bfig
\putsquare<1`-1`-1`0;1800`400>(100,800)[\Hom_D(\Ind(V \o W),U)`
\Hom_H(V \o W,U)`\Hom_D(\Ind(V) \o W,U)`\Hom_H(V,U \o W^*);`
\circ r_{V,W}``]
\putVtriangle<0`-1`-1;800>(200,000)[``\Hom_D(\Ind(V),U \o W^*);``]
\efig
$$
commutes.
\end{prop}
\begin{pf}
Suppose that $f \in \Hom_D(\Ind(V),U \o W^*)$, and consider the image
$g \in \Hom_H(V \o W,U)$ of~$f$ under the left path. Since
$$r_{V,W}(\D \o_H (v \o w)) = (\D \o_H v) \o w$$
we have
$g(v \o w) = \sum_j \psi_j(w_j) u_j$
if $f(\D \o_H v) = \sum_j u_j \o \psi_j$.
But this is exactly the map that arises as the image of~$f$ under the right path.
\qed
\end{pf}

\subsection[More coherence properties]{} \label{MoreCohProp}
As a comparison shows, the coherence condition stated at the beginning
of Paragraph~\ref{RelCoherProp} corresponds to only one of the two conditions that appear in the definition of a quasisymmetry.\endnote{\cite{Kas}, Def.~XIII.1.1, p.~315.}
From the point of view of the center construction, the second condition enters into the definition of the tensor product of two objects. We therefore expect that there is another relation between~$\beta$ and~$c$ that can be deduced from this second condition by arguing as in Paragraph~\ref{RelCoherProp}. Before we state this relation, we recall the relation between braiding and duality:
\begin{lemma}
For an $H$-module~$V$ and a $D$-module~$U$, the diagram

$$
\bfig
\putsquare<1`1`1`1;1500`1000>(100,100)[V^* \o U^*`
(U \o V)^*`U^* \o V^*`(V \o U)^*;`c_{V^*,U^*}`\circ c_{V,U}`]
\efig
$$
commutes.
\end{lemma}
\begin{pf}
Recall\endnote{\cite{Kas}, Prop.~XIV.2.2c, p.~343f.} that the 
top horizontal arrow maps $\varphi \o \psi \in V^* \o U^*$ to the 
linear form $u \o v \mapsto \varphi(v) \psi(u)$, and the horizontal arrow at the bottom is defined similarly. We therefore have for 
$\varphi \in V^*$, $\psi \in U^*$, $v \in V$, and $u \in U$ that
\begin{align*}
c_{V^*,U^*}(\varphi \o \psi)(v \o u) 
&= \sum_{i=1}^n ((b_i^* \o \H).\psi \o b_i.\varphi)(v \o u) \\
&= \sum_{i=1}^n \psi(\sd(b_i^* \o \H).u) \; \varphi(\sh(b_i).v) \\
&= \sum_{i=1}^n \psi((b_i^* \o \H).u) \; \varphi(b_i.v) 
= (\varphi \o \psi)(c_{V,U}(v \o u)) 
\end{align*}
where we have used the notation from Paragraph~\ref{DrinfDouble} resp.~Paragraph~\ref{RelCoherProp} for the R-matrix.
\qed
\end{pf}

With the help of this lemma, we now derive the following additional
relation between~$\beta$ and~$c$:
\begin{prop}
Suppose that~$V$ and~$W$ are $H$-modules and that~$U$ is a $D$-module.
We assume that~$V$ and~$U$ are finite-dimensional. Then the diagram

$$
\bfig
\putsquare<1`1`-1`0;1800`400>(100,800)[\Ind(V \o W \o U)`
\Ind(W \o U \o V)`\Ind(V \o W) \o U`\Ind(W \o V \o U);\beta_{V, W \o U}`r_{V \o W,U}`\Ind(\id_W \o c_{V,U})`]
\putVtriangle<0`1`-1;800>(200,0)[``\Ind(W\o V) \o U;`
\beta_{V,W} \o \id_U`r_{W \o V,U}^{-1}]
\efig
$$
commutes.
\end{prop}
\begin{pf}
By the Yoneda lemma, it suffices to prove the commutativity of the diagram after the application of the contravariant functor~$\Hom_D(-,X)$, where~$X$ is another $D$-module. After this application, the diagram takes the form
$$
\bfig
\putsquare<-1`-1`1`0;1800`400>(100,800)[\Hom_D(\Ind(V \o W \o U),X)`
\Hom_D(\Ind(W \o U \o V),X)`\Hom_D(\Ind(V \o W) \o U,X)`
\Hom_D(\Ind(W \o V \o U),X);\circ \beta_{V, W \o U}`
\circ r_{V \o W,U}`\circ\Ind(\id_W \o c_{V,U})`]
\putVtriangle<0`-1`-1;800>(200,0)[``\Hom_D(\Ind(W\o V) \o U,X);`
\circ (\beta_{V,W} \o \id_U)` \circ r_{W \o V,U}]
\efig
$$
Using the defining property of the maps~$\beta'$ from Paragraph~\ref{RelCent} together with Proposition~\ref{AdjFunct}, we see that the commutativity of this diagram follows from the commutativity of the diagram
$$
\bfig
\putsquare<-1`-1`1`0;1800`400>(100,800)[\Hom_H(V \o W \o U,X)`
\Hom_H(W \o U \o V,X)`\Hom_H(V \o W ,X \o U^*)`
\Hom_H(W \o V \o U,X);\beta'_{V, W \o U;X}``\circ(\id_W \o c_{V,U})`]
\putVtriangle<0`-1`-1;800>(200,0)[``\Hom_H(W\o V,X \o U^*);`
\beta'_{V,W;X \o U^*}`]
\efig
$$
By taking~$V$ and~$U$ to the right side in this diagram, we get from  Proposition~\ref{RelCent} that our assertion is equivalent to the
commutativity of
$$
\bfig
\putsquare<-1`-1`1`0;1800`400>(100,800)[\Hom_H(W,V^* \o X \o U^*)`
\Hom_H(W,X \o V^* \o U^*)`\Hom_H(W ,V^* \o X \o U^*)`
\Hom_H(W,X \o U^* \o V^*);
(c^{-1}_{V^*,X} \o \id_{U^*}) \circ`\id`(\id_X \o c_{V^*,U^*})\circ`]
\putVtriangle<0`-1`-1;800>(200,0)[``\Hom_H(W,X \o U^*\o V^*);`
c^{-1}_{V^*,X \o U^*} \circ`\id]
\efig
$$
where we have used in addition the preceding lemma for the right vertical arrow. But this is now by the Yoneda lemma equivalent to the equation
$$c_{V^*,X \o U^*} = 
(\id_X \o c_{V^*,U^*}) \circ(c_{V^*,X} \o \id_{U^*})$$
that is used to define the tensor product in the center construction.\endnote{\cite{Kas}, Thm.~XIII.4.2, Eq.~(4.4), p.~330.}
\qed
\end{pf}

It is of course also possible to prove this proposition by direct computation: For $\varphi \in H^*$, $v \in V$, $w \in W$, and~$u \in U$,
we have on the one hand
\begin{align*}
&(r_{W \o V,U}^{-1} \circ (\beta_{V,W} \o \id_U) \circ r_{V \o W,U})
(\varphi \o v \o w \o u) \\
&= (r_{W \o V,U}^{-1} \circ (\beta_{V,W} \o \id_U))
(\varphi_\2 \o v \o w \o (\varphi_\1 \o \H).u) \\
&= \sum_{i=1}^n r_{W \o V,U}^{-1}
(\varphi_\2 b^*_i \o w \o b_i.v \o (\varphi_\1 \o \H).u) \\
&= \sum_{i=1}^n 
\varphi_\3 b^*_{i\2} \o w \o b_i.v \o 
(\sh^{-1*}(\varphi_\2 b^*_{i\1}) \o \H)(\varphi_\1 \o \H).u \\
&= \sum_{i=1}^n 
\varphi b^*_{i\2} \o w \o b_i.v \o 
(\sh^{-1*}(b^*_{i\1}) \o \H).u
\end{align*}
and on the other hand
\begin{align*}
&(\Ind(\id_W \o c^{-1}_{V,U}) \circ (\beta_{V,W \o U}))
(\varphi \o v \o w \o u) \\
&= \sum_{i=1}^n \Ind(\id_W \o c^{-1}_{V,U})
(\varphi b^*_{i} \o w \o u \o b_i.v)\\
&= \sum_{i,j=1}^n \varphi b^*_{i} \o w  \o \sh(b_j)b_i.v 
\o (b^*_{j} \o \H).u
\end{align*}
The assertion therefore would follow from the equation
\begin{align*}
\sum_{i=1}^n b^*_{i\2} \o b_i \o \sh^{-1*}(b^*_{i\1}) 
= \sum_{i,j=1}^n b^*_{i} \o \sh(b_j)b_i \o b^*_{j}
\end{align*}
But as we have
\begin{align*}
\sum_{i=1}^n b^*_{i\2}(h) \; b_i \; \sh^{-1*}(b^*_{i\1})(h')
&= \sum_{i=1}^n b^*_{i}(\sh^{-1}(h')h) \; b_i \\
&= \sh^{-1}(h')h
= \sum_{i,j=1}^n b^*_{i}(h) \; \sh(b_j)b_i \; b^*_{j}(h')
\end{align*}
this equation holds.

Besides being substantially simpler, the proof by direct computation also shows that the requirement that~$V$ and~$U$ be finite-dimensional is unnecessary. We have nonetheless chosen to give the proof above because it
exhibits the relation to the second condition in the definition of a quasisymmetry. Note that these conditions also correspond to the equations
for~$(\da \o \id)(R)$ and~$(\id \o \da)(R)$ that appear in the definition of a quasitriangular Hopf algebra stated in Paragraph~\ref{QuasitriHopf}.

\newpage
\section{Equivariant Frobenius-Schur indicators} \label{Sec:EquiFrobSchur}
\subsection[Equivariant Frobenius-Schur indicators]{} \label{EquiFrobSchur}
We continue to work in the setting of Section~\ref{Sec:IndMod}, which was described in Paragraph~\ref{Induct}. So, $H$ is a semisimple Hopf algebra over an algebraically closed field~$K$ of characteristic zero, and $D=D(H)$ is its Drinfel'd double. For a finite-dimensional $H$-module~$V$ and a positive integer~$m$, we can of course form the \mbox{$m$-th} tensor power~$V^{\o m}$ of~$V$, and the Drinfel'd double~$D$ acts on its induced module~$\Ind(V^{\o m})$. We denote the corresponding representation by
$$\rho_m: D \rightarrow \End(\Ind(V^{\o m}))$$
A further endomorphism of~$\Ind(V^{\o m})$ is 
$\beta := \beta_{V,V^{\o (m-1)}}$. Using these ingredients, we can now define the following quantities:
\begin{defn}
For integers $m,l \in \Z$ with $m>1$ and a central element $z \in Z(D)$, we define the $(m,l)$-th equivariant Frobenius-Schur indicator of~$V$ and~$z$ as
$$I_V((m,l),z) := \Tr(\beta^l \circ \rho_m(z))$$
We extend this definition to all integers~$m$ as follows:
If $m=1$, we define the indicator by setting
$I_V((m,l),z) := \Tr(\rho_1(\ud^{-l}) \circ \rho_1(z))$.
If $m=0$, we write $z = \sum_j \varphi_j \o h_j$, and define for $l>0$
$$I_V((0,l),z) :=
\dim(H) \sum_j \eh(h_j) \varphi_{j}(\Lmh_\1)
\chi_{V^{\o l}}(\Lmh_\2)$$
where $\Lmh \in H$ is the integral that satisfies $\eh(\Lmh)=1$.
For $l=0$, we define 
$I_V((0,0),z) := \dim(H) \sum_j \eh(h_j) \varphi_{j}(\Lmh)$, whereas  we define
$$I_V((0,l),z) := I_V((0,-l),\sd(z))$$
for $l<0$. In the last case where $m<0$, we similarly define
$$I_V((m,l),z) = I_V((-m,-l),\sd(z))$$
\end{defn}

In the main case where $m>1$, it should be noted that we have
$\beta^l = \beta_{V^{\o l},V^{\o (m-l)}}$
for $l=1,2,\ldots,m-1$. This follows inductively from the coherence property given in Proposition~\ref{Induct}, because, if we set 
$U=V^{\o l}$ and $W=V^{\o (m-l-1)}$ there, we obtain the equation
$\beta_{V^{\o (l+1)},V^{\o (m-l-1)}} = 
\beta_{V, V^{\o (m-1)}} \circ \beta_{V^{\o l}, V^{\o (m-l)}}$.
If we interpret the $0$-th tensor power as the trivial module 
$K \cong V^{\o 0}$, then this formula also extends to the cases
$l=0$ and~$l=m$, because $\beta^0 = \id$ corresponds 
to~$\beta_{K,V^{\o m}}$ by Proposition~\ref{Induct}, and
$$\beta^m(x) = \beta(\beta^{m-1}(x)) = 
\beta_{V,V^{\o (m-1)}}(\beta_{V^{\o (m-1)},V}(x)) = \ud^{-1}.x$$
by Lemma~\ref{Induct}, which corresponds by Proposition~\ref{Induct}
to~$\beta_{V^{\o m},K}$. From this viewpoint, the case $m=1$ can also be subsumed under the case~$m>1$, because then $\beta$ coincides with the action of~$\ud^{-1}$. It should also be noted that the
formula $I_V((-m,-l),z) = I_V((m,l),\sd(z))$ holds for all integers~$m$ and~$l$ by definition; if $m=l=0$, one needs Lemma~\ref{RoleEval} in addition to see this.

An easy consequence of this definition is the following formula for the indicators of a tensor power:
\begin{lemma}
$I_{V^{\o q}}((m,l),z) = I_{V}((qm,ql),z)$
\end{lemma}
\begin{pf}
It is understood here that $q>0$ is a natural number. We consider the case 
$m>1$ first. As just explained, the $q$-th power of 
$\beta := \beta_{V,V^{\o (qm-1)}}$ is
$\beta^q = \beta_{V^{\o q},V^{\o (qm-q)}}$, so that
$$I_{V^{\o q}}((m,l),z) = 
\Tr(\beta^{ql} \circ \rho_{qm}(z)) 
= I_{V}((qm,ql),z)$$
The formula also holds in the case $m=1$ by the explanations above, and in the cases~$m=0$ and~$m<0$ it follows directly from the definitions.
\qed
\end{pf}

Part of the name given above to the quantities $I_V((m,l),z)$ is explained
by the following proposition, which relates them to the higher Frobenius-Schur indicators:\endnote{\cite{YYY2}, Def.~2.3, p.~15.}
\begin{prop}
Suppose that $\Lmh \in H$ and $\lmh \in H^*$ are integrals that are normalized so that $\eh(\Lmh)=\lmh(\H)=1$, and set $\Lmd := \lmh \o \Lmh$. If $\chi_V$ denotes the character of the $H$-module~$V$, then we have for its $m$-th Frobenius-Schur indicator~$\nu_m(\chi_V)$ that
$$\nu_m(\chi_V) = I_V((m,1),\Lmd)$$
for all integers~$m>0$.
\end{prop}
\begin{pf}
We treat the case $m=1$ separately. We have seen in Paragraph~\ref{IntDrinfDouble} that~$\Lmd$ is an integral of~$D$; however, the normalization here is different from the one in Paragraph~\ref{RoleEval}. By definition, we therefore have
$$I_V((1,1),\Lmd) = \Tr(\rho_1(\ud^{-1}) \circ \rho_1(\Lmd)) = \Tr(\rho_1(\Lmd))$$
Now we have $\rho_1(\Lmd)(\varphi \o v) = \varphi(\H) \lmh \o \Lmh.v$, so that $\Tr(\rho_1(\Lmd))=\chi_V(\Lmh)$, which is the assertion.

In the case $m>1$, note that the map
$$(V^{\o m})^H \rightarrow \Ind(V^{\o m})^D,~w \mapsto \lmh \o w$$
is an isomorphism between the spaces of invariants,\endnote{\cite{M}, \S~1.7, Def.~1.7.1, p.~13.} because~$\rho_m(\Lmd)$ is a projection to
$\Ind(V^{\o m})^D$ and we have
$$\rho_m(\Lmd)(\varphi \o v_1 \o \ldots \o v_m) = 
\lmh \o \Lmh.(\varphi(\H) v_1 \o \ldots \o v_m)$$
Because $\beta$ is $D$-linear, it commutes with~$\rho_m(\Lmd)$, and therefore preserves the space~$\Ind(V^{\o m})^D$ of invariants.
Similarly, the map
$$\alpha: V^{\o m} \rightarrow V^{\o m},~v_1 \o \ldots \o v_m 
\mapsto v_2 \o v_3 \o \ldots \o v_m \o v_1$$
preserves\endnote{\cite{YYY2}, Par.~2.3, p.~17.} the space~$(V^{\o m})^H$, and the diagram 
$$\square<1`1`1`1;1000`700>[(V^{\o m})^H`(V^{\o m})^H`\Ind(V^{\o m})^D`\Ind(V^{\o m})^D;\alpha`w \mapsto \lmh \o w`w \mapsto \lmh \o w`\beta]$$
is commutative, since we have
\begin{align*}
&\beta(\lmh \o v_1 \o v_2 \o \ldots \o v_m) = 
\sum_{i=1}^n \lmh b^*_i \o v_2 \o \ldots \o v_m \o b_i.v_1 \\
&= \lmh \o v_2 \o \ldots \o v_m \o v_1 =
\lmh \o \alpha(v_1 \o v_2 \o \ldots \o v_m) 
\end{align*}
Since the restriction of~$\beta$ to~$\Ind(V^{\o m})^D$ is therefore conjugate to the restriction of~$\alpha$ to~$(V^{\o m})^H$, the traces of these two maps have to coincide, which yields
$$I_V((m,1),\Lmd) = \Tr(\beta \circ \rho_m(\Lmd)) = \Tr(\beta \mid_{\Ind(V^{\o m})^D}) 
= \Tr(\alpha \mid_{(V^{\o m})^H}) = \nu_m(\chi_V)$$
by the first formula for the Frobenius-Schur indicators.\endnote{\cite{YYY2}, Cor.~2.3, p.~17.} 
\qed
\end{pf}

It should be noted that the normalization for the integral of the Drinfel'd double in the preceding proposition is different from the one used in Paragraph~\ref{RoleEval}; we have chosen the normalization in the proposition to avoid the appearance of another proportionality factor.
Furthermore, it should be noted that, as a consequence of the preceding argument, the restriction of a power~$\beta^l$ to~$\Ind(V^{\o m})^D$ is also conjugate to the restriction of the corresponding power~$\alpha^l$ 
to~$(V^{\o m})^H$, so that we get
$$I_V((m,l),\Lmd) = \Tr(\alpha^l \mid_{(V^{\o m})^H})$$
for all $l \in \Z$. This means\endnote{\cite{YYY2}, Prop.~2.3, p.~17.} that
$I_V((m,l),\Lmd) = \chi_V(\Lmh^{[m,l]})$ if~$l$ is relatively prime to~$m$.

\subsection[Indicators and duality]{} \label{IndicatorDual}
It is possible to express the equivariant Frobenius-Schur indicators in terms of the pairing between induced modules that we introduced in Paragraph~\ref{IndDual}. So, let $V$ be a finite-dimensional $H$-module with dual~$V^*$. Applying repeatedly the construction described in Paragraph~\ref{IndDual}, we get from the natural pairing 
$\langle \cdot, \cdot \rangle: V \times V^* \rightarrow K$
a pairing between~$V^{\o m}$ and~$V^{* \o m}$ that is given by
$$\langle v_1 \o v_2 \o \ldots \o v_m, 
\psi_m \o \ldots \o \psi_2 \o \psi_1 \rangle
= \langle v_1, \psi_1 \rangle
\langle v_2, \psi_2 \rangle \cdots \langle v_m, \psi_m  \rangle$$
and this pairing leads, after we choose a nonzero integral~$\Lmh$, to
a pairing $\langle \cdot, \cdot \rangle_\Lmh$ between~$\Ind(V^{\o m})$ and~$\Ind(V^{* \o m})$. We choose an integral satisfying $\eh(\Lmh)=1$. Then, if $v_1,\ldots,v_d \in V$ is a basis of~$V$ with dual basis $v^*_1,\ldots,v^*_d \in V^*$, we get the following formula for the equivariant Frobenius-Schur indicators:
\begin{prop}
For all integers $m,l \in \Z$ with $m>1$ and all central elements 
$z \in Z(D)$, we have
\begin{align*}
&I_V((m,l),z) = 
\dim(H)
\sum_{i_1,\ldots,i_m=1}^d 
\langle (\beta^l \circ \rho_m(z))(\eh \o v_{i_1} \o v_{i_2} \o \ldots \o v_{i_m}), \\
&\mspace{410mu}
\eh \o v^*_{i_m} \o \ldots \o v^*_{i_2} \o v^*_{i_1} \rangle_\Lmh
\end{align*}
\end{prop}
\begin{pf}
If $z = \sum_j \varphi_j \o h_j$, we have because $z$ is central that
$$\rho_m(z)(\varphi \o v_1 \o v_2 \o \ldots \o v_m) = 
\sum_{j} \varphi \varphi_j \o h_{j\1}.v_1 \o \ldots \o h_{j(m)}.v_m $$
and
$\beta(\varphi \o v_1 \o v_2 \o \ldots \o v_m) = 
\sum_{i=1}^n \varphi b^*_i \o v_2 \o \ldots \o v_m \o b_i.v_1$.
This implies that, under the isomorphism 
$\End(H^* \o V^{\o m}) \cong \End(H^*) \o \End(V^{\o m})$, 
\mbox{$\beta^l \circ \rho_m(z)$} decomposes into a sum of tensor products of right multiplications and endomorphisms of~$V^{\o m}$.
But the trace of the right multiplication by~$\varphi$ on~$H^*$ is given by
$\dim(H) \varphi(\Lmh)$, and the trace of an endomorphism of~$V^{\o m}$
can be found by dual bases. Since~$\sh(\Lmh)=\Lmh$, this gives the assertion.
\qed
\end{pf}

As a consequence, we can give a formula for the equivariant indicators of the dual module:
\begin{corollary}
For all $m,l \in \Z$ and all $z \in Z(D)$, we have
$$I_{V^*}((m,l),z) = I_V((m,l),\sd(z))$$
\end{corollary}
\begin{pf}
In the case where $m>1$, we can argue as in the proof of the preceding proposition to obtain the formula
\begin{align*}
&I_{V^*}((m,l),\sd(z)) = 
n \mspace{-7mu} \sum_{i_1,\ldots,i_m} \mspace{-10mu}
\langle \eh \o v_{i_1} \o \! \ldots \!\o v_{i_m}, 
\beta^l(\sd(z).(\eh \o v^*_{i_m} \o \! \ldots \! \o v^*_{i_1})) \rangle_\Lmh
\end{align*}
But by Proposition~\ref{IndDual} and the discussion after Lemma~\ref{IndDual}, the right-hand side of this formula is equal to the right-hand side of the formula in the preceding proposition, establishing the case where $m>1$. In the case~$m=1$ the assertion follows directly from the definition, since~$\ud$ is invariant under the antipode, and the case $m<0$ reduces to the cases already treated.

Now suppose that $m=0$. For $l>0$, we have
\begin{align*}
I_{V^*}((0,l),z) 
&= \dim(H) \sum_j \eh(h_j) \varphi_{j}(\Lmh_\1) 
\chi_{V^{* \o l}}(\Lmh_\2)\\
&= \dim(H) \sum_j \eh(h_j) \varphi_{j}(\Lmh_\1)
\chi_{V^{\o l}}(\sh(\Lmh_\2))
\end{align*}
Because $\Lmh$ is cocommutative and invariant under the antipode,\endnote{\cite{R4}, Prop.~3, p.~590; Thm.~3, p.~594.}
we can rewrite this as
\begin{align*}
I_{V^*}((0,l),z) 
&= \dim(H) \sum_j \eh(h_j) \varphi_{j}(\sh(\Lmh_\1))
\chi_{V^{\o l}}(\Lmh_\2)
\end{align*}
But it follows from Lemma~\ref{RoleEval} that
$\sd(z) = \sum_j \sh^*(\varphi_j) \o \sh(h_j)$, so that the
last expression is~$I_{V}((0,l),\sd(z))$. The case $l=0$ can be established by a very similar reasoning, and the case $l<0$ reduces as above to the cases already treated.
\qed
\end{pf}

It should be noted that, by comparison with Lemma~\ref{EquiFrobSchur} and in view of our definitions, this corollary shows that the dual space behaves with respect to the indicators like a $-1$-st tensor power of~$V$.

\subsection[The equivariance theorem]{} \label{EquivarThm}
The other part of the name of the quantities~$I_V((m,l),z)$ is explained
by the following equivariance theorem:
\begin{thm}
For all $g \in \SL(2,\Z)$, we have
$I_V((m,l)g,z) = I_V((m,l),g.z)$.
\end{thm}
\begin{pf}
\begin{list}{(\arabic{num})}{\usecounter{num} \leftmargin0cm \itemindent5pt}
\item
It suffices to check this on the generators given in Paragraph~\ref{GenRel}, and begin with the generator~$\gt$.
As this generator acts via~$\t$, the assertion then is that
$I_V((m,m+l),z) = I_V((m,l),\ud^{-1}z)$. For $m>1$, this follows from the fact that $\beta^m = \rho_m(\ud^{-1})$ as endomorphisms of 
$H^* \o V^{\o m}$, a fact that was already discussed after Definition~\ref{EquiFrobSchur}. The case $m=1$ follows directly from the definition, and for $m<0$ we have
\begin{align*}
I_V((m,m+l),z) &= I_V((-m,-m-l),\sd(z)) = I_V((-m,-l),\ud^{-1} \sd(z)) \\
&= I_V((-m,-l), \sd(\ud^{-1} z)) = I_V((m,l), \ud^{-1} z)
\end{align*}
In the case $m=0$ and $l>0$, write $z = \sum_j \varphi_j \o h_j$.
Because~$\ud$ is central, we then have 
$\ud^{-1}z = \sum_{i=1}^n \sum_j \varphi_j b^*_i \o b_i h_j$
and therefore
\begin{align*}
I_V((0,l),\ud^{-1}z) &= 
\dim(H) \sum_{i=1}^n \sum_j 
\eh(b_i h_j) (\varphi_{j} b^*_i)(\Lmh_\1) \chi_{V^{\o l}}(\Lmh_\2)\\
&=\dim(H) \sum_j \eh(h_j) \varphi_{j}(\Lmh_\1) \chi_{V^{\o l}}(\Lmh_\2) 
= I_V((0,l),z)
\end{align*}
The case $m=0$ and $l=0$ can be established by a very similar reasoning, and the case $m=0$ and $l<0$ reduces as above to the cases already treated.

\item
As the generator~$\gr^{-1}$ acts via~$\r^{-1}$, the assertion in this case says that $I_V((m+l,l),z) = I_V((m,l),\r^{-1}(z))$. It follows
from Proposition~\ref{NewMaps} and the fact, explained in Paragraph~\ref{CharRing}, that~$\Phi$ and~$\bar{\Phi}$
agree on the character ring, that we have
$\r^{-1}(z) = \bar{\r}^{-1}(z) = e^{-1}(z_\2) z_\1$ for every central element~$z$. The assertion therefore can also be written in the form $I_V((m+l,l),\bar{\r}(z)) = I_V((m,l),z)$, which we now establish in the case~$m>0$ and~$l>0$. For this, we write~$l=pm+q$,
where $0 \le q < m$. If $z=\sum_j \varphi_j \o h_j$, we have
\begin{align*}
&(\rho_{l+m}(\bar{\r}(z)) \circ \beta^l)
(\varphi \o v_1 \o \ldots \o v_{l+m}) = \\
&\sum_{i_1,\ldots,i_l=1}^n \sum_j \varphi_{j\1}(h_{j\2}) \;
\varphi b^*_{i_1} \cdots b^*_{i_l} \varphi_{j\2} \o \\
&\mspace{200mu}
h_{j\1}.(v_{l+1} \o \ldots \o v_{l+m} \o b_{i_1}.v_1 \o \ldots \o b_{i_l}.v_{l})
\end{align*}
By the dual basis formulas stated in the introduction, this means that we have
\begin{align*}
(\rho_{l+m}(\bar{\r}(z)) \circ \beta^l)(\varphi \o x \o w) = 
\sum_{i=1}^n \sum_j \varphi_{j\1}(h_{j\2}) \;
\varphi b^*_{i} \varphi_{j\2} \o h_{j\1}.(w \o b_{i}.x)
\end{align*}
for all $x \in V^{\o l}$ and $w \in V^{\o m}$.
This is a sum of tensor products of right multiplications on~$H^*$ and endomorphisms of~$V^{\o (l+m)}$. If $\Lmh \in H$ is an integral satisfying
$\eh(\Lmh)=1$, the right multiplication by~$\varphi \in H^*$ has the trace $n\varphi(\Lmh)$, so that the equivariant Frobenius-Schur indicator~$I_V((m+l,l),\bar{\r}(z))$, which is the trace of this map, is $n$-times the trace of the endomorphism~$f$ of~$V^{\o (l+m)}$ given by
\begin{align*}
&f(x \o w) = \sum_{i=1}^n \sum_j \varphi_{j\1}(h_{j\2}) 
(b^*_{i} \varphi_{j\2})(\Lmh) \; h_{j\1}.(w \o b_{i}.x)
\end{align*}
for $x \in V^{\o l}$ and $w \in V^{\o m}$, 
which can be rewritten in the form
\begin{align*}
f(x \o w) &=  \sum_j \varphi_{j\1}(h_{j\3}) 
\varphi_{j\2}(\Lmh_\2) \; (h_{j\1}.w \o h_{j\2}\Lmh_\1.x) \\
&= \sum_j \varphi_{j\1}(h_{j\3}) 
\varphi_{j\2}(\sh(h_{j\2})\Lmh_\2) \; (h_{j\1}.w \o \Lmh_\1.x) \\
&= \sum_j \varphi_{j}(\Lmh_\2) \; (h_{j}.w \o \Lmh_\1.x)
\end{align*}

\item
As we have discussed in Paragraph~\ref{EquiFrobSchur}, we have on the right-hand side of the assertion that
$$\rho_{m}(z) \circ \beta^l = 
\rho_{m}(z) \circ \beta^{pm} \circ \beta^{q} = 
\rho_{m}(z \ud^{-p}) \circ \beta^{q} $$
so that
\begin{align*}
&(\rho_{m}(z) \circ \beta^l)
(\varphi \o v_1 \o \ldots \o v_{m}) = \\
&\sum_{i_1,\ldots,i_q=1}^n \sum_{j_1,\ldots,j_p=1}^n \sum_j 
\varphi b^*_{i_1} \cdots b^*_{i_q} 
b^*_{j_1} \cdots b^*_{j_p} \varphi_{j} \o \\
&\mspace{200mu}
h_{j}b_{j_p} \cdots b_{j_1}.(v_{q+1} \o \ldots \o v_{m} \o b_{i_1}.v_1 \o \ldots \o b_{i_q}.v_{q})
\end{align*}
Using the dual basis formulas as before, we can write this as 
\begin{align*}
&(\rho_{m}(z) \circ \beta^l)(\varphi \o y \o t) = \\
&\sum_{i=1}^n \sum_{j_1,\ldots,j_p=1}^n \sum_j 
\varphi b^*_{i} b^*_{j_1} \cdots b^*_{j_p} \varphi_{j} \o 
h_{j}b_{j_p} \cdots b_{j_1}.(t \o b_{i}.y)
\end{align*}
for $y \in V^{\o q}$ and $t \in V^{\o (m-q)}$.
This is again a sum of tensor products of right multiplications on~$H^*$ and endomorphisms of~$V^{\o m}$, so that we see as before that the equivariant Frobenius-Schur indicator~$I_V((m,l),z)$, which is the trace of this map, is $n$-times the trace of the endomorphism~$g$ of~$V^{\o m}$
given by
\begin{align*}
&g(y \o t) = 
\sum_{i=1}^n \sum_{j_1,\ldots,j_p=1}^n \sum_j 
(b^*_{i} b^*_{j_1} \cdots b^*_{j_p} \varphi_{j})(\Lmh)
h_{j}b_{j_p} \cdots b_{j_1}.(t \o b_{i}.y)
\end{align*}
for $y \in V^{\o q}$ and $t \in V^{\o (m-q)}$, which can be rewritten  in the form
\begin{align*}
&g(y \o t) = 
\sum_j \varphi_{j}(\Lmh_{(p+2)}) \;
h_{j} \Lmh_{(p+1)} \cdots \Lmh_{(2)}.(t \o \Lmh_{\1}.y) 
\end{align*}

\item
The assertion therefore now is that 
$\Tr_{V^{\o (m+l)}}(f) = \Tr_{V^{\o m}}(g)$.
This will hold if we can show that~$g$ is the partial trace of~$f$
over the last $l$ tensor factors. Let us explain in greater detail what this means. Choose a basis $v_1,\ldots,v_d$ of~$V$ with dual basis
$v^*_1,\ldots,v^*_d$ of~$V^*$. The assertion then is that
$$g(w) = \sum_{i_1,\ldots,i_l=1}^d 
(\id_{V^{\o m}} \o v^*_{i_1} \o \ldots \o v^*_{i_l})
f(w \o v_{i_1} \o \ldots \o v_{i_l})$$
for all $w \in V^{\o m}$. To establish this, it is better to use a 
basis $w_1,\ldots,w_{dm}$ of~$W:=V^{\o m}$ with dual basis
$w^*_1,\ldots,w^*_{dm}$ of~$W^*$ as well as a basis $y_1,\ldots,y_{dq}$ of~$Y:=V^{\o q}$ with dual basis $y^*_1,\ldots,y^*_{dq}$ of~$Y^*$, and also to decompose~$w$ in the form $w = y \o t$ for $y \in V^{\o q}$
and $t \in V^{\o (m-q)}$. The assertion then becomes
$$g(y \o t) = \sum_{j_1,\ldots,j_p=1}^{dm} \sum_{i=1}^{dq} 
(\id_{V^{\o m}} \o y^*_i \o w^*_{j_1} \o \ldots \o w^*_{j_p})
f(y \o t \o y_i \o w_{j_1} \o \ldots \o w_{j_p})$$
To see this, we start at the right-hand side:
\begin{align*}
&\sum_{j_1,\ldots,j_p=1}^{dm} \sum_{i=1}^{dq} 
(\id_{V^{\o m}} \o y^*_i \o w^*_{j_1} \o \ldots \o w^*_{j_p})
f(y \o t \o y_i \o w_{j_1} \o \ldots \o w_{j_p}) = \\
&\sum_{j_1,\ldots,j_p=1}^{dm} \sum_{i=1}^{dq} \sum_{j} \varphi_j(\Lmh_{(p+3)}) (\id_{V^{\o m}} \o y^*_i \o w^*_{j_1} \o \ldots \o w^*_{j_p}) \\
&\mspace{50mu}
(h_j.w_{j_p} \o \Lmh_{(1)}.y \o \Lmh_{(2)}.t \o \Lmh_{(3)}.y_i \o \Lmh_{(4)}.w_{j_1} \o \ldots \o \Lmh_{(p+2)}.w_{j_{p-1}}) 
\end{align*}
In this expression, we can carry out the summation over~$i$,
in which case it becomes
\begin{align*}
&\sum_{j_1,\ldots,j_p=1}^{dm} \sum_{j} \varphi_j(\Lmh_{(p+3)}) 
(\id_{V^{\o m}}  \o w^*_{j_1} \o \ldots \o w^*_{j_p}) \\
&\mspace{50mu}
(h_j.w_{j_p} \o \Lmh_{(2)}.t \o \Lmh_{(3)}\Lmh_{(1)}.y \o \Lmh_{(4)}.w_{j_1} \o \ldots \o \Lmh_{(p+2)}.w_{j_{p-1}}) 
\end{align*}
Next, we carry out the summation over~$j_p$ to get
\begin{align*}
&\sum_{j_1,\ldots,j_{p-1}=1}^{dm} \sum_{j}\varphi_j(\Lmh_{(p+3)}) (\id_{V^{\o m}} \o w^*_{j_1} \o \ldots \o w^*_{j_{p-1}}) \\
&\mspace{50mu}
(h_j\Lmh_{(p+2)}.w_{j_{p-1}} \o \Lmh_{(2)}.t \o \Lmh_{(3)}\Lmh_{(1)}.y \o \Lmh_{(4)}.w_{j_1} \o \ldots \o \Lmh_{(p+1)}.w_{j_{p-2}}) 
\end{align*}
Continuing to carry out the summations up to~$j_2$, this becomes
\begin{align*}
&\sum_{j_1=1}^{dm} \sum_{j}\varphi_j(\Lmh_{(p+3)})
(\id_{V^{\o m}} \o w^*_{j_1})
(h_j (\Lmh_{(p+2)} \cdots \Lmh_{\4}).w_{j_{1}} \o \Lmh_{(2)}.t \o \Lmh_{(3)}\Lmh_{(1)}.y) 
\end{align*}
We can even carry out the summation over~$j_1$ to get
\begin{align*}
\sum_{j}\varphi_j(&\Lmh_{(p+3)})
h_j (\Lmh_{(p+2)} \cdots \Lmh_{\4}).(\Lmh_{(2)}.t \o \Lmh_{(3)}\Lmh_{(1)}.y) \\
&= \sum_{j}\varphi_j(\Lmh_{(p+2)})  
h_j (\Lmh_{(p+1)} \cdots \Lmh_{\2}).(t \o \Lmh_{(1)}.y) 
= g(y \o t)
\end{align*}
It should be pointed out that this argument needs to be slightly modified in the case $p=0$, where the summation over~$j_1,\ldots,j_p$ is empty. In fact, this implies that the computation simplifies substantially in this case. Furthermore, the reader is urged to check that this argument also covers the case~$m=1$.

\item
Next, we establish the formula $I_V((m+l,l),z) = I_V((m,l),\r^{-1}(z))$ in the case~$m>0$ and~$l=0$, in which it asserts that 
$\Tr(\rho_m(z)) = \Tr(\rho_m(\r^{-1}(z)))$.
For this, we show that~$z$ and~$\r^{-1}(z)$ have the same trace 
on every induced module, which corresponds in fact to the case~$m=1$. Now suppose that $z=\sum_j \varphi_j \o h_j \in Z(D)$,
and that~$\chi_V$ is the character of~$V$. We can write~$\rho_1(z)$ as before as a tensor product of right multiplications on~$H^*$ and an endomorphism of~$V$ and get that
$\Tr(\rho_1(z)) = n \sum_j \varphi_j(\Lmh) \chi_V(h_j)$. Since
$$\r^{-1}(z) = e^{-1}(z_\1) z_\2 = 
\sum_j \varphi_{j\2}(\sh(h_{j\1})) \; \varphi_{j\1} \o h_{j\2}$$
this implies also that
\begin{align*}
\Tr(\rho_1(&\r^{-1}(z))) = 
n \sum_j \varphi_{j\2}(\sh(h_{j\1})) \varphi_{j\1}(\Lmh) \chi_V(h_{j\2}) \\
&= n \sum_j \varphi_{j}(\Lmh \sh(h_{j\1})) \chi_V(h_{j\2}) 
= n \sum_j \varphi_{j}(\Lmh) \chi_V(h_{j})
= \Tr(\rho_1(z))
\end{align*}
as asserted.

\item
To establish the assertion $I_V((m+l,l),z) = I_V((m,l),\r^{-1}(z))$
in the case $m=0$ and $l>0$, we argue similarly: We have
\begin{align*}
I_V((0,l),\r^{-1}(z)) &= 
n \sum_j \varphi_{j\2}(\sh(h_{j\1})) \eh(h_{j\2}) 
\varphi_{j\1}(\Lmh_\1) \chi_{V^{\o l}}(\Lmh_\2) \\
&= n \sum_j \varphi_{j}(\Lmh_\1 \sh(h_{j})) \chi_{V^{\o l}}(\Lmh_\2)
\end{align*}
on the one hand and 
\begin{align*}
I_V((l,l),z) &= \Tr(\rho_l(\ud^{-1} z))  
= n \sum_{i=1}^n \sum_j (\varphi_j b_i^*)(\Lmh) \chi_{V^{\o l}}(b_i h_j)\\
&= n \sum_j \varphi_j(\Lmh_\1) \chi_{V^{\o l}}(\Lmh_\2 h_j)
\end{align*}
on the other hand. Both expressions are equal by the basic Casimir properties of the integral.\endnote{\cite{LR2}, Lem.~1.2, p.~270.}
A very similar reasoning establishes the formula in the case~$m=0$ and~$l=0$, and for $m=0$ and~$l<0$ we have
\begin{align*}
I_V((l,l),z) &= I_V((-l,-l),\sd(z)) \\
&= I_V((0,-l),\r^{-1}(\sd(z))) = I_V((0,l),\r^{-1}(z))
\end{align*}
because~$\r$ and~$\sd$ commute. 

\item
We have now established that
$I_V((m+l,l),z) = I_V((m,l),\r^{-1}(z))$
whenever $m \ge 0$ and $l \ge 0$. Instead of establishing the 
remaining cases, we use this fact to prove that 
$I_V((l,m),\v(z)) = I_V((m,-l),z)$ if~$m>0$ and~$l \ge 0$.
For this, we write $l=am+b$, where~$a \ge 0$ and~$0 \le b < m$, and argue by induction on~$a$. The induction beginning is the case~$a=0$, in which we have $l=b<m$.
We have seen in Paragraph~\ref{GenRel} that
$\gv = \gt^{-1} \gr^{-1} \gt^{-1}$, so that by the first step we get
\begin{align*}
I_V((l,m),\v(z)) &= I_V((l,m),(\t^{-1} \circ \r^{-1} \circ \t^{-1})(z)) \\
&= I_V((l,m-l),(\r^{-1} \circ \t^{-1})(z)) 
\end{align*}
Because $m-l > 0$, we can apply the identity established above to rewrite this further as
\begin{align*}
I_V((l,m),\v(z)) &= I_V((m,m-l),\t^{-1}(z)) = I_V((m,-l),z)
\end{align*}
where we have applied the first step again.

For the induction step, note that it follows from the discussion in Paragraph~\ref{GenRel} that~$\gr \gv =\gv \gt$. By the induction assumption, we have
$$I_V(((a-1)m+b,m),\v(z)) = I_V((m,-(a-1)m-b),z)$$
which means that
$I_V((am+b,m),\r(\v(z))) = I_V((m,-am-b),\t(z))$.
By the preceding commutation relation, this asserts that
$$I_V((l,m),\v(\t(z))) = I_V((m,-l),\t(z))$$
so that the assertion now follows by substituting~$\t^{-1}(z)$ for~$z$.

\item
Inspection of the preceding argument shows that it also proves the
formula $I_V((l,m),\v(z)) = I_V((m,-l),z)$ in the case~$m=l=0$. To establish it if~$m \le 0$ and~$l>0$, note that it asserts in this case that
$$I_V((l,m),\v(z)) = I_V((-m,l),\sd(z))$$
Since $\sd(z)=\v^2(z)$, this is equivalent to
$I_V((l,m),z) = I_V((-m,l),\v(z))$, a fact that we have just 
established.

The proof that $I_V((l,m),\v(z)) = I_V((m,-l),z)$ if~$m \ge 0$ and~$l<0$
is similar: The assertion then is that
$$I_V((-l,-m),\sd(\v(z))) = I_V((m,-l),z)$$
If we substitute~$\v(z)$ for~$z$, this becomes
$I_V((-l,-m),z) = I_V((m,-l),\v(z))$, which we have obtained 
already.

Finally, if~$m<0$ and~$l \le 0$, the assertion is that
$$I_V((-l,-m),\sd(\v(z))) = I_V((-m,l),\sd(z))$$
which upon substituting~$\sd(z)$ for~$z$ also reduces to the established case.
\qed
\end{list}
\end{pf}

\subsection[The orbit theorem]{} \label{OrbitThm}
In our situation, the Drinfel'd element~$\ud$ has finite order.\endnote{\cite{EG4}, Thm.~2.5, p.~133.} This order is called the exponent of~$H$; we denote it by~$N$. It is known that~$N$ divides~$\dim(H)^3$; however, the original conjecture of Y.~Kashina, namely that~$N$ divides~$\dim(H)$, is still open.\endnote{\cite{EG4}, Thm.~4.3, p.~136; \cite{KashAnti}, p.~1261; \cite{KashPow}, p.~159.}

We now consider the cyclotomic field~$\Q_N \subset K$ that arises by adjoining to the prime field~$\Q \subset K$ all $N$-th roots of unity that are contained in~$K$. We denote by~$Z_{\Q_N}(D)$ the span of the centrally primitive idempotents $e_1,\ldots,e_k$ introduced in Paragraph~\ref{CharRing} over the subfield~$\Q_N$ of~$K$. This space has the following property:
\begin{lemma}
$Z_{\Q_N}(D)$ is invariant under the action of the modular group.
\end{lemma}
\begin{pf}
It suffices to show that it is invariant under~$\t$ and~$\v$. The fact that~$\ud$ has order~$N$ means for the expansion~$\ud = \sum_{i=1}^k u_i e_i$
considered in Paragraph~\ref{VerlMat} that the coefficients~$u_i$ are $N$-th roots of unity. Since~$\t$ is the multiplication by~$\ud^{-1}$, we see that $Z_{\Q_N}(D)$ is invariant under~$\t$.

To see that $Z_{\Q_N}(D)$ is invariant under~$\v$, recall\endnote{\cite{YYY2}, Prop.~6.2, p.~44.}
that the entries~$s_{ij}$ of the Verlinde matrix are contained in~$\Q_N$. 
Therefore, the assertion follows from the formula
$\v(e_j) = 
\frac{1}{\dim(H)} \sum_{i=1}^k \frac{n_j}{n_i} s_{j^*i} e_i$
established, taking Paragraph~\ref{RoleEval} into account, in Corollary~\ref{VerlMat}.
\qed
\end{pf}

This lemma has the following consequence for the indicators:
\begin{prop}
For $z \in Z_{\Q_N}(D)$, we have $I_V((m,l),z) \in \Q_N$.
\end{prop}
\begin{pf}
In the case~$(m,l)=(0,0)$, it follows easily from the definition that
$I_V((m,l),z)$ is the trace of the action of~$z$ on the induced module~$\Ind(K)$ of the trivial module, which is in~$\Q_N$ if $z \in Z_{\Q_N}(D)$. If $(m,l) \neq (0,0)$, we set $t := \gcd(m,l)>0$. By Corollary~\ref{CongSubgr}, we can find~$g \in \Gamma$ such that
$(m,l)=(t,0)g$. By Theorem~\ref{EquivarThm}, we then have
$$I_V((m,l),z) = I_V((t,0),g.z) = \Tr(\rho_t(g.z))$$
which is in~$\Q_N$ since $g.z \in Z_{\Q_N}(D)$ by the preceding lemma.
\qed
\end{pf}

We now consider the principal congruence subgroup~$\Gamma(N)$ corresponding to the exponent~$N$. The following orbit theorem asserts that the indicator depends only on the $\Gamma(N)$-orbit of the lattice point:
\begin{thm}
Suppose that two lattice points $(m,l)$ and $(m',l')$ are in the same 
$\Gamma(N)$-orbit. Then we have $I_V((m,l),z) = I_V((m',l'),z)$
for every $H$-module~$V$ and every $z \in Z(D)$.
\end{thm}
\begin{pf}
\begin{list}{(\arabic{num})}{\usecounter{num} \leftmargin0cm \itemindent5pt}
\item
We fix an $H$-module~$V$, and introduce an equivalence relation~$\approx$ on the lattice~$\Z^2$ by defining $(m,l) \approx (m',l')$ if and only if
$$I_{V}((m,l),z) = I_{V}((m',l'),z)$$
for all $z \in Z(D)$. Then~$\approx$ is a congruence relation, since, for $g \in \Gamma$, $(m,l) \approx (m',l')$ implies in particular that 
$I_{V}((m,l),g.z) = I_{V}((m',l'),g.z)$, which yields
$I_{V}((m,l)g,z) = I_{V}((m',l')g,z)$
by Theorem~\ref{EquivarThm}, so that~$(m,l)g \approx (m',l')g$.
Note that there is a slight adaption necessary: In Section~\ref{Sec:ModGroup}, we have considered the left action of the modular group on the lattice points, considered as columns, whereas we consider here the transposed right action, where the lattice points are considered as rows.

\item
We now want to check that~$\approx$ satisfies the two defining properties of the congruence relation~$\sim$ listed in Paragraph~\ref{OrbCongRel}. Although the transpose of~$\gt^N$ is~$\gr^{-N}$, we can also work with~$\gt^N$ in the transposed situation, since~$\gt^N$ and~$\gr^{N}$ are conjugate and the inverse sign does not matter. For the first property, we therefore have to check that $(m,l) \approx (m,l)\gt^N$. But this is immediate, since we have $\t^N(z) = \ud^{-N}z = z$ and therefore
\begin{align*}
I_{V}((m,l)\gt^N,z) = I_{V}((m,l),\t^N(z)) = I_{V}((m,l),z)
\end{align*}
by Theorem~\ref{EquivarThm}.

\item
For the second property, we are given $q \in \Z$ that satisfies
$q \equiv 1 \pmod{N}$ and~$t:=\gcd(m,l) = \gcd(m,ql)$, and have to establish that 
$(m,l) \approx (m,ql)$, in other words, that
$$I_{V}((m,l),z) = I_{V}((m,ql),z)$$
for all $z \in Z(D)$. We treat the case~$m>0$ first, where we also have~$t>0$. It is sufficient to establish this in the case where $z=e_i$ is a centrally primitive idempotent. 

If we write $m=tm'$, $l=tl'$, we have that~$q$ is relatively prime to~$m'$.
Consider the cyclotomic field~$\Q_{Nm'} \subset K$. Because~$q$ is relatively prime to~$Nm'$, there is a unique 
automorphism~$\sigma_q \in \Gal(\Q_{Nm'}/\Q)$ with the property that
$$\sigma_q(\zeta) = \zeta^{q}$$
for every $Nm'$-th root of unity~$\zeta$. Because $q \equiv 1 \pmod{N}$, we have 
$\sigma_q(\zeta) = \zeta$ if~$\zeta$ is an $N$-th root of unity, and therefore even~$\sigma_q \in \Gal(\Q_{Nm'}/\Q_N)$. By the preceding proposition, this means that
$\sigma_q(I_{V}((m,l),z)) = I_{V}((m,l),z)$. On the other hand, we have by definition that
$I_{V}((m,l),z) = \Tr(\beta^{l} \circ \rho_{m}(z))$,
where $\beta = \beta_{V,V^{\o (m-1)}}$, properly understood in the case~$m=1$. Now we have
$$(\beta^{l})^{m'} = \beta^{tl'm'} = (\beta^{m})^{l'} = \rho_m(\ud^{-1})^{l'}$$
so that $(\beta^{l})^{m'N} = \id_{\Ind(V^{\o m})}$.
Since~$z=e_i$, $\rho_{m}(z)$ is the projection to the isotypical component of~$V_i$ in~$\Ind(V^{\o m})$, so that~$\beta^l \circ \rho_{m}(z)$ coincides with~$\beta^l$ on this isotypical component and vanishes on the other isotypical components. In particular, the eigenvalues 
of~$\beta^l \circ \rho_{m}(z)$ are $Nm'$-th roots of unity, and the eigenvalues of its $q$-th power $\beta^{lq} \circ \rho_{m}(z)$ are the $q$-th powers of its eigenvalues, so that for the trace we get the formula
$$I_V((m,lq),z) = \Tr(\beta^{lq} \circ \rho_{m}(z))
= \sigma_q(\Tr(\beta^{l} \circ \rho_{m}(z)))
= \sigma_q(I_V((m,l),z))$$
Combining this with our earlier observation, this establishes the assertion in the case~$m>0$. 

\item
The case $m<0$ reduces immediately to the case just treated, since
$$I_{V}((m,l),z) = I_{V}((-m,-l),\sd(z)) 
= I_{V}((-m,-ql),\sd(z)) =I_{V}((m,ql),z)$$
Now suppose that~$m=0$. If also $l=0$, the assertion is obvious, so that we can assume that~$l \neq 0$. In this case, the conditions that
$q \equiv 1 \pmod{N}$ and~$\gcd(m,l) = \gcd(m,ql)$ imply
that~$ql= \pm l$, so that $q= \pm 1$. The case $q=1$ is obvious, so that we now assume that~$q=-1$, which can only happen if~$N=1$ or~$N=2$. A Hopf algebra of exponent~$1$ is one-dimensional, and for a Hopf algebra of exponent~$2$ we have
$h = h_\1 h_\2 \sh(h_\3) = \sh(h)$,
so that the antipode of~$H$ is the identity. Then the antipode of~$H^*$
is also the identity, so that~$H$ is commutative and cocommutative, which implies that~$D$ is commutative and cocommutative,\endnote{\cite{R2}, Prop.~6, p.~302.} so that its antipode is again the identity. 
We therefore have
$$I_{V}((0,-l),z) = I_{V}((0,l),\sd(z)) = I_{V}((0,l),z)$$
in the case $N=1$ as well as in the case~$N=2$, and the second defining property is completely established.

\item
In Paragraph~\ref{OrbCongRel}, we have defined the relation~$\sim$ as
the intersection of all congruence relations that satisfy the two defining properties just verified. Therefore $(m,l) \sim (m',l')$ implies that 
$(m,l) \approx (m',l')$. But by Theorem~\ref{OrbCongRel},  $(m,l) \sim (m',l')$ means that~$(m,l)$ and~$(m',l')$ are in the same $\Gamma(N)$-orbit, so that this is exactly the assertion.
\qed
\end{list}
\end{pf}

We put down one easy special case of this theorem that will be needed later:
\begin{corollary}
For an $H$-module~$V$, let $\chi$ be the character of~$\Ind(V)$. Then we have $\chi(g.z) = \chi(z)$ for all $g \in \Gamma(N)$ and all~$z \in Z(D)$.
\end{corollary}
\begin{pf}
We have $\chi(z) = \Tr(\rho_1(z)) = I_V((1,0),z)$, so that by Theorem~\ref{EquivarThm} the assertion is
equivalent to $I_V((1,0)g,z) = I_V((1,0),z)$. But since $(1,0)$ and
$(1,0)g$ are obviously in the same $\Gamma(N)$-orbit, this follows directly from the preceding theorem.
\qed
\end{pf}

\newpage
\section{The congruence subgroup theorem} \label{Sec:CongSubgrThm}
\subsection[The dual projective representation of the modular group]{} \label{DualRepMod}
We now will apply the results of Section~\ref{Sec:EquiFrobSchur} to prove that the kernel of the projective representation of the modular group on the center of a semisimple factorizable Hopf algebra is a congruence subgroup. Note that the kernel of a group homomorphism consists of those elements that are mapped to the unit element, which in the projective linear group consists of all nonzero scalar multiples of the identity. If therefore a projective representation is induced from an ordinary linear representation, the kernel of the projective representation is in general larger than the kernel of the linear representation.

So, let~$A$ be a semisimple factorizable Hopf algebra over our algebraically closed base field~$K$ of characteristic zero. The exponent of~$A$ will be denoted by~$N$. Otherwise, we will use the notation introduced in Section~\ref{Sec:FaktSemisim}; in particular, we will use
the inverse Drinfel'd element~$\ua^{-1}$ as our ribbon element. However, the whole discussion in Section~\ref{Sec:FaktSemisim} depended on a parameter~$\kappa$ that was introduced in Paragraph~\ref{VerlMat}; we will now dispose of this parameter in the following intricate way:
If $\chi_R(\ua)=\chi_R(\ua^{-1})$, in which case, as discussed in Paragraph~\ref{RibEl}, our representation is linear, we set 
$\kappa = \frac{1}{\chi_R(\ua)}$, so that our integral~$\rha \in A^*$ satisfies $\rha(\ua)=\rha(\ua^{-1})=1$. By Lemma~\ref{RibEl}, this integral also satisfies $(\rha \o \rha)(R'R)=1$. Note that in the case where $A=D(H)$ is the Drinfel'd double of a semisimple Hopf algebra~$H$, we therefore pick here the integral~$\rhd$ from Paragraph~\ref{RoleEval}.

If $\chi_R(\ua) \neq \chi_R(\ua^{-1})$, in which case our representation is not linear, we choose~$\kappa$ so that~$\kappa^2 = \frac{1}{\dim(A)}$.
This obviously only determines~$\kappa$ up to a sign, but in view of the formula $\chi_R(\ua) \chi_R(\ua^{-1})=\dim(A)$ observed in Paragraph~\ref{MatIdent}, implies that 
$(\rha \o \rha)(R'R) = \rha(\ua)\rha(\ua^{-1}) = 1$. In particular, the choice of~$\rha$ in both cases is compatible, the only difference is that in the first, linear case even the sign of~$\kappa$ is determined. In any case, this choice of~$\kappa$ is the one that makes~$\Psi$ equivariant, as explained in Paragraph~\ref{ModGrDoub}, as it is compatible with the condition that~$\rhd = \Psi^*(\rha \o \rha)$.

Next, we discuss the dual of our projective representation. We denote the morphism that a linear map~$f$ induces between the corresponding projective spaces by~$P(f)$. By considering the character ring as dual to the center, we can dualize the projective representation of the modular group on the center to a projective representation of the modular group on the character ring as follows: If the group element $g \in \SL(2,\Z)$ is represented by the equivalence 
class~$P(f) \in \PGL(Z(A))$, consider for $\chi \in \Ch(A)$  the character $\chi' \in \Ch(A)$ that satisfies
$$\chi'(z) = \chi(f^{-1}(z))$$
for all $z \in Z(A)$, and set $g.\bar{\chi} = \bar{\chi}'$. This does not depend on the choice of the representative~$f$ and gives a projective representation of~$\SL(2,\Z)$ on~$P(\Ch(A))$.

This construction raises the question whether the isomorphism~$P(\iota)$ that~$\iota$ induces 
between the projective spaces of the center and the character ring is equivariant with respect to the corresponding actions. This is the case only after the action
on the center is modified with the help of the automorphism introduced in Definition~\ref{GenRel}:
\begin{prop}
For all $g \in \Gamma$, the diagram
$$\square<1`-1`-1`1;1400`700>[P(\Ch(A))`P(\Ch(A))`P(Z(A))`P(Z(A));
\bar{\chi} \mapsto g.\bar{\chi}`P(\iota)`P(\iota)`\bar{z} \mapsto \tilde{g}.\bar{z}]$$
commutes.
\end{prop}
\begin{pf}
It suffices to check this for the generators~$\gv$ and~$\gt$, for which we have seen in Paragraph~\ref{GenRel} that $\tilde{\gv}= \gv^{-1}$ and $\tilde{\gt}= \gt^{-1}$. In the case of~$\gv$, the assertion therefore follows from Proposition~\ref{RoleInt}, and in the case of~$\gt$ it follows from the corresponding diagram given in Paragraph~\ref{RibEl}.
\qed
\end{pf}

The analogous question for~$P(\Phi)$ can be deduced from this proposition:
\begin{corollary}
For all $g \in \Gamma$, the diagram
$$\square<1`1`1`1;1400`700>[P(\Ch(A))`P(\Ch(A))`P(Z(A))`P(Z(A));
\bar{\chi} \mapsto g.\bar{\chi}`P(\Phi)`P(\Phi)`\bar{z} \mapsto (\gv\tilde{g}\gv^{-1}).\bar{z}]$$
commutes.
\end{corollary}
\begin{pf}
This will follow from the preceding proposition by reversing the
vertical arrows if we can verify that $\v \circ \Phi = \iota^{-1}$
on the character ring, or equivalently that $\Phi \circ \iota = \v^{-1}$.
But as~$\Phi$ and~$\bar{\Phi}$ agree on the character ring, we have
from the definition of~$\Phi$ that $\Phi \circ \iota = \sa \circ \v$
on the center, which in view of Corollary~\ref{InvSigma} implies the assertion.
\qed
\end{pf}

Note that in the case where we have a linear representation of~$\SL(2,\Z)$
on the center~$Z(A)$, we will also get in this way a linear representation on~$\Ch(A)$, and~$\iota$ and~$\Phi$ will then be equivariant with respect to the linear representations.

\subsection[Induction and multiplicities]{} \label{IndMult}
If $\Lma \in A$ is an integral which satisfies $\ea(\Lma)=1$, we define the bilinear form
$\langle \cdot, \cdot \rangle_*$ on $\Ch(A)$ by the equation
$$\langle \chi, \chi' \rangle_* := (\chi \sa^*(\chi'))(\Lma)$$
Then we have 
$\langle \chi_i, \chi_j \rangle_* = \delta_{ij} = \dim \Hom_A(V_i,V_j)$
for the irreducible characters,\endnote{\cite{So4}, Par.~3.5, p.~211.}
which shows that this bilinear form is nondegenerate and symmetric.

For $l=1,\ldots,k$, we denote the character of the induced
$D(A)$-module~$\Ind(V_l)$ by~$\eta_l$. Using this, we can express this bilinear form as follows:
\begin{prop}
For all $\chi, \chi' \in \Ch(A)$, we have
$$\langle \chi \chi', \chi_l \rangle_*
= \frac{1}{\dim(A)} \;
\eta_l(\Psi^{-1}(\iota^{-1}(\chi) \o \iota^{-1}(\chi')))$$
\end{prop}
\begin{pf}
Because both sides of the equation are bilinear in~$\chi$ and~$\chi'$, we can assume that both characters are irreducible, so that
$\chi = \chi_i$ and~$\chi' = \chi_j$ for some $i,j \le k$.
The vector space $V_i \otimes V_j$ can be considered as an 
$A \o A$-module by the componentwise action; it is then simple. We can turn it into a $D(A)$-module by pullback along~$\Psi$. We denote this $D(A)$-module by~$U$; since~$\Psi$ is an isomorphism, this module is also simple, and the centrally primitive idempotent in $Z(D(A))$ corresponding to~$U$ is~$\Psi^{-1}(e_i \o e_j)$. Because $\Psi(\ea \o a) = \da(a)$, the restriction of~$U$ to~$A \subset D(A)$ is the same as the pullback of the $A \o A$-module structure along~$\da$, which is exactly how the tensor product~$V_i \o V_j$ of $A$-modules is formed.

We now have by the Frobenius reciprocity theorem\endnote{\cite{Lang}, Chap.~XVIII, \S~7, p.~689.} that
\begin{align*}
\langle \chi_i \chi_j, \chi_l \rangle_* &= 
\dim \Hom_{A}(V_l,V_i \o V_j)
= \dim \Hom_{D(A)}(\Ind(V_l),U) \\
&= \dim \Hom_{D(A)}(U,\Ind(V_l))
= \frac{1}{n_i n_j}\eta_l(\Psi^{-1}(e_i \o e_j))
\end{align*}
By Proposition~\ref{VerlMat}, we have 
$\iota(e_i) = \kappa n_i \chi_i$,
so that 
$$\frac{e_i}{n_i} = \kappa \mspace{2mu} \iota^{-1}(\chi_i)$$
Inserting this into the preceding formula and using 
$\kappa^2 = \frac{1}{\dim(A)}$, we get the assertion.
\qed
\end{pf}

This proposition has the following consequence:
\begin{corollary}
Suppose that $\chi_R(u) = \chi_R(u^{-1})$. Then we have
$$(\tilde{g}.\chi)(g.\chi') = \chi \chi'$$
for all $g \in \Gamma(N)$ and all $\chi,\chi' \in \Ch(A)$.
\end{corollary}
\begin{pf}
Recall that the assumption implies that the representation of the modular group is linear. By the nondegeneracy of the bilinear form above, it suffices to show that
$\langle (\tilde{g}.\chi)(g.\chi'),\chi_l \rangle_* = 
\langle \chi \chi', \chi_l \rangle_*$
for all $l=1,\ldots,k$. Now we get from the preceding proposition for $g \in \Gamma(N)$ that
\begin{align*}
\langle (\tilde{g}.\chi)(g.\chi'),\chi_l \rangle_* 
&= \frac{1}{\dim(A)} \; \eta_l(\Psi^{-1}(\iota^{-1}(\tilde{g}.\chi) \o 
\iota^{-1}(g.\chi'))) \\
&= \frac{1}{\dim(A)} \; 
\eta_l(\Psi^{-1}(g.\iota^{-1}(\chi) \o \tilde{g}.\iota^{-1}(\chi')))\\
&= \frac{1}{\dim(A)} \; 
\eta_l(g.\Psi^{-1}(\iota^{-1}(\chi) \o \iota^{-1}(\chi'))) \\
&= \frac{1}{\dim(A)} \; 
\eta_l(\Psi^{-1}(\iota^{-1}(\chi) \o \iota^{-1}(\chi')))  
= \langle \chi \chi', \chi_l \rangle_*
\end{align*}
where the second equality follows from Proposition~\ref{DualRepMod}, the third from Proposition~\ref{ModGrDoub}, and the fourth from Corollary~\ref{OrbitThm}.
\qed
\end{pf}

\subsection[The congruence subgroup theorem]{} \label{CongDrinfDoubl}
We now turn for a moment to the special case where $A=D(H)$, the Drinfel'd double of a semisimple Hopf algebra~$H$, which is denoted by~$D$ to distinguish it from the general case. Note that~$H$ and~$D$ have the same exponent~$N$.\endnote{\cite{KashPow}, Sec.~3, Thm.~3.4, p.~170; \cite{EG4}, Cor.~3.4, p.~135.} 
In this case, we now prove the following congruence subgroup theorem:
\begin{thm}
The kernel of the representation of the modular group on the center of~$D$ is a congruence subgroup of level~$N$.
\end{thm}
\begin{pf}
It follows from Corollary~\ref{OrbitThm} that the character of every $D$-module that is induced from an $H$-module is invariant under~$\Gamma(N)$. The regular representation of~$D$ is induced from the regular representation of~$H$, and therefore its character is invariant under~$\Gamma(N)$. This now implies that the counit is also invariant under~$\Gamma(N)$. To see this, recall that Corollary~\ref{VerlMat} gives in particular that 
$\v_*(p_1) = \frac{1}{\dim(H)} \chi_{1}$, which means that
$\v_*(\chi_R) = \dim(H) \ed$. In view of Proposition~\ref{RoleInt}, this in turn says that $\gv^{-1}.\chi_R = \dim(H) \ed$. For an element 
$g \in \Gamma(N)$, this gives
$$g.\ed = \frac{1}{\dim(H)} g \gv^{-1}.\chi_R
= \frac{1}{\dim(H)} \gv^{-1} (\gv g \gv^{-1}).\chi_R
= \frac{1}{\dim(H)}  \gv^{-1}.\chi_R = \ed$$
since $\Gamma(N)$ is a normal subgroup.

If we now substitute~$\ed$ for~$\chi'$ in Corollary~\ref{IndMult}, we get
that $\tilde{g}.\chi = \chi$ for every character~$\chi$ of~$D$ and every
$g \in \Gamma(N)$, and since conjugation by~$\ga$ restricts to an automorphism of~$\Gamma(N)$, we see that every character is invariant under~$\Gamma(N)$. But considering how we defined this action in Paragraph~\ref{DualRepMod}, this implies that every central element is invariant under~$\Gamma(N)$, since the pairing between the character ring and the center is nondegenerate. In other words, the kernel of the representation contains~$\Gamma(N)$, and therefore is a congruence subgroup.

It remains to be proved that the level of the kernel is exactly~$N$. But if there were some $N' < N$ with the property that~$\Gamma(N')$ would also be contained in the kernel, this would in particular imply that~$\gt^{-N'}$
acts trivially on the center, which would mean that~$\ud^{N'}=1$. But 
as~$N$ is by definition the order of~$\ud$, this cannot be the case.
\qed
\end{pf}

\subsection[The projective congruence subgroup theorem]{} \label{ProjCong}
Returning to the general case of an arbitrary factorizable semisimple Hopf algebra~$A$, in which the action of the modular group in general is only projective, we can still look at the kernel of the corresponding group homomorphism to~$\PGL(Z(A))$. For this kernel, the following analogue of Theorem~\ref{CongDrinfDoubl} holds:
\begin{thm}
The kernel of the projective representation of the modular group on the center of~$A$ is a congruence subgroup of level~$N$.
\end{thm}
\begin{pf}
To show that~$\Gamma(N)$ is contained in the kernel, suppose that this is not the case, and choose $g \in \Gamma(N)$ that is not mapped to the identity in~$\PGL(Z(A))$. Choose a representative~$f \in \GL(Z(A))$
for the action of~$g$, and also a representative~$\tilde{f}$
for the action of~$\tilde{g}$. Because~$f$ is not a scalar multiple of the identity, there exists an element~$z \in Z(A)$ such that~$z$ and~$f(z)$ are not proportional, and therefore linearly independent.
This implies that $z \o z$ and $f(z) \o \tilde{f}(z)$ are not proportional.
But we saw in Paragraph~\ref{ModGrDoub} that the tensor product of our projective representation with its conjugate under~$\ga$ is induced by a linear representation, and that this linear representation is via~$\Psi$ isomorphic to the representation on~$Z(D(A))$, on which~$g$ acts trivially by Theorem~\ref{CongDrinfDoubl}. But this means that it acts trivially on
$\bar{z} \o \bar{z}$, too, contradicting the fact that 
$z \o z$ and $f(z) \o \tilde{f}(z)$ are not proportional.

As in the proof of Theorem~\ref{CongDrinfDoubl}, it remains to be proved that the level of the kernel is exactly~$N$. Now if there were some 
$N' < N$ with the property that~$\Gamma(N')$ would also be contained in the kernel, this would in this case only imply that~$\ud^{N'}$ acts on the center by multiplication by a scalar. But as~$\ud$ always preserves the integral, this scalar has to be~$1$, which contradicts the definition of~$N$ as in the previous case.
\qed
\end{pf}

\newpage
\section{The action of the Galois group} \label{Sec:GaloisGroup}
\subsection[The Galois group and the character ring]{} \label{GaloisChar}
In this section, we will introduce an action of the Galois group of the cyclotomic field determined by the exponent that will turn out to be intimately connected to the action of the modular group.
As in Section~\ref{Sec:CongSubgrThm}, we consider a semisimple factorizable Hopf algebra~$A$ over our algebraically closed base field~$K$ of characteristic zero. The exponent of~$A$ will be denoted by~$N$. Otherwise, we will use the notation introduced in Section~\ref{Sec:FaktSemisim}; in particular, we will use the inverse Drinfel'd element~$\ua^{-1}$ as our ribbon element. The constant~$\kappa$, which determines the normalization of the integral~$\rha$, is chosen as in Paragraph~\ref{DualRepMod}, and~$\iota$ is defined using this integral~$\rha$.

In Paragraph~\ref{CharRing}, we have constructed the algebra homomorphisms
$\xi_1,\ldots,\xi_k$ from the character ring to the center. If we denote by~$\Ch_\Q(A)$ the span of the irreducible characters~$\chi_1,\ldots,\chi_k$ not over~$K$, but over the rational numbers~$\Q \subset K$, it can be shown\endnote{\cite{YYY2}, Prop.~6.2, p.~44.} as in the case of the Drinfel'd double that the images $\xi_i(\Ch_\Q(A))$ are contained in the cyclotomic field~$\Q_N$ determined by the exponent. The restriction of~$\xi_i$ to~$\Ch_\Q(A)$ is a $\Q$-algebra homomorphism to~$\Q_N$, and clearly~$\xi_i$ is uniquely determined by this restriction. For $\sigma \in \Gal(\Q_N/\Q)$, the map 
$\sigma \circ \xi_i$ is again a $\Q$-algebra homomorphism from~$\Ch_\Q(A)$ to~$\Q_N$, which must coincide with one of the restrictions of
$\xi_1,\ldots,\xi_k$. In this way, we get an action of~$\Gal(\Q_N/\Q)$
on the set~$\{1,\ldots,k\}$ so that
$$\sigma \circ \xi_i \mid_{\Ch_\Q(A)} = \xi_{\sigma.i} \mid_{\Ch_\Q(A)}$$
From this permutation representation, we get an action of the Galois group on the character ring~$\Ch(A)$ over the full base field~$K$ by permuting the characters accordingly, in the sense that 
$$\sigma.\chi_i := \chi_{\sigma.i}$$
and extending this action $K$-linearly.

The following lemma lists some basic properties of this action:
\begin{lemb}
\begin{enumerate}
\item $n_{\sigma.i} = n_i$

\item $\sigma(s_{ij}) = s_{\sigma.i,j} = s_{i,\sigma.j}$

\item $\sigma.\chi_1 = \chi_1$
\end{enumerate}
\end{lemb}
\begin{pf}
For the first assertion, recall the formula
$\xi_i(\chi_A) = \frac{\dim(A)}{n_i^2}$ established in Corollary~\ref{MatIdent}. Using this, we get
\begin{align*}
\frac{\dim(A)}{n_{\sigma.i}^2} = \xi_{\sigma.i}(\chi_A)
= \sigma(\xi_i(\chi_A)) = \sigma(\frac{\dim(A)}{n_i^2}) 
= \frac{\dim(A)}{n_i^2}
\end{align*}
from which the first assertion is immediate.

For the second assertion, recall the formula
$s_{ij} = n_i \xi_i(\chi_j)$
from Lemma~\ref{VerlMat}, from which we get
\begin{align*}
\sigma(s_{ij}) = n_i \sigma(\xi_i(\chi_j))
= n_{\sigma.i} \xi_{\sigma.i}(\chi_j) = s_{\sigma.i,j}
\end{align*}
Furthermore, we observed in Lemma~\ref{VerlMat} that the Verlinde matrix is symmetric, from which we see that
$\sigma(s_{ij}) = \sigma(s_{ji}) = s_{\sigma.j,i} = s_{i,\sigma.j}$.

For the third assertion, note that we have 
$$\xi_1(\chi_i) = \omega_1(\Phi(\chi_i)) = \ea(\Phi(\chi_i)) = n_i$$
so that $\sigma(\xi_1(\chi_i)) = \xi_1(\chi_i)$. This implies
$\sigma.1 = 1$.
\qed
\end{pf}

The action on the character ring can also be viewed in a different way:
$\Ch_\Q(A)$ is a commutative semisimple $\Q$-algebra, and therefore by Wedderburn's theorem\endnote{\cite{FarbDennis}, Thm.~1.11, p.~40.} isomorphic to a direct sum of fields. As in the case of the Drinfel'd double,\endnote{\cite{YYY2}, Prop.~6.2, p.~44.} it can be shown that these fields are subfields of the cyclotomic field~$\Q_N$. Since the Galois group of~$\Q_N$ is abelian, every subfield of the cyclotomic field is normal, and therefore preserved by the action of the Galois group of~$\Q_N$. We therefore get an action of this Galois group on~$\Ch_\Q(A)$ as the sum of the actions on the Wedderburn components. The action of the Galois group constructed above is exactly the $K$-linear extension of this action to 
$\Ch(A) \cong \Ch_\Q(A) \o_\Q K$. To see this, note that the restrictions of~$\xi_1,\ldots,\xi_k$ to~$\Ch_\Q(A)$ arise by projecting to some Wedderburn component and then embedding it into~$\Q_N \subset K$. Because the Galois group is abelian, it does not matter whether we first act on the Wedderburn component and then embed into~$\Q_N$, or first embed and then act. In other words, the action on the Wedderburn components
satisfies
$$\xi_i(\sigma.\chi_j) = \sigma(\xi_i(\chi_j))$$
for all $\sigma \in \Gal(\Q_N/\Q)$. But on the other hand it follows from the preceding lemma that we have 
$\xi_i(\chi_{\sigma.j}) = \sigma(\xi_i(\chi_j))$,
so that the action constructed before also satisfies this equation, which means that the two actions have to coincide.

If we consider the cyclotomic field~$\Q_N$ as a subfield of the complex numbers, complex conjugation restricts to an automorphism of the cyclotomic field, which we denote by $\gamma \in \Gal(\Q_N/\Q)$. It does not depend on the way how the cyclotomic field is embedded into the complex numbers,
as it can be characterized by the property that it maps any
$N$-th root of unity to its inverse. As proved in several places in the literature,\endnote{\cite{Beauv}, \S~6, Cor.~6.2; \cite{NR2}, Sec.~2, Rem.~11, p.~1088; \cite{With}, Sec.~3, Prop.~3.1, p.~885.} it acts on~$\Ch(A)$ via the antipode:
\begin{prop}
For all $\chi \in \Ch(A)$, we have $\gamma.\chi = \sa^*(\chi)$.
\end{prop}
Stated differently, this asserts that $\gamma.i=i^*$, which in particular implies that $\sigma.(i^*) = (\sigma.i)^*$ for all 
$\sigma \in \Gal(\Q_N/\Q)$, as the Galois group~$\Gal(\Q_N/\Q)$ is abelian. Using this, we can deduce further properties of our action:
\begin{corb}
\begin{enumerate}
\item 
$\displaystyle
\sigma.p_i = p_{\sigma^{-1}.i}$

\item
$\displaystyle
\v_*(\sigma(\chi)) = \sigma^{-1}(\v_*(\chi))$
\end{enumerate}
\end{corb}
\begin{pf}
By Proposition~\ref{CharRing} and Corollary~\ref{MatIdent}, we have
$$p_i = \frac{n_i^2}{\dim(A)} \sum_{j=1}^k \xi_i(\chi_j) \chi_{j^*}$$
Using the facts just proved, we therefore get
\begin{align*}
\sigma.p_i 
&= \frac{n_i^2}{\dim(A)} \sum_{j=1}^k \xi_i(\chi_j) \chi_{\sigma.j^*}
= \frac{n_i^2}{\dim(A)} \sum_{j=1}^k \xi_i(\chi_{\sigma^{-1}.j}) \chi_{j^*} \\
&= \frac{n_{\sigma^{-1}.i}^2}{\dim(A)} \sum_{j=1}^k \xi_{\sigma^{-1}.i}(\chi_{j}) \chi_{j^*}
= p_{\sigma^{-1}.i}
\end{align*}
which is the first assertion.

It suffices to check the second assertion on a basis, so that we can assume that
$\chi = \chi_j$. We then have by Corollary~\ref{VerlMat} that
\begin{align*}
\v_*(\sigma(\chi_j)) = \v_*(\chi_{\sigma.j}) 
&= \kappa n_{\sigma.j} \xi_{\sigma.j}(\chi_A)) p_{\sigma.j} \\
&= \kappa n_{j} \xi_{j}(\chi_A)) \sigma^{-1}.p_{j}
= \sigma^{-1}(\v_*(\chi_j))
\end{align*}
by the first assertion and the proof of the preceding lemma.
\qed
\end{pf}

This corollary shows in particular that the Galois group permutes the idempotents, which means that it acts via algebra automorphisms. This fact, however, is also obvious from the second description via the Wedderburn decomposition that we gave above.

\subsection[The semilinear actions]{} \label{SemiLin}
Besides the spaces~$\Ch(A)$ and~$\Ch_\Q(A)$ that we have considered above, we need to consider a third space that lies in between, namely the space
$\Ch_{\Q_N}(A)$, which we define to be the span of the irreducible characters with coefficients in the cyclotomic field~$\Q_N$. From the form of the base change matrix between the irreducible characters and the primitive idempotents of the character ring given in Proposition~\ref{CharRing}, we see that we could alternatively have defined it as the span of~$p_1,\ldots,p_k$ over the cyclotomic field~$\Q_N$. If we therefore, as in Paragraph~\ref{OrbitThm}, denote by~$Z_{\Q_N}(A)$ the span of the centrally primitive idempotents $e_1,\ldots,e_k$ with coefficients in the cyclotomic field, we have that 
$\Phi(\Ch_{\Q_N}(A))=Z_{\Q_N}(A)$.

We use these two different bases of the space to define two semilinear actions of the Galois group on~$\Ch_{\Q_N}(A)$ as follows:
\begin{defn}
For $\sigma \in \Gal(\Q_N/\Q)$, we define automorphisms
$\pi_*(\sigma)$ and~$\pi'_*(\sigma)$ of~$\Ch_{\Q_N}(A)$ by
$$\pi_*(\sigma)(\sum_{i=1}^k \lambda_i \chi_i)
:= \sum_{i=1}^k \sigma(\lambda_i) \chi_i \qquad \quad 
\pi'_*(\sigma)(\sum_{i=1}^k \lambda_i p_i)
:= \sum_{i=1}^k \sigma(\lambda_i) p_i$$
\end{defn}
In other words, $\pi_*(\sigma)$ acts on the coefficients in an expansion in  terms of the irreducible characters, and $\pi'_*(\sigma)$ acts on the coefficients in an expansion in  terms of the primitive idempotents.
Both of the automorphisms are semilinear in the sense that
$$\pi_*(\sigma)(\lambda \chi) = \sigma(\lambda)  \; \pi_*(\sigma)(\chi)
\qquad \qquad 
\pi'_*(\sigma)(\lambda \chi) = \sigma(\lambda) \; \pi'_*(\sigma)(\chi)$$
for $\lambda \in \Q_N$ and $\chi \in \Ch_{\Q_N}(A)$. 
Moreover, we have $\pi_*(\sigma)(\chi_i) = \chi_i$ as well as
$\pi'_*(\sigma)(p_i) = p_i$ for all $i=1,\ldots,k$.

The connection with the action of the Galois group considered in Paragraph~\ref{GaloisChar} is given by the following formula:
\begin{prop}
For all $\chi \in \Ch_{\Q_N}(A)$, we have
$$\sigma.\chi = (\pi'_*(\sigma) \circ \pi_*(\sigma)^{-1})(\chi)$$
Moreover, $\pi_*(\sigma)$ and~$\pi'_*(\tau)$ commute for all 
$\sigma, \tau \in \Gal(\Q_N/\Q)$. If $\chi_R(u)$ is rational, we have furthermore that
$\pi'_*(\sigma) = \v_* \circ \pi_*(\sigma) \circ \v_*^{-1}$.
\end{prop}
\begin{pf}
For the first assertion, note that 
$\pi'_*(\sigma) \circ \pi_*(\sigma)^{-1}$ is actually $\Q_N$-linear, so that it suffices to prove that
$\sigma.\chi_i = (\pi'_*(\sigma) \circ \pi_*(\sigma)^{-1})(\chi_i)$.
But we have
\begin{align*}
(\pi'_*(\sigma) \circ \pi_*(\sigma)^{-1})(\chi_i) &= 
\pi'_*(\sigma)(\chi_i) = \pi'_*(\sigma)(\sum_{j=1}^k \xi_j(\chi_i) p_j)\\
&= \sum_{j=1}^k \sigma(\xi_j(\chi_i)) p_j 
= \sum_{j=1}^k \xi_j(\chi_{\sigma.i}) p_j 
= \chi_{\sigma.i} = \sigma.\chi_i
\end{align*}
by Lemma~\ref{GaloisChar} and the discussion in Paragraph~\ref{CharRing}.

For the commutativity assertion, note that~$\pi_*(\sigma)$ obviously commutes with the action of~$\tau$, because the action of~$\tau$ is linear and permutes the characters. Also, it clearly commutes with~$\pi_*(\tau)$,
and therefore also with~$\pi'_*(\tau)$ by the result just proved.

For the third assertion, note that the assumption that $\chi_R(u)$ is rational implies that $\chi_R(u) = \gamma(\chi_R(u)) = \chi_R(u^{-1})$,
so that by our convention also~$\kappa=1/\chi_R(\ua)$ is rational.
To prove that $\pi'_*(\sigma)  \circ \v_* = \v_* \circ \pi_*(\sigma)$,
we also use that both sides are semilinear, so that it again suffices to check that both sides give the same result on~$\chi_i$.
But here we have by Corollary~\ref{VerlMat} that
\begin{align*}
\pi'_*(\sigma)(\v_*(\chi_i)) = 
\pi'_*(\sigma)(\kappa n_i \xi_i(\chi_A) p_i)
= \kappa n_i \xi_i(\chi_A) p_i 
= \v_*(\chi_i) = \v_*(\pi_*(\sigma)(\chi_i))
\end{align*}
which implies the assertion.
\qed
\end{pf}

We give a second proof of the fact that 
$\sigma.\chi = (\pi'_*(\sigma) \circ \pi_*(\sigma)^{-1})(\chi)$
from the point of view of the second construction of the action via the Wedderburn decomposition of the character ring, discussed after Lemma~\ref{GaloisChar}. From Proposition~\ref{CharRing}, we know that the primitive idempotents~$p_i$ are already contained in~$\Ch_{\Q_N}(A)$.
As we discussed above, a simple ideal of~$\Ch_\Q(A)$ is isomorphic to a subfield~$L$ of the cyclotomic field~$\Q_N$, and the action on the character ring restricts on these
Wedderburn components to the action of the Galois group. In other words, with respect to the isomorphism 
$\Ch_{\Q_N}(A) \cong \Ch_{\Q}(A) \o_\Q \Q_N$, we have
$$\sigma.(\chi \o \lambda) = \sigma(\chi) \o \lambda
\qquad \qquad
\pi_*(\sigma)(\chi \o \lambda) = \chi \o \sigma(\lambda)$$
for $\chi \in L$ and $\lambda \in \Q_N$. This shows that the formula that we have to prove 
is~$\pi'_*(\sigma)(\chi \o \lambda) = \sigma(\chi) \o \sigma(\lambda)$.

Because $L$ is a Galois extension of the rationals, the map
$$L \o_\Q \Q_N \rightarrow \Q_N^{\Gal(L/\Q)},~\chi \o \lambda \mapsto
(\sigma(\chi) \lambda)_{\sigma \in \Gal(L/\Q)}$$
is an algebra isomorphism,\endnote{\cite{FarbDennis}, Chap.~4, Exerc.~30, p.~140.}
where the right-hand side is an algebra with respect to componentwise multiplication. Therefore, for every~$\tau \in \Gal(\Q_N/\Q)$ there is 
a unique element $\sum_j l_j \o \lambda_j$ such that
$\sum_j \sigma(l_j) \lambda_j = \delta_{\sigma,\tau}$, corresponding to a primitive idempotent of~$\Q_N^{\Gal(L/\Q)}$. Because of its uniqueness,
we have 
$$\sum_j \sigma(l_j) \o \sigma(\lambda_j) = \sum_j l_j \o \lambda_j$$
for all $\sigma \in \Gal(\Q_N/\Q)$. But this shows that the endomorphism
$\chi \mapsto \pi_*(\sigma)(\sigma.\chi)$ of~$\Ch_{\Q_N}(A) $
is a semilinear map that preserves primitive idempotents, which is the defining property of~$\pi'_*(\sigma)$, establishing the assertion.

\subsection[The action on the center]{} \label{ActCent}
As we have $\Phi(\Ch(A)) = Z(A)$, we can use~$\Phi$ to transfer the action of the Galois group on the character ring to an action on the center. In other words, we define an action of~$\Gal(\Q_N/\Q)$ on~$Z(A)$ by requiring that the diagram
$$\square<1`1`1`1;1400`700>[\Ch(A)`\Ch(A)`Z(A)`Z(A);\chi \mapsto \sigma.\chi`\Phi`\Phi`z \mapsto \sigma.z]$$
is commutative. With respect to a smaller base field, we can also define
representations 
$\pi: \Gal(\Q_N/\Q) \rightarrow \GL(Z_{\Q_N}(A))$
and 
$\pi': \Gal(\Q_N/\Q) \rightarrow \GL(Z_{\Q_N}(A))$
by requiring that the diagrams
$$
\square<1`1`1`1;1000`700>[\Ch_{\Q_N}(A)`\Ch_{\Q_N}(A)`Z_{\Q_N}(A)`
Z_{\Q_N}(A);
\pi_*(\sigma)`\Phi`\Phi`\pi(\sigma)]
\qquad \qquad
\square<1`1`1`1;1000`700>[\Ch_{\Q_N}(A)`\Ch_{\Q_N}(A)`Z_{\Q_N}(A)`
Z_{\Q_N}(A);
\pi'_*(\sigma)`\Phi`\Phi`\pi'(\sigma)]
$$
commute. It is then a direct consequence of Proposition~\ref{SemiLin} that 
$$\sigma.z = (\pi'(\sigma) \circ \pi(\sigma)^{-1})(z)$$
Furthermore, $\pi(\sigma)$ and~$\pi'(\tau)$ commute for all 
$\sigma, \tau \in \Gal(\Q_N/\Q)$, and if~$\chi_R(u)$ is rational, we get from Proposition~\ref{RoleInt} that
$\pi'(\sigma) = \v \circ \pi(\sigma) \circ \v^{-1}$.
We can also deduce immediately from Proposition~\ref{VerlMat}, Corollary~\ref{GaloisChar}, and the equation $\Phi(p_i)=e_i$ that we have
$$\sigma.z_i = z_{\sigma.i} \qquad \qquad \sigma.e_i = e_{\sigma^{-1}.i}$$
Similarly, we have for the semilinear representations that
$$\pi(\sigma)(\sum_{i=1}^k \lambda_i z_i)
= \sum_{i=1}^k \sigma(\lambda_i) z_i \qquad \quad 
\pi'(\sigma)(\sum_{i=1}^k \lambda_i e_i)
= \sum_{i=1}^k \sigma(\lambda_i) e_i$$
for $\lambda_i \in \Q_N$.

Let us list some basic properties of this action:
\begin{prop}
For $\sigma \in \Gal(\Q_N/\Q)$, $\chi \in \Ch(A)$, and $z \in Z(A)$,
we have
\begin{enumerate}
\item 
$\iota(\sigma.z) = \sigma^{-1}.\iota(z) $

\item 
$\chi(\sigma.z) = (\sigma.\chi)(z)$

\item 
$\v(\sigma.z) = \sigma^{-1}.\v(z)$
\end{enumerate}
\end{prop}
\begin{pf}
For the first assertion, we have by Proposition~\ref{VerlMat} and Lemma~\ref{GaloisChar} that
$$\iota(\sigma.e_i) = \iota(e_{\sigma^{-1}.i}) =
\kappa n_{\sigma^{-1}.i} \chi_{\sigma^{-1}.i} =
\kappa n_{i} \sigma^{-1}.\chi_i = \sigma^{-1}.\iota(e_i) $$
For the second assertion, we can assume that $\chi = \chi_j$ and
$z=z_i$. We then have by Lemma~\ref{VerlMat} and Lemma~\ref{GaloisChar} that
\begin{align*}
n_i \chi_j(\sigma.z_i) = n_{\sigma.i} \chi_j(z_{\sigma.i})
&= n_{\sigma.i} \xi_{\sigma.i}(\chi_j) = s_{\sigma.i,j} \\
&= s_{i,\sigma.j} = n_i \xi_i(\chi_{\sigma.j})
= n_i \chi_{\sigma.j}(z_i)
\end{align*}
Alternatively, one can deduce this from the equation
$\chi_j(e_i) = n_i \delta_{ij}$. The third assertion follows from
Corollary~\ref{GaloisChar} by applying~$\Phi$ and using Proposition~\ref{RoleInt}.~\qed
\end{pf}

\subsection[Representations of the Drinfel'd double]{} \label{RepDouble}
Our next goal is to investigate the equivariance properties of the action of the Galois group with respect to the isomorphism~$\Psi$ introduced in Paragraph~\ref{DoubleQuasitri}. If~$V_i$ and~$V_j$ are any two of our simple $A$-modules, $V_i \o V_j$
can be considered as an $A \o A$-module  by the componentwise action.
As $\Psi$ is an algebra isomorphism between~$D(A)$ and~$A \o A$, we can introduce a $D(A)$-module structure on $V_i \o V_j$ by pullback via~$\Psi$.
We denote $V_i \o V_j$ by~$V_{ij}$ if endowed with this $D(A)$-module structure; note that this module was denoted by~$U$ in Paragraph~\ref{IndMult}. If~$\chi_{ij}$ denotes the character of~$V_{ij}$, we have $\chi_{ij} = \Psi^*(\chi_i \o \chi_j)$. Its degree is
$n_{ij} = n_i n_j$, and the corresponding centrally primitive idempotent is $e_{ij} = \Psi^{-1}(e_i \o e_j)$.

As described in Paragraph~\ref{CharRing}, from~$\chi_{ij}$ we can derive several additional quantities: the central characters
$$\omega_{ij}: Z(D(A)) \rightarrow K,~z \mapsto 
\frac{1}{n_{ij}} \chi_{ij}(z)$$
the idempotents of the character ring~$p_{ij} := \Phi^{-1}(e_{ij})$,
and the corresponding characters
$$\xi_{ij}: \Ch(D(A)) \rightarrow K,~\chi \mapsto \omega_{ij}(\Phi(\chi))$$
which in turn are used to define the class sums $z_{ij} \in Z(D(A))$ via the requirement that $\chi(z_{ij}) = \xi_{ij}(\chi)$. The following proposition describes how these quantities compare to the corresponding quantities for~$A$:
\begin{prop}
For $z \in Z(D(A))$ and $\chi,\chi' \in \Ch(A)$, we have
\begin{enumerate}
\item 
$\omega_{ij}(z) = (\omega_{i} \o \omega_{j})(\Psi(z))$

\item 
$p_{i,j^*} = \Psi^*(p_{i} \o p_{j})$

\item 
$\xi_{i,j^*}(\Psi^*(\chi \o \chi')) = \xi_{i}(\chi) \xi_{j}(\chi')$

\item 
$z_{i,j^*} = \Psi^{-1}(z_i \o z_j)$
\end{enumerate}
\end{prop}
\begin{pf}
The first assertion follows directly from the definitions; however, it should be noted that~$\Psi$ as an algebra isomorphism induces an isomorphism between the centers~$Z(D(A))$ and $Z(A \o A) \cong Z(A) \o Z(A)$.\endnote{\cite{FarbDennis}, Chap.~3, Exerc.~22, p.~103.}
The second assertion is equivalent to the equation
$\Phi(p_{i,j^*}) = \Phi(\Psi^*(p_{i} \o p_{j}))$, which by Proposition~\ref{DoubleTens} is equivalent to
$$e_{i,j^*} = \Psi^{-1}(\Phi(p_{i}) \o \sa(\Phi(p_{j}))) = 
\Psi^{-1}(e_{i} \o \sa(e_{j}))$$
which we have established above. Recall in this context that we have
$$\Ch((A \o A)_F) = \Ch(A \o A) \cong \Ch(A) \o \Ch(A)$$
by Lemma~\ref{DoubleTens}.

The third assertion follows from the second, because it suffices to check it in the case where $\chi = p_m$ and $\chi' = p_l$. But in view of the definition of the class sums, the third assertion can also be written as
$$\Psi^*(\chi \o \chi')(z_{i,j^*}) = \chi(z_{i}) \chi'(z_{j})$$
which implies that $\Psi(z_{i,j^*}) = z_{i} \o z_{j}$, which is the fourth assertion.
\qed
\end{pf}

\subsection[The equivariance of the isomorphism]{}  \label{EquivGalois}
In Paragraph~\ref{GaloisChar}, we have defined the action of the Galois group by first requiring that
$\xi_{\sigma.i}(\chi_j) = \sigma(\xi_{i}(\chi_j))$
and then defining it on~$\Ch(A)$ by setting 
$\sigma.\chi_i = \chi_{\sigma.i}$ and extending linearly. This action is also defined in exactly the same way for~$D(A)$; the only difference now is that we have indexed the corresponding quantities for~$D(A)$ by pairs. In terms of these pairs, the first equation above reads 
$\xi_{\sigma.(i,j)}(\chi_{ml}) = \sigma(\xi_{i,j}(\chi_{ml}))$.
But by Proposition~\ref{RepDouble}, this can be rewritten as 
\begin{align*}
\xi_{\sigma.(i,j)}(\chi_{ml}) &= \sigma(\xi_{i,j}(\chi_{ml}))
= \sigma(\xi_{i,j}(\Psi^*(\chi_{m} \o \chi_{l})))
= \sigma(\xi_{i}(\chi_{m}) \xi_{j^*}(\chi_{l})) \\
&= \sigma(\xi_{i}(\chi_{m})) \sigma(\xi_{j^*}(\chi_{l}))
= \xi_{\sigma.i}(\chi_{m}) \xi_{\sigma.j^*}(\chi_{l})
= \xi_{(\sigma.i,\sigma.j)}(\chi_{ml})
\end{align*}
where we have used the fact that $\sigma.(j^*) = (\sigma.j)^*$
discussed after Proposition~\ref{GaloisChar}. This means that we have
$\sigma.(i,j) = (\sigma.i,\sigma.j)$, which can be restated as follows:
\begin{prop}
The map
$$\Ch(A) \o \Ch(A) \rightarrow \Ch(D(A)),~\chi \o \chi' \mapsto 
\Psi^*(\chi \o \chi')$$
is $\Gal(\Q_N/\Q)$-equivariant if $\Ch(A) \o \Ch(A)$ is endowed with the diagonal action.
\end{prop}
\begin{pf}
For $\sigma \in \Gal(\Q_N/\Q)$, we have
$$\Psi^*(\sigma.\chi_i \o \sigma.\chi_j) =
\Psi^*(\chi_{\sigma.i} \o \chi_{\sigma.j}) 
= \chi_{\sigma.i,\sigma.j}= \chi_{\sigma.(i,j)}
= \sigma.\chi_{ij} = \sigma.\Psi^*(\chi_i \o \chi_j)$$
As the characters $\chi_i \o \chi_j$ form a basis of 
$\Ch(A) \o \Ch(A)$, this is sufficient.~\qed
\end{pf}

The preceding result can also be understood from the point of view of the Wedderburn decomposition of the character ring, as described after Lemma~\ref{GaloisChar}. As we pointed out there, the character rings~$\Ch_\Q(A)$ as well as~$\Ch_\Q(D(A))$ decompose into direct sums of subfields of the cyclotomic field~$\Q_N$, and the Galois group preserves the Wedderburn components and acts there via restriction to the corresponding subfield. Now~$\Psi$ is a Hopf algebra isomorphism between~$D(A)$ and~$(A \o A)_F$, so that~$\Psi^*$ restricts to an isomorphism between the character rings and therefore maps Wedderburn components to Wedderburn components. By Lemma~\ref{DoubleTens}, we have
$\Ch_\Q((A \o A)_F) = \Ch_\Q(A \o A)$, so that the assertion now will follow if we can justify that the isomorphism
$\Ch_\Q(A \o A) \cong \Ch_\Q(A) \o \Ch_\Q(A)$ is equivariant with respect to the diagonal action on the right-hand side.

This now follows from an argument that is similar to the one used 
at the end of Paragraph~\ref{SemiLin}. Let~$L$ and~$M$ be two subfields of~$\Q_N$ that appear as Wedderburn components of~$\Ch_\Q(A)$, and let~$P$ be a subfield of~$\Q_N$ that appears as a Wedderburn component
of~$\Ch_\Q(A \o A)$. We then have a commutative diagram of the form
$$\square<1`1`-1`1;1000`500>[L \o M`P`\Ch_\Q(A) \o \Ch_\Q(A)`\Ch_\Q(A \o A);f```]$$
where the left vertical arrow is the tensor product of the injections of the Wedderburn components, and the right vertical arrow is the projection
to the Wedderburn component. The resulting multiplicative map
$f: L \o M \rightarrow  P$ may be zero, in which case it is equivariant.
If it is not zero, then $f(1 \o 1)$ is a nonzero idempotent in~$P$, which implies that $f(1 \o 1)=1$. Then the map
$$L \rightarrow P,~x \mapsto f(x \o 1)$$
is a field homomorphism, and since $\Gal(\Q_N/\Q)$ preserves all subfields and acts on them in a way that is independent of the embedding into~$\Q_N$, we get that 
$$\sigma(f(x \o 1)) = f(\sigma(x) \o 1)$$
for all $\sigma \in \Gal(\Q_N/\Q)$. Similarly, we 
$\sigma(f(1 \o y)) = f(1 \o \sigma(y))$
and therefore
\begin{align*}
\sigma(f(x \o y)) &= \sigma(f(x \o 1)f(1 \o y)) = 
\sigma(f(x \o 1)) \sigma(f(1 \o y)) \\
&= f(\sigma(x) \o 1)f(1 \o \sigma(y)) = f(\sigma(x) \o \sigma(y))
\end{align*}
Pasting all the Wedderburn components together, we see that 
the isomorphism 
$\Ch_\Q(A \o A) \cong \Ch_\Q(A) \o \Ch_\Q(A)$ is equivariant with respect to the diagonal action on the right-hand side.

In Paragraph~\ref{ActCent}, we have transferred the action of the Galois group from the character ring~$\Ch(A)$ to the center~$Z(A)$ by requiring that~$\Phi$ be equivariant. As~$D(A)$ is also a semisimple factorizable Hopf algebra, this whole discussion applies to~$D(A)$ as well, so that we
also have an action of~$\Gal(\Q_N/\Q)$ on~$Z(D(A))$. The formulas obtained in Paragraph~\ref{ActCent} then give in particular that
$$\sigma.z_{ij} = z_{\sigma.i, \sigma.j} \qquad  \qquad
\sigma.e_{ij} = e_{\sigma^{-1}.i, \sigma^{-1}.j}$$
where we have used the formula $\sigma.(i,j)=(\sigma.i,\sigma.j)$
obtained earlier.

Now we have already pointed out in the proof of Proposition~\ref{RepDouble} that~$\Psi$ induces an isomorphism between~$Z(D(A))$ and $Z(A) \o Z(A)$. The above formulas now imply that this isomorphism is equivariant if we endow $Z(A) \o Z(A)$ with the diagonal action of the Galois group:
\begin{corollary}
For $\sigma \in \Gal(\Q_N/\Q)$ and $z \in Z(D(A))$, we have
$\Psi(\sigma.z) = \sigma.\Psi(z)$.
\end{corollary}
\begin{pf}
It suffices to check that we have
$$\Psi(\sigma.e_{ij}) = \Psi(e_{\sigma^{-1}.i,\sigma^{-1}.j})
= e_{\sigma^{-1}.i} \o e_{\sigma^{-1}.j}
= \sigma.e_{i} \o \sigma.e_{j} = \sigma.\Psi(e_{ij})
$$
since the centrally primitive idempotents form a basis of~$Z(D(A))$.
\qed
\end{pf}

\newpage
\section{Galois groups and indicators} \label{Sec:GalInd}
\subsection[A digression on Frobenius algebras]{} \label{DigFrob}
Before we really begin, we present a little background from the theory of Frobenius algebras. We therefore defer the discussion of the setup of this section to Paragraph~\ref{InvInd}.

Recall that the ring $M(r \times r, K)$ of $r \times r$-matrices is a Frobenius algebra with respect to the ordinary matrix trace function~$\Tr$ as Frobenius homomorphism. The dual basis of matrix units~$E_{ij}$ with respect to the bilinear form arising from the trace is again formed by the matrix units~$E_{ji}$, with the indices reversed. The corresponding Casimir element therefore is 
$$\sum_{i,j=1}^r  E_{ij} \o E_{ji}$$
Note that this element is symmetric under interchange of the tensorands.
The following lemma states that this and its Casimir property characterize it up to proportionality: 
\begin{lemma}
Suppose that $\sum_{i}  A_{i} \o B_{i} \in M(r \times r, K)^{\o 2}$
satisfies
\begin{enumerate}
\item 
$\sum_{i} A A_{i} \o B_{i} = \sum_{i} A_{i} \o B_{i} A$

\item 
$\sum_{i}  A_{i} \o B_{i} = \sum_{i}  B_{i} \o A_{i}$
\end{enumerate}
Then there is a number $\mu \in K$ such that
$\sum_{i}  A_{i} \o B_{i} = \mu \sum_{i,j=1}^r  E_{ij} \o E_{ji}$.
\end{lemma}
\begin{pf}
This verification is left to the reader.
\qed
\end{pf}

If we multiply the tensorands of our Casimir element together, we get a multiple of the unit matrix~$E_r$:
$$\sum_{i,j=1}^r  E_{ij} E_{ji} = r E_r =
r \sum_{i,j=1}^r  \Tr(E_{ij}) E_{ji}$$
Multiplying this equation by~$\mu$, we see that an element of the form considered in the lemma will also satisfy this equation:
$$\sum_{i}  A_{i} B_{i} = r \sum_{i}  \Tr(A_{i}) B_{i}$$
Note that this discussion applies directly to the Wedderburn components of an arbitrary semisimple algebra, where, however, the number~$r$ varies with the Wedderburn component. The algebra that we have in mind is the 
character ring~$\Ch(H)$ of a semisimple Hopf algebra, which is a Frobenius algebra.\endnote{\cite{So4}, Prop.~3.5, p.~211.} In this application, the element~$A_i$ that appears in the lemma will be the Wedderburn component of an irreducible character, and~$B_i$ will be the Wedderburn component of the corresponding dual character, so that $\sum_{i}  A_{i} B_{i}$ is
the Wedderburn component of the character of the adjoint representation.

\subsection[The invariance of the induced trivial module]{} \label{InvInd}
It is the aim of this section to discuss how the action of the Galois group relates to the action of the modular group, and again the equivariant Frobenius-Schur indicators that we introduced in Paragraph~\ref{EquiFrobSchur} will be our main tool. We assume throughout the section that $A=D:=D(H)$, the Drinfel'd double of a semisimple Hopf algebra~$H$. Recall from Paragraph~\ref{RoleEval} that $\chi_R(\ud)=\dim(H)$ is a rational number in this case; furthermore, as pointed out in Paragraph~\ref{CongDrinfDoubl}, $D$ and~$H$ have the same exponent~$N$. We will use the notation of Section~\ref{Sec:GaloisGroup} throughout.

We need a preparatory result about the induced module of the trivial module. As discussed in Paragraph~\ref{Induct}, the induced $D$-module
$\Ind(K)$ of the trivial $H$-module~$K$ can be realized on the underlying vector space~$H^*$, and this is the way in which we will look at it in this paragraph. If we denote its character by~$\eta$, the result that we will need is the following:
\begin{prop}
For all $z \in Z(D)$ and all $\sigma \in \Gal(\Q_N/\Q)$, we have
$$\eta(\sigma.z) = \eta(z)$$
\end{prop}
\begin{pf}
\begin{list}{(\arabic{num})}{\usecounter{num} \leftmargin0cm \itemindent5pt}
\item
By linearity, it suffices to check this on a basis of the center, so that we can assume that $z=e_i$ is a centrally primitive idempotent. If~$m_i$ denotes the multiplicity of the simple module~$V_i$ in~$\Ind(K)$, we have
$\eta(e_i) = n_i m_i$. We first treat the case where $\eta(e_i) \neq 0$; i.e., the case in which~$V_i$ is really a constituent of~$\Ind(K)$.

Recall\endnote{\cite{YYY2}, Par.~6.3, p.~45.} that
$$\Ch(H) \rightarrow \End_D(H^*),~\chi \mapsto 
(\varphi \mapsto \varphi \chi)$$
is an algebra anti-isomorphism. Since~$e_i$ is central, the action of~$e_i$ is also a $D$-endomorphism of~$H^*$, and therefore must be given by right multiplication with a centrally primitive idempotent~$q \in Z(\Ch(H))$. Further discussion\endnote{\cite{YYY2}, Par.~6.3, p.~46.}
shows that~$q :=\Res(p_i)$, where 
$$\Res: \Ch(D) \rightarrow Z(\Ch(H))$$
is the restriction map. The fact that $e_i.\varphi = \varphi q$ for all
$\varphi \in H^*$ implies in particular that $\eta(e_i) = \dim(H^*q)$.
Furthermore, the general theory of endomorphism rings of semisimple modules\endnote{\cite{FarbDennis}, Prop.~1.8, p.~36.}
implies that the Wedderburn component of~$\Ch(H)$ corresponding to~$q$ has dimension~$\dim(\Ch(H)q) = m_i^2$.

\item
Let
$$\xi: Z(\Ch(H)) \rightarrow K$$
be the central character corresponding to~$q$; i.e., the algebra homomorphism from~$Z(\Ch(H))$ to~$K$ that maps~$q$ to~$1$ and vanishes on all other centrally primitive idempotents of~$\Ch(H)$. Note that we then have $\xi_i = \xi \circ \Res$. Now we know from Lorenz' proof of the class equation,\endnote{\cite{Lo2}, p.~2843.} combined with the discussion at the end of Paragraph~\ref{DigFrob}, that 
$$\frac{m_i \dim(H)}{\dim(H^* q)} = \frac{\xi(\chi'_A)}{m_i}$$
where~$\chi'_A$ denotes the character of the adjoint representation of~$H$.
As we have $\dim(H^* q) = \eta(e_i) = n_i m_i$, we can rewrite this as
$$m_i \dim(H) = n_i \xi(\chi'_A)$$

\item
If $\chi \in \Ch_\Q(D)$, we have
\begin{align*}
\xi_{\sigma.i}(\chi) = \sigma(\xi_{i}(\chi)) =
\sigma(\xi(\Res(\chi)))
\end{align*}
This shows that $\xi_{\sigma.i} = \xi' \circ \Res$, where
$\xi': \Ch(H) \rightarrow K$ arises by scalar extension of 
$\sigma \circ \xi$ from~$\Ch_\Q(H)$ to~$\Ch(H)$. 
If $q' \in Z(\Ch(H))$ is the centrally primitive idempotent that corresponds to~$\xi'$, in the sense that $\xi'(q')=1$, it follows exactly as above that $\eta(e_{\sigma.i}) = \dim(H^* q') \neq 0$. Even stronger,
by applying the equation obtained in the preceding step to these
idempotents and using Lemma~\ref{GaloisChar}, we get 
\begin{align*}
m_{\sigma.i} \dim(H) = n_{\sigma.i} \xi'(\chi'_A)
= n_{i} \sigma(\xi(\chi'_A)) 
= n_{i} \xi(\chi'_A) = m_{i} \dim(H)
\end{align*}
from which our assertion follows immediately, as we have
$\eta(e_{\sigma.i}) = n_{\sigma.i} m_{\sigma.i} = n_i m_i = \eta(e_i)$
by Lemma~\ref{GaloisChar} again.

\vspace{1mm}
\item
Finally, it remains to consider the case $\eta(e_i) = 0$. But then we must also have $\eta(e_{\sigma.i}) = 0$, because otherwise we can replace~$i$ by~$\sigma.i$ and~$\sigma$ by its inverse in the preceding discussion to get~$\eta(e_i) = \eta(e_{\sigma^{-1}\sigma.i}) \neq 0$, which is obviously a contradiction.
\qed
\end{list}
\end{pf}
\vspace{1mm}
In view of Proposition~\ref{ActCent}, this result can be restated by saying that the character~$\eta$ of the induced trivial module is invariant under the action of the Galois group:
\begin{corollary}
For all $\sigma \in \Gal(\Q_N/\Q)$, we have $\sigma.\eta = \eta$.
\end{corollary}

\subsection[The action and the indicators]{} \label{GaloisInd}
We now want to relate the action of the Galois group to the equivariant Frobenius-Schur indicators that we introduced in Paragraph~\ref{EquiFrobSchur}. Recall that for any cyclotomic field~$\Q_m$ and an integer~$q$ relatively prime to~$m$, we have an 
automorphism~$\sigma_q \in \Gal(\Q_m/\Q)$ with the property that 
$\sigma_q(\zeta) = \zeta^q$ for every $m$-th root of unity~$\zeta$. 
Every element of the Galois group is of this form. Although $\sigma_q$
depends on the field and therefore on~$m$, this dependence is diminished by the fact that when~$m'$ divides~$m$, so that~$\Q_{m'} \subset \Q_m$
and~$q$ is also relatively prime to~$m'$, the restriction of~$\sigma_q$
to~$\Q_{m'}$ coincides with the automorphism defined for this field.

Using this notation, we can now relate the action of the Galois group to the equivariant Frobenius-Schur indicators in the following way:
\pagebreak
\begin{prop}
Consider an $H$-module~$V$, a central element $z \in Z_{\Q_N}(D)$, and three integers $m,l,q \in \Z$.
\begin{enumerate}
\item 
If $q$ is relatively prime to~$N$ and~$m$, we have
$$\sigma_q(I_V((m,l),z)) = I_V((m,lq),\pi'(\sigma_q)(z))$$

\item 
If $q$ is relatively prime to~$N$ and~$l$, we have
$$\sigma_q(I_V((m,l),z)) = I_V((mq,l),\pi(\sigma_q)(z))$$

\item 
If $q$ is relatively prime to~$N$, $m$, and~$l$, we have
$$I_V((m,lq), \sigma_q.z) = I_V((mq,l),z)$$
\end{enumerate}
\end{prop}
\begin{pf}
\begin{list}{(\arabic{num})}{\usecounter{num} \leftmargin0cm \itemindent5pt}
\item
Recall that, by Proposition~\ref{OrbitThm}, we have 
$I_V((m,l),z) \in \Q_N$, so that the expressions considered are well-defined. We begin by proving the first assertion in the case $m>0$. 
For this, we note that both sides of the equation are semilinear in the variable~$z$, so that we can assume that
$z=e_i$ for some~$i$.  Then we have by definition that $I_V((m,l),z) = \Tr(\beta^l \circ \rho_m(e_i))$, 
where $\beta = \beta_{V,V^{\o (m-1)}}$, properly interpreted in the case $m=1$. The endomorphism $\beta^l \circ \rho_m(e_i)$ coincides with~$\beta^l$ on the isotypical component corresponding to~$e_i$ and is zero otherwise. Since $\beta^{mN} = \id$, we see that the eigenvalues of
$\beta^l \circ \rho_m(e_i)$ are $mN$-th roots of unity, which are raised to their $q$-th power by the action of~$\sigma_q$. But these $q$-th powers are exactly the eigenvalues of~$\beta^{lq} \circ \rho_m(e_i)$, which implies the first assertion in the case~$m>0$. 

\item
The first assertion in the case $m<0$ follows from the case $m>0$, because we then have by definition that 
\begin{align*}
\sigma_q(I_V((m,l),z)) &= \sigma_q(I_V((-m,-l),\sd(z))) 
= I_V((-m,-lq),\pi'(\sigma_q)(\sd(z))) \\
&= I_V((-m,-lq),\sd(\pi'(\sigma_q)(z)))
= I_V((m,lq),\pi'(\sigma_q)(z))
\end{align*}
where the fact that $\pi'(\sigma_q)$ and~$\sd$ commute follows from the fact that~$\sd$ permutes the centrally primitive idempotents. The remaining case of the first assertion therefore is the case $m=0$; however, we leave this case open for a moment.

\item
Instead, we now prove the second assertion in the case $m>0$ and $q=-1$.
In this case, we have that $\sigma_{-1} = \gamma$ is the restriction of complex conjugation. In view of the assertion already established, we have to show that  
$$I_V((m,-l),\pi'(\gamma)(z)) = I_V((-m,l),\pi(\gamma)(z))$$
Replacing $z$ by $\pi(\gamma)^{-1}(z)$, we get the equivalent assertion that
$$I_V((m,-l),\gamma.z) = I_V((m,-l),(\pi'(\gamma) \circ \pi(\gamma)^{-1})(z)) = I_V((-m,l),z)$$
which follows from the fact that $\gamma.z = \sa(z)$ by Proposition~\ref{GaloisChar}.

\item
From this, we now deduce the first assertion in the case $m=0$ and~$l>0$.
The condition that $q$ is relatively prime to~$m=0$ then forces that
$q=\pm 1$. As the case $q=1$ is obvious, we can assume that $q=-1$, and this reduces to the case just treated since
\begin{align*}
I_V((0,-l),\pi'(\sigma_{-1})(z)) &= 
I_V((0,-l),(\v \circ \pi(\sigma_{-1}) \circ \v^{-1})(z)) \\
&= I_V((-l,0),(\pi(\sigma_{-1}) \circ \v^{-1})(z)) \\
&= \sigma_{-1}(I_V((l,0),\v^{-1}(z)))= \sigma_{-1}(I_V((0,l),z))
\end{align*}

\item
As in the second step, the case of the first assertion where 
$m=0$ and~$l<0$ follows from the case $m=0$ and~$l>0$ just established
by using the antipode. It therefore remains to establish the first assertion in the case where $m=l=0$. Exploiting semilinearity as in the first step, we can again assume that~$z=e_i$. Note that in general 
$I_V((0,0),z)$ is the trace of the action of~$z$ on the induced module
$\Ind(K)$; in case $z=e_i$, this is an integer. It is therefore invariant under every Galois automorphism, which establishes the first assertion in this case and therefore completely.

\item
The second assertion follows from the first by a variant of the one we have used in the fourth step:
\begin{align*}
I_V((mq,l),\pi(\sigma_q)(z)) &= 
I_V((mq,l),(\v^{-1} \circ \pi'(\sigma_q) \circ \v)(z)) \\
&= I_V((-l,mq),(\pi'(\sigma_q) \circ \v)(z)) \\
&= \sigma_q(I_V((-l,m),\v(z)))= \sigma_q(I_V((m,l),z))
\end{align*}

\item
By comparing the first and the second assertion, we get that
$$I_V((m,lq),\pi'(\sigma_q)(z)) = I_V((mq,l),\pi(\sigma_q)(z))$$
Replacing~$z$ by~$\pi(\sigma_q)^{-1}(z)$, this becomes
$$I_V((m,lq), \sigma_q.z)) = I_V((m,lq),\pi'(\sigma_q)(\pi(\sigma_q)^{-1}(z))) = I_V((mq,l),z)$$
which is the third assertion.
\qed
\end{list}
\end{pf}

For an integer~$q$ that is relatively prime to~$N$, we can find another integer~$q'$ such that $qq' \equiv 1 \pmod{N}$, which describes the inverse of the residue class of~$q$ in the group of units~$\Z_N^\times$.
Using it, we can derive the following corollary, which should be viewed as a kind of adjunction relation between the Galois action and an action on the lattice points:
\begin{corollary}
Suppose that~$m$ and~$l$ are nonzero integers. Furthermore, suppose that~$q$ and~$q'$ are relatively prime integers that are both relatively prime to~$ml$. If $qq' \equiv 1 \pmod{N}$, we have
$$I_V((m,l),\sigma_{q}.z) = I_V((mq,lq'),z)$$
\end{corollary}
\begin{pf}
Since~$q$ is relatively prime to~$lq'$, it follows from the preceding proposition that
$$I_V((m,lqq'),\sigma_{q}.z) =  I_V((mq,lq'),z)$$
As we have $qq' \equiv 1 \pmod{N}$ and 
$\gcd(m,lqq') = \gcd(m,l)$, it follows from Proposition~\ref{CongSubgr} that~$(m,lqq')$ and~$(m,l)$ are in the same $\Gamma(N)$-orbit, so that the assertion follows from Theorem~\ref{OrbitThm}.
\qed
\end{pf}

\subsection[Diagonal matrices]{} \label{DiagMat}
For integers~$q$ and~$q'$ such that $qq' \equiv 1 \pmod{N}$ as above, we denote the residue classes in~$\Z_N$ by~$\bar{q}$ resp.~$\bar{q}'$. With these numbers, we have associated in Paragraph~\ref{PresFac} the matrix
$$\gd(q) = \begin{pmatrix} \bar{q} & 0 \\ 0 & \bar{q}' \end{pmatrix} 
\in \SL(2,\Z_N)$$
Because the principal congruence subgroup~$\Gamma(N)$ acts trivially by Theorem~\ref{CongDrinfDoubl}, the action of the modular group factors over
the quotient group~$\SL(2,\Z_N)$, so that in particular the action of~$\gd(q)$ on the center is defined. It has the following basic property:
\begin{prop}
$I_V((m,l),\sigma_q.z) = I_V((m,l),\gd(q).z)$
\end{prop}
\begin{pf}
\begin{list}{(\arabic{num})}{\usecounter{num} \leftmargin0cm \itemindent5pt}
\item
We first prove this in the case where both~$m$ and~$l$ are nonzero.
The numbers~$q$, $ml\neq 0$, and~$N$ are relatively prime, because already~$q$ and~$N$ are relatively prime, and therefore we get from Lemma~\ref{OrbCongRel} that there is an integer~$c$ such that 
$q+cN$ is relatively prime to~$ml$. Note that~$q+cN$ is necessarily nonzero. As the asserted equation only depends on the residue class of~$q$ modulo~$N$, we can replace~$q$ by~$q+cN$
if necessary to achieve that~$q$ is relatively prime to~$ml$.

Similarly, since~$q'$, $mlq \neq 0$, and~$N$ are relatively prime, we can again by Lemma~\ref{OrbCongRel} find an 
integer~$c'$ such that~$q'+c'N$ is relatively prime to~$mlq$. If necessary, we can replace~$q'$ by~$q'+c'N$ to achieve that on the one hand~$q$ and~$q'$ are relatively prime, on the other hand both of them
are relatively prime to~$ml\neq 0$.

\item
If~$q$ and~$q'$ are chosen so that they have these additional properties, it follows from Corollary~\ref{GaloisInd} that
$I_V((m,l),\sigma_{q}.z) = I_V((mq,lq'),z)$. Therefore,
our claim will follow if we can establish that
$$I_V((mq,lq'),z) = I_V((m,l),\gd(q).z)$$
For this, suppose that 
$\begin{pmatrix} a & b \\ c & d \end{pmatrix} \in \SL(2,\Z)$
is a lift 
of~$\gd(q) \in \SL(2,\Z_N)$. By Theorem~\ref{EquivarThm} and Theorem~\ref{OrbitThm},
it then suffices to show that $(mq,lq')$ and 
$(am + cl, bm +dl)$ are in the same $\Gamma(N)$-orbit. But this follows again from Proposition~\ref{CongSubgr}, as we have
$$t:= \gcd(mq,lq') = \gcd(m,l) = \gcd(am + cl, bm +dl)$$
and the two lattice points are by construction componentwise congruent modulo~$N$ after division by~$t$.

\item
The case $m \neq 0$, $l=0$ can be reduced to the case just treated by observing that the lattice points $(m,0)$ and
$$(m, mN) = (m,0)\begin{pmatrix} 1 & N \\ 0 & 1 \end{pmatrix}$$
are in the same $\Gamma(N)$-orbit, so that we get
\begin{align*}
I_V((m,0),\sigma_q.z) &= I_V((m,mN),\sigma_q.z) = 
I_V((m,mN),\gd(q).z) \\
&= I_V((m,0),\gd(q).z)
\end{align*}
by Theorem~\ref{OrbitThm}. Similarly, the 
case $m = 0$, $l \neq 0$ can be reduced to the previous case, as the lattice points $(0,l)$ and
$$(lN, l) = (0,l)\begin{pmatrix} 1 & 0 \\ N & 1 \end{pmatrix}$$
are in the same $\Gamma(N)$-orbit.

\item
It remains to consider the case $m=l=0$. In this case, we have in view of Theorem~\ref{EquivarThm} that
$I_V((0,0),\gd(q).z) = I_V((0,0),z)$
But on the other hand, this indicator is by definition equal to the character of the induced trivial module, so that the
equation $I_V((0,0),\sigma_q.z) = I_V((0,0),z)$ is exactly Proposition~\ref{InvInd}.~\qed
\end{list}
\end{pf}

A special case of this proposition leads to an invariance property that will be important later. Recall from Paragraph~\ref{CharRing} that~$e_1$ is a normalized integral. We now show that it is invariant under the
action of the diagonal matrices~$\gd(q)$ considered above. Before we derive this, note a difference in the dualization of the Galois group and the modular group: While the action of the Galois group was carried over from the character ring to the center in Paragraph~\ref{ActCent} by requiring that~$\Phi$ is equivariant, the action of the modular group on the character ring was introduced in Paragraph~\ref{DualRepMod} by regarding the character ring as dual to the center via the canonical pairing.
\begin{corollary}
If $q \in \Z$ is relatively prime to~$N$, we have
$\gd(q).1 = 1$ as well as $\gd(q).e_1 = e_1$.
\end{corollary}
\begin{pf}
As we have already used in the proof of Corollary~\ref{OrbitThm}, the indicators corresponding to the lattice points~$(1,0)$ are exactly the characters of the induced modules. The regular representation of~$D$ is induced from the regular representation of~$H$, so that we get as a special case of the above proposition that
$\chi_R(\sigma_q.z) = \chi_R(\gd(q).z)$.
If~$q' \in \Z$ satisfies $qq' \equiv 1 \pmod{N}$,
this equation, in terms of the action of the modular group on the character ring introduced in Paragraph~\ref{DualRepMod}, reads
$$\sigma_q.\chi_R = \gd(q').\chi_R$$
where we have used Proposition~\ref{ActCent} in addition. 
As we have $\iota(\D) = \kappa \chi_R$ by construction, it follows from Proposition~\ref{DualRepMod} and Proposition~\ref{ActCent} that
$\sigma_q^{-1}.\D = \gd(q').1$.
Now it follows from Lemma~\ref{GaloisChar} that $\sigma_q.e_1 = e_1$. 
Applying~$\v$ and using Proposition~\ref{ActCent}, we get
$\sigma_q^{-1}.\v(e_1) = \v(e_1)$. But $\v(e_1) = \kappa \D$ by Corollary~\ref{VerlMat}, which shows that $\sigma_q^{-1}.\D = \D$ and establishes the first assertion. This equation also shows that
the second assertion follows from the first by applying~$\v^{-1}$ and using
the commutation relation
\begin{align*}
\begin{pmatrix} \bar{q} & 0 \\ 0 & \bar{q}' \end{pmatrix}
\begin{pmatrix} 0 & -1 \\ 1 & 0 \end{pmatrix}
= \begin{pmatrix} 0 & -\bar{q} \\ \bar{q}' & 0 \end{pmatrix}
= \begin{pmatrix} 0 & -1 \\ 1 & 0 \end{pmatrix}
\begin{pmatrix} \bar{q}' & 0 \\ 0 & \bar{q} \end{pmatrix}
\end{align*}
between~$\gv$ and the diagonal matrix.
\qed
\end{pf}

\subsection[The Galois group and the modular group]{} \label{GalMod}
Proposition~\ref{DiagMat} can be substantially strengthened: The following theorem, which is the main result of this section, asserts that the two actions do not only give the same result inside the indicators, but are just equal:
\begin{thm}
If~$q \in \Z$ is relatively prime to~$N$, then we have
$\sigma_q.z = \gd(q).z$
for all elements~$z$ in the center of~$D=D(H)$.
\end{thm}
\begin{pf}
For this, it suffices to show that $\chi(\sigma_q.z) = \chi(\gd(q).z)$
for all $\chi \in \Ch(D)$. If we rewrite this equation in terms of the action of the Galois group on the character ring introduced in Paragraph~\ref{GaloisChar} and of the action of the modular group on the character ring introduced in Paragraph~\ref{DualRepMod}, it takes by Proposition~\ref{ActCent} the form
$\sigma_q.\chi =  \gd(q').\chi$,
where $q' \in \Z$ satisfies $qq' \equiv 1 \pmod{N}$. This is equivalent to the condition
$$\langle \sigma_q.\chi, \chi_l \rangle_* = 
\langle \gd(q').\chi, \chi_l \rangle_*$$
for all~$l$, using the bilinear form introduced in Paragraph~\ref{IndMult}. If~$\eta_l$ denotes the character of the induced $D(D)$-module~$\Ind(V_l)$, we can by Proposition~\ref{IndMult}
rewrite this equation in the equivalent form
$$\eta_l(\Psi^{-1}(\iota^{-1}(\sigma_q.\chi) \o \iota^{-1}(\ea))) = 
\eta_l(\Psi^{-1}(\iota^{-1}(\gd(q').\chi) \o \iota^{-1}(\ea)))$$
As we have $\iota^{-1}(\ea) = \frac{1}{\kappa} e_1$ by Proposition~\ref{VerlMat}, we can rewrite this equation further in the form
$$\eta_l(\Psi^{-1}(\sigma_q^{-1}.\iota^{-1}(\chi) \o e_1)) = 
\eta_l(\Psi^{-1}(\gd(q').\iota^{-1}(\chi) \o e_1))$$
where we have used Proposition~\ref{DualRepMod} and Proposition~\ref{ActCent}. By Lemma~\ref{GaloisChar}, we have
$\sigma_q(e_1) = e_1$, which can be used together with the analogous equation in Corollary~\ref{DiagMat} to give with
$$\eta_l(\Psi^{-1}(\sigma_q^{-1}.\iota^{-1}(\chi) \o \sigma_q^{-1}.e_1)) = 
\eta_l(\Psi^{-1}(\gd(q').\iota^{-1}(\chi) \o \gd(q').e_1))$$
still another equivalent version of our condition. This in turn is by Corollary~\ref{EquivGalois} and Proposition~\ref{ModGrDoub} equivalent to
$$\eta_l(\sigma_q^{-1}.\Psi^{-1}(\iota^{-1}(\chi) \o e_1)) = 
\eta_l(\gd(q').\Psi^{-1}(\iota^{-1}(\chi) \o e_1))$$
But as we have $\eta_l(z') = I_{V_l}((1,0),z')$ for all $z' \in Z(D(D))$,
this is a special case of Proposition~\ref{DiagMat}.
\qed
\end{pf}

This theorem has an interesting consequence for the components~$u_i$ of the Drinfel'd element that we introduced in Paragraph~\ref{VerlMat}:
\begin{corollary}
For all $\sigma \in \Gal(\Q_N/\Q)$, we have
$\sigma^2(u_i) = u_{\sigma.i}$. 
\end{corollary}
\begin{pf}
Recall that all component~$u_i$ of the Drinfel'd element are $N$-th roots of unity. If $\sigma = \sigma_q$ and $qq' \equiv 1 \pmod{N}$, it follows from the relation
$$\begin{pmatrix} \bar{q}' & 0 \\ 0 & \bar{q} \end{pmatrix}
\begin{pmatrix} 1 & 1 \\ 0 & 1 \end{pmatrix}
\begin{pmatrix} \bar{q} & 0 \\ 0 & \bar{q}' \end{pmatrix}.e_i
=\begin{pmatrix} 1 & \bar{q}'^2 \\ 0 & 1 \end{pmatrix}.e_i$$
by the preceding theorem that 
$\sigma_{q}^{-1}.(u^{-1}\sigma_q.e_i) = u^{-q'^2} e_i$,
or alternatively 
$$\frac{1}{u_{\sigma_q^{-1}.i}} e_{i} = \sigma_{q}^{-1}.(\frac{1}{u_{\sigma_q^{-1}.i}}e_{\sigma_q^{-1}.i}) = \sigma_{q}^{-1}.(u^{-1}e_{\sigma_q^{-1}.i}) = u_i^{-q'^2} e_i
= \sigma_q^{-2}(\frac{1}{u_i}) e_i$$
which establishes the assertion.
\qed
\end{pf}

\newpage
\section{Galois groups and congruence subgroups} \label{Sec:GalCong}
\subsection[The Hopf symbol]{} \label{HopfSymb}
In this section, we leave the case of a Drinfel'd double and reconsider an arbitrary semisimple factorizable Hopf algebra~$A$. We will use the notation introduced in Section~\ref{Sec:FaktSemisim}, and for the double~$D:=D(A)$ we will use the notation introduced in Paragraph~\ref{RepDouble}. We first consider a new quantity that will play an important role in the sequel:
\begin{defn}
For~$q \in \Z$ , we define the Hopf symbol
$$\jac{q}{A} := 
\begin{cases}
\dfrac{\chi_R(u^{-q})}{\chi_R(u^{-1})} &: \gcd(q,N) = 1 \\
0 &: \gcd(q,N) \neq 1
\end{cases}$$
where $N$ is the exponent of~$A$.
\end{defn}

The Hopf symbol generalizes the Jacobi symbol:
If we take Radford's example, i.e., $A=K[\Z_n]$, the group ring of a cyclic group of odd order~$n$ endowed with a modified R-matrix, then it follows from the discussion in Paragraph~\ref{ExamRadf} that
$$\jac{q}{K[\Z_n]} = \jac{q}{n}$$

The Hopf symbol should be viewed as a 1-cocycle
in the following way:\endnote{\cite{Serre2}, Chap.~VII, \S~3, p.~113.}
The Galois group~$\Gal(\Q_N/\Q)$ acts on
the multiplicative group~$\Q_N^\times$ of nonzero elements in the cyclotomic field. With the element $\chi_R(u^{-1}) \in \Q_N^\times$,
we can therefore associate a 1-coboundary
$$f: \Gal(\Q_N/\Q) \rightarrow \Q_N^\times,~\sigma 
\mapsto \frac{\sigma(\chi_R(u^{-1}))}{\chi_R(u^{-1})}$$
This is essentially the Hopf symbol: If $q$ is relatively prime to~$N$, we have
$$f(\sigma_q) = \frac{\sigma_q(\chi_R(u^{-1}))}{\chi_R(u^{-1})}
= \frac{\chi_R(u^{-q})}{\chi_R(u^{-1})} = \jac{q}{A}$$
As a 1-coboundary, it is in particular a 1-cocycle, and therefore satisfies
$$\jac{qq'}{A} = \jac{q}{A} \sigma_q(\jac{q'}{A})$$
Actually, it follows from a version of Hilbert's theorem~90 that every cocycle is a coboundary in this situation.\endnote{\cite{Serre2}, Chap.~X, \S~1, Prop.~2, p.~150.}

From the cocycle equation, we immediately obtain the following:
\pagebreak
\begin{prop}
The following assertions are equivalent:
\begin{enumerate}
\item 
The Hopf symbol is a Dirichlet character.

\item
$\jac{q}{A} \in \Q$ for all $q \in \Z$.
\end{enumerate}
In this case, we even have that $\jac{q}{A} \in \{0,1,-1\}$ for all 
$q \in \Z$.
\end{prop}
\begin{pf}
That the Hopf symbol is a Dirichlet character means by definition\endnote{\cite{Wash}, Chap.~3, p.~19.} that we have
$$\jac{qq'}{A} = \jac{q}{A} \jac{q'}{A}$$
for all $q,q' \in \Z$. Note that this equation is satisfied automatically
if~$q$ or~$q'$ are not relatively prime to~$N$. Comparing this equation to the 1-cocycle equation above, we see that it is equivalent to the condition
$\sigma_q(\jac{q'}{A}) = \jac{q'}{A}$. But that the Hopf symbol is invariant under the Galois group just means that it is a rational number.

If it now happens that the Hopf symbol is a Dirichlet character, then its image
$$\{\jac{q}{A} \mid q \in \Z, \gcd(q,N)=1\}$$
is a finite subgroup of the multiplicative group~$\Q^\times$.
As $\{1,-1\}$ is the largest finite subgroup of~$\Q^\times$, it must  contain all the Hopf symbols.\endnote{\cite{Lan}, Chap.~II.III.2, Def.~25, p.~114.}
\qed
\end{pf}

\subsection[Properties of the Hopf symbol]{} \label{PropHopfSymb}
To find out more about the Hopf symbol, the first step is to note that Corollary~\ref{GalMod} still holds in this more general situation:
\begin{lemma}
For all $\sigma \in \Gal(\Q_N/\Q)$, we have
$\sigma^2(u_i) = u_{\sigma.i}$. 
\end{lemma}
\begin{pf}
The expansion of the Drinfel'd element
considered in Paragraph~\ref{VerlMat} takes in the case of~$D(A)$ the form
$$u_D = \sum_{i,j=1}^k u_{ij} e_{ij}$$
In terms of these components, the equation $\Psi(u_D) = u \o u^{-1}$ obtained in Lemma~\ref{DoubleQuasitri} takes the form
$$\sum_{i,j=1}^k u_{ij} e_{i} \o e_j = 
\sum_{i,j=1}^k \frac{u_{i}}{u_{j}} e_{i} \o e_j$$
so that $u_{ij} = \frac{u_{i}}{u_{j}}$.
From Corollary~\ref{GalMod}, we get that 
$\sigma^2(u_{ij}) = u_{\sigma.i,\sigma.j}$,
which translates into
$$\frac{\sigma^2(u_{i})}{\sigma^2(u_{j})} = \frac{u_{\sigma.i}}{u_{\sigma.j}}$$
so that 
$\sigma^2(u_{i})/u_{\sigma.i} = \sigma^2(u_{j})/u_{\sigma.j}$.
As we have $u_1 = \ea(u) = 1$ and $\sigma.1=1$ by Lemma~\ref{GaloisChar},
we can insert~$j=1$ into this equation to get 
$\sigma^2(u_{i})/u_{\sigma.i} = 1$, from which the assertion follows immediately. 
\qed
\end{pf}

It may be noted that inserting~$\gamma$ for~$\sigma$ and using Proposition~\ref{GaloisChar}, we recover the fact that $u_{i^*} = u_i$, which expresses that~$\sa(u)=u$, a fact already pointed out in Paragraph~\ref{CharRing}.

As a consequence, we can derive several facts about the Hopf symbol:
\begin{prop}
The Hopf symbol is a root of unity. If~$N$ is odd, its order divides~$6$ and $2N$. If~$N$ is even, its order divides~$24$ and~$N$. Furthermore, we have
$$\jac{q}{A} = 1$$
if~$q$ is a square modulo~$N$.
\end{prop}
\begin{pf}
\begin{list}{(\arabic{num})}{\usecounter{num} \leftmargin0cm \itemindent5pt}
\item
From Proposition~\ref{MatIdent}, we get 
$\V \T \V \T \V \T = \chi_R(\ua^{-1}) \dim(A) \; \C^2 $.
As $\C^2$ is the unit matrix, this implies by taking determinants that\endnote{\cite{BakKir}, Thm.~3.1.19, p.~57f.}
$$\det(\V)^3 \det(\T)^3  = \chi_R(\ua^{-1})^k \; \dim(A)^k  $$
Proposition~\ref{MatIdent} also yields that
$\det(\V)^2 = \dim(A)^k \det(\C) = \pm \dim(A)^k$. 
Therefore, if~$q$ is relatively prime to~$N$, we have 
$\sigma_q(\det(\V)) = \pm \det(\V)$. If we now apply~$\sigma_q$ to the above equation, we get
$$\pm \det(\V)^3 \sigma_q(\det(\T))^{3}  = \chi_R(\ua^{-q})^k \; \dim(A)^k  $$
If we divide the two equations by each other, we therefore get that
$$\jac{q}{A}^k = \frac{\chi_R(\ua^{-q})^{k}}{\chi_R(\ua^{-1})^k} 
= \pm \frac{\sigma_q(\det(\T))^{3}}{\det(\T)^{3}}$$
Since $\det(\T)$ is an $N$-th root of unity, we see that~$\jac{q}{A}$ is a root of unity. Moreover, it is clear from its definition that
$\jac{q}{A} \in \Q_N$, which implies that $\jac{q}{A}^{2N} =1$, and actually 
$\jac{q}{A}^{N} =1$ if~$N$ is even.\endnote{\cite{Wash}, Exerc.~2.3, p.~17.}

\item
It now follows from the preceding lemma and Lemma~\ref{GaloisChar} that we have
\begin{align*}
\chi_R(u^{-q^2 l}) = \sigma_q^2(\chi_R(u^{-l}))
= \sum_{i=1}^k n_i \sigma_q^2(\frac{1}{u_i^l})
= \sum_{i=1}^k n_{\sigma_q.i} \frac{1}{u_{\sigma_q.i}^l}
= \chi_R(u^{-l})
\end{align*}
Dividing this equation by~$\chi_R(u^{-1})$, this shows that
$\jac{lq^2}{A} = \jac{l}{A}$, which shows for $l=1$ that $\jac{q^2}{A} =1$.
But the computation also shows that
\begin{align*}
\sigma_q^2(\jac{l}{A}) = \frac{\sigma_q^2(\chi_R(u^{-l}))}{\sigma_q^2(\chi_R(u^{-1}))}
= \frac{\chi_R(u^{-l})}{\chi_R(u^{-1})} = \jac{l}{A}
\end{align*}
Now if $\zeta$ is a primitive $N$-th root of unity, we have already shown in the first step that we can write $\jac{l}{A} = \pm \zeta^m$ for some~$m$, so that the preceding equation becomes
$\pm \zeta^{m q^2} = \pm \zeta^m$, which implies 
$m q^2 \equiv m \pmod{N}$ for all~$q$ that are relatively prime to~$N$, which means that~$N$ divides $m(q^2-1)$. If $N$ is odd, we see by taking~$q=2$ that~$N$ divides~$3m$, so that
$\jac{l}{A}^3 = \pm \zeta^{3m} = \pm 1$ and therefore~$\jac{l}{A}^6 = 1$.

If~$N$ is even, we have seen that we actually have 
$\jac{l}{A} = \zeta^m$ for some~$m$. If $N=\prod_i p_i^{m_i}$ is the prime factorization of~$N$ into powers of distinct primes, we can find by the Chinese remainder theorem a unit~$q$ modulo~$N$ that satisfies $q \equiv 3 \pmod{p_i^{m_i}}$ if~$p_i=2$ and  
$q \equiv 2 \pmod{p_i^{m_i}}$ if~$p_i \neq 2$. If~$p_i=2$, we therefore get that~$p_i^{m_i}$ divides~$8m$, and if~$p_i \neq 2$, $p_i^{m_i}$ divides~$3m$. In any case, $p_i^{m_i}$ divides~$24m$, so that~$N$
divides~$24m$, showing\endnote{\cite{CostGann}, \S~2.4, Prop.~3.b, p.~9.}
 that $\jac{l}{A}^{24} = 1$.
\qed
\end{list}
\end{pf}

It may be noted that this proposition asserts in particular that
$\chi_R(u) = \chi_R(u^{-1})$ if~$-1$ is a square modulo~$N$. If this happens, we can say even more about the Hopf symbol:
\begin{corollary}
Suppose that $\chi_R(u) = \chi_R(u^{-1})$. Then the Hopf symbol is a Dirichlet character, and we have $\jac{q}{A} \in \{0,1,-1\}$ for all 
$q \in \Z$.
\end{corollary}
\begin{pf}
Let $\gamma \in \Gal(\Q_N/\Q)$ be the restriction of complex conjugation considered in Paragraph~\ref{GaloisChar}. Using it, we can write the assumption in the form $\gamma(\chi_R(u^{-1})) = \chi_R(u^{-1})$.
If~$q$ is relatively prime to~$N$, we then also have
$$\gamma(\chi_R(u^{-q})) = \gamma(\sigma_q(\chi_R(u^{-1}))) =
\sigma_q(\gamma(\chi_R(u^{-1}))) = \sigma_q(\chi_R(u^{-1})) =\chi_R(u^{-q})$$
because~$\Gal(\Q_N/\Q)$ is abelian. Dividing this equation by the preceding one, we obtain $\gamma(\jac{q}{A}) = \jac{q}{A}$.
But since~$\jac{q}{A}$ is a root of unity by the preceding result, we also have $\gamma(\jac{q}{A}) = 1/\jac{q}{A}$. Therefore $\jac{q}{A}^2=1$
and $\jac{q}{A} \in \{1,-1\}$. The remaining assertions follow from Proposition~\ref{HopfSymb}.
\qed
\end{pf}

We have already pointed out that the condition $\chi_R(u) = \chi_R(u^{-1})$
means for the Hopf symbol that~$\jac{-1}{A} = 1$. A Dirichlet character with this property is called even.\endnote{\cite{Wash}, Chap.~3, p.~19.}

\subsection[Hopf symbols and the congruence subgroup theorem]{} \label{HopfSymbCong}
We have seen in Paragraph~\ref{ProjCong} that the kernel of the projective representation of the modular group is a congruence subgroup of level~$N$.
However, in the case where $\chi_R(u) = \chi_R(u^{-1})$, we have also seen that this projective representation comes from a linear representation.
This raises the question whether in this case also the kernel of the linear representation is a congruence subgroup of level~$N$. As we will see now, this in fact holds. Our reasoning is based on the following lemma, which is an adaption of an argument by A.~Coste and T.~Gannon to our situation:\endnote{\cite{CostGann}, \S~2.3, Thm.~2, p.~7.}
\begin{lemma}
Suppose that~$q, q' \in \Z$ satisfy $qq' \equiv 1 \pmod{N}$. 
Then we have
\begin{align*}
(\v \circ \t^{q'} \circ \v^{-1} \circ \t^{q} \circ \v \circ \t^{q'})(e_{m})
= \kappa \chi_R(\ua^{-q})  \; e_{\sigma_q^{-1}.m}
\end{align*}
\end{lemma}
\begin{pf}
It follows from Proposition~\ref{MatIdent} that
$$\V \T \V \T \V \T = \chi_R(\ua^{-1}) \dim(A) \; \C^2$$
and~$\C^2$ is the unit matrix by construction. If we write this out
in components, it means that
$$\sum_{j,l=1}^k \frac{s_{ij}}{u_j} \frac{s_{jl}}{u_l} \frac{s_{lm}}{u_m}
= \chi_R(\ua^{-1}) \dim(A) \delta_{im}$$
If we apply~$\sigma_q$ to this equation, it becomes
$$\sum_{j,l=1}^k \frac{s_{i,\sigma_q.j}}{u_j^q} \frac{s_{\sigma_q.j,l}}{u_l^q} \frac{s_{l,\sigma_q.m}}{u_m^q}
= \chi_R(\ua^{-q}) \dim(A) \delta_{im}$$
If we replace~$j$ by~$\sigma_q^{-1}.j$ and~$m$ by~$\sigma_q^{-1}.m$,
this becomes
$$\sum_{j,l=1}^k \frac{s_{ij}}{u_{\sigma_q^{-1}.j}^q} \frac{s_{jl}}{u_l^q} \frac{s_{lm}}{u_{\sigma_q^{-1}.m}^q}
= \chi_R(\ua^{-q}) \dim(A) \delta_{i,\sigma_q^{-1}.m}$$
But this can by the preceding proposition be written in the form
$$\sum_{j,l=1}^k \frac{s_{ij}}{u_{j}^{q'}} \frac{s_{jl}}{u_l^q} \frac{s_{lm}}{u_{m}^{q'}}
= \chi_R(\ua^{-q}) \dim(A) \delta_{i,\sigma_q^{-1}.m}$$
Multiplying~$e_i/n_i$ by this scalar and summing over~$i$, this becomes 
$$\frac{1}{\kappa} \sum_{j,l=1}^k  
\frac{s_{jl}}{u_l^q} \frac{s_{lm}}{u_{m}^{q'}} 
(\v \circ \t^{q'})(\frac{e_{j^*}}{n_j}) = 
\sum_{i,j,l=1}^k \frac{s_{ij}}{u_{j}^{q'}} \frac{s_{jl}}{u_l^q} \frac{s_{lm}}{u_{m}^{q'}} \frac{e_i}{n_i} = 
\chi_R(\ua^{-q}) \dim(A) \frac{1}{n_m} e_{\sigma_q^{-1}.m} $$
by Corollary~\ref{VerlMat}, and by repeating this argument we get
\begin{align*}
\chi_R(\ua^{-q}) \frac{\dim(A)}{n_m} e_{\sigma_q^{-1}.m} &=
\frac{1}{\kappa^2} \sum_{l=1}^k \frac{s_{lm}}{u_{m}^{q'}} 
(\v \circ \t^{q'} \circ \v \circ \t^{q})(\frac{e_{l}}{n_l}) \\
&= \frac{1}{\kappa^3}  
(\v \circ \t^{q'} \circ \v \circ \t^{q} \circ 
\v \circ \t^{q'})(\frac{e_{m^*}}{n_m})
\end{align*}
which gives  
\begin{align*}
(\v \circ \t^{q'} \circ \v \circ \t^{q} \circ \v \circ \t^{q'})(e_{m})
= \kappa^3 \chi_R(\ua^{-q}) \dim(A) \; e_{\sigma_q^{-1}.m^*}
\end{align*}
after substituting~$m^*$ for~$m$. Applying the antipode, which commutes with~$\v$ and~$\t$, we get
\begin{align*}
(\v \circ \t^{q'} \circ \sa\v \circ \t^{q} \circ \v \circ \t^{q'})(e_{m})
= \kappa \chi_R(\ua^{-q}) (\rha \o \rha)(R'R) \; e_{\sigma_q^{-1}.m}
\end{align*}
where we have used that 
$\kappa^2 \dim(A) = \rha(u) \rha(u^{-1}) = (\rha \o \rha)(R'R)$, as observed in Paragraph~\ref{MatIdent}. Now the assertion follows from Corollary~\ref{InvSigma}.
\qed
\end{pf}

From this lemma, we can now deduce the result indicated above:
\begin{thm}
Suppose that $\chi_R(u) = \chi_R(u^{-1})$. Then the kernel of the representation of the modular group on the center of~$A$ is a congruence subgroup of level~$N$.
\end{thm}
\begin{pf}
\begin{list}{(\arabic{num})}{\usecounter{num} \leftmargin0cm \itemindent5pt}
\item
To begin, recall that we saw in Paragraph~\ref{RibEl} that the condition
$\chi_R(u) = \chi_R(u^{-1})$, which for the Hopf symbol means that
$\jac{-1}{A} = 1$, ensures that the representation of the modular group is linear, and not only projective, so that it is meaningful to talk about the kernel. Recall also our convention from Paragraph~\ref{DualRepMod} that
$\kappa = \frac{1}{\chi_R(u)}$ in this case. Furthermore, we have just seen
in Corollary~\ref{PropHopfSymb} that the Hopf symbol is then a Dirichlet character and takes only the values~$0$, $1$, and~$-1$.

We have to verify the relations listed in Proposition~\ref{PresFac}.
The relations \mbox{$\bv^4 = 1$}, $(\bt\bv)^3 = \bv^2$, and $\bt^N = 1$
that are listed there first are clearly satisfied. 
Next, we verify the second relation, i.e., the relation
$\bt^{2^e} (\bv \bt^m \bv^{-1}) = (\bv \bt^m \bv^{-1}) \bt^{2^e}$,
where we have factored the exponent in the form $N=2^e m$ for~$m$ odd.
Now, as 
$\gt^{2^e} (\gv \gt^m \gv^{-1}) \gt^{-2^e} (\gv \gt^{-m} \gv^{-1})
\in \Gamma(N)$, we know from Theorem~\ref{ProjCong} that there is a scalar~$\mu \in K$ such that
$$\t^{2^e} (\v \t^m \v^{-1}) \t^{-2^e} (\v \t^{-m} \v^{-1})(z)
= \mu  z$$
for all $z \in Z(A)$. Inserting 
$z = (\v \t^m \v^{-1} \t^{2^e})(e_1)$, this becomes
\begin{align*}
\mu (\v \t^m \v^{-1} \t^{2^e})(e_1) = 
\t^{2^e} (\v \t^m \v^{-1}) (e_1)
\end{align*}
As $e_1$ is an integral, and we have
$\v(\A) = \v(z_1) = \kappa \dim(A) e_1$ by Corollary~\ref{VerlMat},
this equation can be rewritten as
\begin{align*}
\mu (\v \t^m \v^{-1})(e_1) = 
\frac{1}{\kappa \dim(A)} \t^{2^e} (\v \t^m) (\A)
\end{align*}
which implies
\begin{align*}
\mu \v(u^{-m}) = 
\mu (\v \t^m)(\A) = \t^{2^e} (\v \t^m) (\A)= \t^{2^e} (\v (u^{-m}))
\end{align*}
Because $\v(u^{-m}) \neq 0$, we see that~$\mu$ is an eigenvalue
for~$\t^{2^e}$, and therefore an $m$-th root of unity.

Similarly, we can insert $z = (\v \t^m \v^{-1} \t^{2^e})(\A)$ into the above equation, which then becomes
\begin{align*}
\mu (\v \t^m \v^{-1} \t^{2^e})(\A) &= 
\t^{2^e} (\v \t^m \v^{-1}) (\A) \\
&= 
\frac{1}{\kappa \dim(A)} \t^{2^e} (\v \t^m) (e_1)
= \frac{1}{\kappa \dim(A)} \t^{2^e} \v(e_1)
= \t^{2^e} (\A)
\end{align*}
Applying~$\v^{-1}$ to both sides and dividing by~$\mu$, this becomes
\begin{align*}
\t^m (\v^{-1}(u^{-2^e})) &= \frac{1}{\mu} \v^{-1} (u^{-2^e})
\end{align*}
Therefore~$1/\mu$ is an eigenvalue for~$\t^{m}$, and therefore a $2^e$-th root of unity. But now~$\mu$ is both a $2^e$-th root of unity and an $m$-th root of unity, which can only be if $\mu=1$, which in turn establishes our relation.

\item
The remaining two relations involve the diagonal matrices~$\gd(q)$.
Now note that with the help of these matrices the formula in the preceding lemma can be expressed as
\begin{align*}
\gd(q).e_{m} = \kappa \chi_R(\ua^{-q})  \; e_{\sigma_q^{-1}.m}
= \frac{\chi_R(\ua^{-q})}{\chi_R(u^{-1})}  \; e_{\sigma_q^{-1}.m}
= \jac{q}{A}  \; e_{\sigma_q^{-1}.m}
\end{align*}
which shows that $\gd(q).z = \jac{q}{A} \sigma_q.z$
for all $z \in Z(A)$. This means that the relation
$\gd(q) \bv  = \bv \gd(q)^{-1}$ listed third in Proposition~\ref{PresFac} reads
$$\jac{q}{A} \sigma_q.\v(z) = 
\frac{1}{\jac{q}{A}} \v(\sigma_q^{-1}.z)$$
But as $\jac{q}{A} = \pm 1$, this amounts to the relation
$\sigma_q.\v(z)  = \v(\sigma_q^{-1}.z)$, which was proved in Proposition~\ref{ActCent}.

\item
Finally, for the fourth relation in Proposition~\ref{PresFac}, we have on the one hand that
\begin{align*}
\gd(q) \bt.e_j = \frac{1}{u_j} \gd(q).e_j 
= \frac{1}{u_j} \jac{q}{A}\; e_{\sigma_q^{-1}.j}
\end{align*}
and on the other hand
\begin{align*}
\bt^{q^2} \gd(q).e_j = 
\jac{q}{A}  \; \bt^{q^2}.e_{\sigma_q^{-1}.j}
= 
\jac{q}{A} \frac{1}{u^{q^2}_{\sigma_q^{-1}.j}}  \; e_{\sigma_q^{-1}.j}
\end{align*}
Both expressions are equal by Lemma~\ref{PropHopfSymb}, so that all the required relations are satisfied. This proves that~$\Gamma(N)$ is contained in the kernel, and that the level of the kernel is exactly~$N$ follows as in Theorem~\ref{CongDrinfDoubl}.
\qed
\end{list}
\end{pf}

The preceding theorem generalizes Theorem~\ref{CongDrinfDoubl}, because in the case where $A=D(H)$ is the Drinfel'd double of a semisimple Hopf algebra~$H$, we saw in Paragraph~\ref{RoleEval} that
$$\chi_R(u_D) = \chi_R(u_D^{-1}) = \dim(H)$$
Applying~$\sigma_q$ to this formula, we see that also 
$\chi_R(u_D^{-q}) = \dim(H)$, so that $\jac{q}{A} = 1$.
We therefore see that the formula
$$\gd(q).z = \jac{q}{A} \sigma_q.z$$
that we obtained in the preceding proof reduces to Theorem~\ref{GalMod}.
However, one has to keep in mind that all these results were used in the preceding proof.

\newpage

\addcontentsline{toc}{section}{Notes} \label{Notes} 
\theendnotes

\newpage

\addcontentsline{toc}{section}{Bibliography} \label{Bibliography} 

\end{document}